# Congestion-Aware Charging Coordination for Electric Ride-Hailing Fleets under Stochastic Demand

Tai-Yu Ma, Luxembourg Institute of Socio-Economic Research (LISER), 11 Porte des Sciences, 4366 Esch-sur-Alzette, Luxembourg
Richard D. Connors, Department of Engineering, University of Luxembourg, Esch-sur-Alzette, Luxembourg
Francesco Viti, Department of Engineering, University of Luxembourg, Esch-sur-Alzette, Luxembourg

**Abstract**

In this study, we develop a dynamic charging scheduling method that integrates day-ahead charging planning and online reactive vehicle dispatch and charger assignment. The day-ahead planning model minimizes expected charging costs by considering anticipated waiting times and overnight energy replenishment costs. The online model then adapts this plan by maximizing the profit from vehicle dispatching and assigning vehicles to chargers through a reservation-based charging request scheduler to minimize the total charging-related operational time. The developed charging policy is tested on a Manhattan-like study area using synthetic data drawn from NYC yellow taxi data with a fleet size of 100 vehicles given the scenarios of 3000 and 4000 customers/day. The computational results show that our method outperforms different benchmark policies with up to +21.68% profit and +22.22% service rate for 4000 customers relative to the weakest benchmark. Sensitivity analyses are conducted for different system parameters, and managerial insights are discussed.

## 1. Introduction

The climate change crisis has motivated governments and Transport Network Companies to accelerate fleet electrification to reduce $CO_2$ emissions. Charging management is a major obstacle to this transition because electric vehicles (EVs) take far longer to charge than internal combustion engine vehicles take to refuel. Paffumi et al. (2015) compared EV characteristics and alternative charging strategies using large-scale EV driving pattern data. They demonstrated the importance of coordinating EV charging decisions in order to meet daily mobility requirements and energy demand. According to a recent TLC Electrification Report (Taxi & Limousine Commission, 2022), TLC plans to transform its licensed fleet to EVs by 2030 to reduce environmental impact (around 600k tons of $CO_2$ in 2022). However, due to the high investment costs of rapid chargers, only a limited number of fast charging stations are available in major cities in the United States. Furthermore, the high daily mileage of TLC's vehicles necessitates drivers to charge their vehicles several times a day, relying mainly on rapid chargers to save charging time (Jenn, 2019). Making good use of scarce fast-charging capacity under uncertain demand is a central challenge for this transition.

Several factors make dynamic charging management of electric ride-hailing systems challenging. First, customer demand is volatile, affecting vehicles' charging needs over time. Second, rapid chargers (charging power ≥50kW) are limited due to their high investment costs. This might result in EVs queuing at charging stations, thereby increasing charging station search costs for the drivers. Furthermore, charging costs might vary according to the time of day due to variations in electricity prices. The decision on when and how much energy to recharge becomes a significant online decision problem for drivers/fleet operators to maximize their profit. However, most studies assume constant energy prices and neglect congestion issues at charging stations (Al-Kanj et al., 2020; Shi et al., 2020).

Methodologies developed in recent years are mainly based on mixed integer programming optimization approaches, where the problem is decomposed into sequential sub-problems over different planning horizons for vehicle dispatching, relocation, and charging scheduling to respond to uncertain customer demand. In the first stage, a long-time horizon planning model is used to determine how many vehicles to charge and at what type of charger for each decision time interval (typically half an hour) for the day by anticipating vehicles' energy demand. The objective is to determine vehicles' charging times and durations for the day to minimize a desired objective function (e.g., minimize the number of total charging trips and rejected customers (Jamshidi et al., 2021) or the shortage of available vehicles (Zalesak and Samaranayake, 2021) or total charging operational costs (Ma, 2021) for the planning horizon). In the second stage, a short-time horizon planning model is applied to determine where to charge to minimize total charging access costs or charging operational costs (including access times, waiting times, and charging durations) under the capacity constraints of charging stations. Three drawbacks can be identified in these existing works. First, when vehicles wait to charge at charging stations, there are no maximum waiting time limits, resulting in unrealistically long waiting times at charging stations. Second, there is no minimum charging duration constraint per each charging operation, which might result in undesirably short charging durations (damaging battery lifespan) and inefficient charging operations (requiring higher charging access times/costs and setup times due to more frequent charging operations). Third, vehicles only need enough energy to reach the end of their working day, when overnight charging is available; ignoring this results in more charging than is needed. Consequently, vehicles might charge longer, resulting in a shortage of available vehicles to serve customers and suboptimal fleet availability. The above identified drawbacks reflect a common structural issue in the existing studies. Addressing these issues is therefore essential to ensure operational realism and to avoid systematically biased performance evaluation.

We address these drawbacks with a more realistic model of charging operations for ride-hailing fleets under stochastic demand. Our main contributions are as follows.

a) We propose a novel two-stage sequential optimization approach that first integrates an approximate day-ahead charging plan and subsequently adapts it with online vehicle dispatch and a reservation-based vehicle-to-charger assignment under stochastic demand. The objective is to maximize the ride-hailing system's profit over the entire planning horizon. In existing studies (e.g., Ma, 2021; Jamshidi et al., 2021; Zalesak and Samaranayake, 2021), charger capacity is represented as an aggregate per-epoch limit on the number of vehicles charging, whereas our online reservation scheduler resolves capacity at the level of individual chargers accounting for charger-level congestion. In contrast to our earlier two-stage model (Ma, 2021), which ignored charging station capacity constraints when deriving day-ahead charging plans for individual vehicles, the proposed approach explicitly incorporates charging station congestion. Furthermore, our day-ahead charging plan minimizes total charging costs over a 24-hour period by accounting for time-varying electricity prices and additional overnight charging costs required to fully recharge the vehicles. These characteristics represent the key methodological differences from our earlier study (Ma, 2021).
b) We develop an online reactive charging decision model that adapts to the current system state and anticipates vehicles' waiting times at charging stations. The remaining energy needs for the day are anticipated, and vehicle charging time, duration, and location are optimized online. Minimum charging duration and maximum waiting time limits are incorporated as configurable operational policy assumptions to discourage long waits and short sessions. Our computational study shows that the proposed charging policy significantly increases the customer service rate and operators' profit by reducing vehicles' waiting times and charging durations at charging stations.
c) This study considers more realistic charging operation modeling and simulation under congested charging facilities. Moreover, vehicle queuing and waiting for charging are explicitly simulated, and more flexible energy prices are considered (time-varying and charging-type-specific), allowing testing the performance of the proposed methodology under different dimensions of uncertainty (demand, charging waiting time, queuing, and energy prices).

d) Numerical experiments are conducted on a Manhattan-like study area with synthetic demand data drawn from NYC yellow taxi data, and the impacts of different system parameters (battery capacity, maximum waiting times at chargers, number of chargers, and customer demand variability) are analyzed. The performance of the proposed approach is compared with several benchmark charging policies, showing its benefit in increasing the customer service rate and total profit of the operator in a stochastic environment.

The organization of this paper is as follows. We first present the related studies. Section 3 presents the problem description, developed models, and benchmark charging policies. Section 4 describes the test instances and reports the computational studies. The impact of different system parameters is analyzed. Finally, the conclusion is drawn, and future extensions are discussed.

## 2. Related studies

Dynamic charging management of ride-hailing systems involves sequential decision making under uncertainty. Recent works for addressing this problem are mainly based on optimization-based approaches and reinforcement learning (RL). Earlier studies can be classified into two categories of methodology: sequential mixed integer programming approaches and RL. As it is not possible to jointly optimize all the decision variables (vehicle dispatching and charging) under uncertainty and system dynamics, the problem can be approximated as a sequential decision-making problem based on a decomposition approach. The sequential mixed integer programming approach is based on this optimization framework. This approach decomposes the dynamic charging scheduling problem into multiple decision horizons, where the long-horizon planning for the day determines the number of vehicles to charge, charging time, and durations for each charging decision time interval (e.g., half an hour). In contrast, short-horizon planning aims to determine where to charge, given the charging plans obtained at the first step. Similar two-stage optimization frameworks have been developed for the integrated planning and operation of EV-sharing systems, which can be used to determine the fleet size and plan the charging infrastructure, as well as evaluate operational performance under stochastic demand (Zhuang et al., 2025). The second approach (RL) considers dynamic vehicle dispatching, relocation, and charging management as sequential decision-making problems under uncertainty and models the problem as a Markov decision process. Under the RL framework, the operator is considered as an agent making these decisions to maximize the total profit of the system. Recent studies extend the single-agent-based approach to the multi-agent-based approach to allow individual vehicles to make their own decisions in response to local information of the system. More detailed reviews are described as follows.

**Optimization-based approaches**: Earlier works on dynamic ride-hailing systems focus on vehicle dispatching and relocation optimization under stochastic demand. Vehicle charging operations are simplified by assuming uncapacitated charging stations or simply neglected by considering internal combustion engine vehicles. Existing methodology can be classified into four categories as follows.

<u>Model predictive control</u>
This approach uses a model to represent the system dynamics based on which a set of control variables is optimized to minimize a cost function over a prediction horizon. The process is repeated until the end of the planning horizon. Zhang et al. (2016) proposed a model predictive control approach for operating a dynamic ride-hailing system using autonomous vehicles. Vehicles' charging scheduling problems are not considered. Iacobucci et al. (2019) extended their work for the fleet management of shared autonomous EVs based on the model predictive control approach. Two different model predictive control schemes interact over different planning horizons: one for vehicle dispatching and relocation decisions, and another for vehicle charging scheduling. However, the applicability of the developed approach is limited to small problem sizes with tens of vehicles, and charging congestion is not considered by assuming uncapacitated charging stations.

Sequential mixed integer programming

This approach decomposes the system dynamics into several sub-processes with different temporal dimensions, usually embedded in a hierarchical structure. Given this underlying structure, different models are developed to pre-allocate limited resources or re-optimize the decision variables based on the current system state or prediction of system changes over a prediction horizon. For example, Jamshidi et al. (2021) proposed a three-stage sequential MILP model to address e-taxi dispatching, relocation, and charging with charging station capacity constraints and time-of-use energy prices. However, waiting times at charging stations are approximated without explicitly modelling charging queuing times on different chargers. Zalesak and Samaranayake (2021) developed a MILP model to optimize the charging schedules of ride-pooling systems using EVs. A two-stage planning framework is proposed: a long-term planning horizon optimization model for determining when to charge and a short-term planning horizon optimization model to determine where to charge. The objective of the long-time horizon planning model is to maintain a sufficient fleet size to meet varying customers' demands and vehicles' energy needs for the planning horizon. The number of vehicles recharged for each decision time interval cannot exceed the capacity of the total charging stations. Heterogeneous charging infrastructure and time-dependent energy prices are not considered. For the short-time horizon planning, a simple vehicle-to-charging-station assignment model is formulated to minimize vehicles' total charging access distance under charging station capacity constraints. Ma (2021) proposed a two-stage MILP model for dynamic charging management of shared ride-hailing services. The developed approach anticipates vehicles' driving needs and waiting time at chargers to determine charging time and durations per half hour to minimize total daily charging operational costs. An online vehicle-to-charger assignment model is applied to minimize total charging operational times under charging station capacity constraints. Several simplifications are made for this study. First, the charging speed of vehicles is assumed linear without distinguishing the fact that vehicles' charging speed slows down when their state-of-charge (SoC) is above a threshold (around 80% of the vehicle's battery capacity). Note that SoC is measured in kWh. Second, the day-ahead plan does not distinguish charger types or charging-station capacity, and no minimum charging duration is imposed. Yang et al. (2019) proposed a two-stage charging coordination approach for optimizing electric taxi fleet charging operations by considering the current queuing state at charging stations. The first-stage problem determines when vehicles go to recharge using a time-dependent charging cost function that considers average revenue loss per kWh charged. The second-stage problem determines where to recharge as a Nash equilibrium problem to model non-cooperative taxi drivers' charging station selections.

Optimization-embedded heuristics

This approach decomposes the sequential decision process under a rolling horizon framework, then solves an optimization model iteratively until the end of the planning period. The considered optimization problem aims to jointly optimize idle vehicle relocation and charging operations to minimize a cost function under vehicle operations and charging capacity constraints. Pantelidis et al. (2022) developed a MILP to optimize vehicle relocation and charging for electric carshare systems jointly. The optimization problem is modelled as a route-capacitated minimum-cost flow relocation problem. A node-charge graph is proposed to allow partial recharges over different discretized recharge levels to cover vehicle energy needs for serving customers. A simulation case study is conducted using realistic carshare data in Brooklyn. Yi and Smart (2021) proposed a heuristic for vehicle repositioning and charging management of ride-hailing systems using autonomous EVs. Vehicles go to the nearest unoccupied charging stations to recharge when their SoCs are below a predefined threshold. Dean et al. (2022) proposed a MILP model to jointly optimize vehicle relocation and charging operations for shared autonomous EVs using zones as the operational units. When assigning vehicles for charging in a zone, the number of vehicles for charging cannot exceed the total charging station capacity. Vehicles' queuing times at charging stations are not explicitly modelled.

Charging demand management using dynamic pricing

The above three approaches do not leverage dynamic pricing to mitigate charging congestion for better utilization of limited charging facilities. Several studies integrate dynamic charging pricing to manage the charging demand of demand-responsive mobility systems to better allocate charging resources and balance charging demand. This congestion pricing approach has been widely applied for road traffic management and charging management of private EVs (Yu and Hong, 2017; Zhang et al., 2021). For shared mobility services, Maljkovic et al. (2023) propose a hierarchical control, game-theoretical-based mechanism to coordinate the charging operations of multiple ride-hailing operators. The problem is formulated as a reverse Stackelberg game with a central authority (e.g., a grid operator) at the upper level (minimizing a cost function) and multiple ride-hailing operators at the lower level (maximizing their own profits). Assuming different charging tariffs on different charging stations, the ride-hailing operators determine their vehicle-to-charging-station assignment to minimize the charging operation costs and maximize the total profit. The central authority wants to balance the load over different charging stations using dynamic pricing. Given the non-cooperative behaviour of different ride-hailing operators, system optimal pricing policies under Nash equilibrium are searched. Another similar study is Laha et al. (2019), who developed dynamic pricing strategies for charging stations to maximize their own utility under limited charging station capacity, considering that EV users select their charging stations based on charging prices and charging station location via real-time information sharing. The problem is modelled as a multi-leader multi-follower Stackelberg game, and a dynamic pricing strategy is proposed to achieve the Stackelberg equilibrium.

Reinforcement-learning-based approaches

RL is a model-free approach that has been successfully applied to solving various sequential decision-making problems under uncertainty (Farazi et al., 2021). For dynamic ride-hailing systems using EVs, Al-Kanj et al. (2020) proposed an approximate dynamic programming approach for ride-hailing fleet management to maximize the operators' profit. Vehicles' dispatching, relocation, and charging decisions (actions) are controlled by a central controller. The system state (vehicles' locations, SoCs, activities, etc.) and time are discretized. Charging stations are assumed uncapacitated, and vehicles are fully recharged for each charging operation. As the possible state-action combination is huge, an approximate dynamic programming approach is applied with hierarchical aggregation for value function approximation. Yan et al. (2023) proposed an online model-based RL algorithm based on State–action–reward–state–action (SARSA) for joint vehicle dispatch and charging optimization of electric ride-hailing systems. Like the previous study, a full-recharge policy is applied. Kullman et al. (2021) developed a deep reinforcement learning (DRL) approach to overcome the curse of dimensionality for electric ride-hailing systems. A hybrid scheme is proposed where the problem of vehicle dispatching is solved by a MILP optimization model, and vehicle relocation and charging decisions are made under the RL framework. Unlike the single-agent-based RL approaches, Ahadi et al. (2023) propose a multiagent-based DRL approach for the fleet management of shared and autonomous EVs. The authors proposed a hierarchical learning and mean-field approximation approach to coordinate vehicles' charging decisions under charging station capacity constraints to maximize the total revenue of the fleet. A methodological comparison of the existing literature and the present study on charging scheduling for electric ride-hailing systems is provided in Appendix A.

A comprehensive review of RL-based approaches for EV charging management can be found in Abdullah et al. (2021). We review reinforcement-learning approaches because they represent the main alternative paradigm for this class of problem; however, they require extensive offline training on a fixed environment and are not directly comparable to the training-free, deployable policies evaluated here, so the benchmarks in this study are drawn from the optimization- and rule-based literature instead.

## 3. Dynamic charging management with charging congestion and time-varying energy prices

In this section, we first present the problem description and the discrete-event simulation framework, then formulate the problem using mixed integer linear programming, which includes three models: a) day-ahead charging planning, b) vehicle dispatching, and c) online vehicle-to-charger assignment. We present the simulation framework to test the proposed approach and compare its performance with five benchmark charging approaches.

***Notation***

| | |
|---|---|
| *Index* | |
| $t, t(i)$ | Time or the time of a decision point $i$ within the planning horizon |
| $h$ | Index of the day-ahead charging schedule planning decision epoch |
| $s$ | Charger |
| $r$ | Request (customer) |
| *Set* | |
| $H^{\mathscr{b}}$ | Set of decision epochs for online dispatch and charging operations with decision time interval $\Delta_t$ one minute, i.e., $H^{\mathscr{b}} = \{1,2,\dots,T-T_0\}$ |
| $\widehat{H}^{\ell}$ | Set of decision epochs for day-ahead charging schedule planning with decision time interval $\Delta_{\ell}$, $\widehat{H}^{\ell} = \{1,2,\dots,\left\lceil\frac{T-T_0}{\Delta_{\ell}}\right\rceil\}$ |
| $H^{\ell}$ | Set of the shifted decision epochs for day-ahead charging schedule planning, i.e., $H^{\ell} = \{1,2,\dots,n_{H^{\ell}}+1\}$ where the first index 1 in $H^{\ell}$ denotes the first charging decision epoch in $\widehat{H}^{\ell}$ with SoC lower than $E_v^{max}$ (80% of battery capacity as vehicles are not to charge above this threshold to save charging times during service hours). The end $n_{H^{\ell}}+1$ corresponds to the end of service hours. |
| $R_t$ | Set of unserved requests (customers) at time $t$ |
| $V$ | Set of vehicles |
| $\tilde{V}_h$ | Sets of vehicles that need to recharge at epoch $h \in H^{\ell}$ (output of day-ahead charging plan model) |
| $V_t$ | Set of idle vehicles at time $t$ (online optimization) |
| $\Omega_t$ | Set of need-to-charge vehicles at time $t$ (online optimization) |
| $\widetilde{\Omega}_t$ | Set of to-charge vehicles at time $t$ (online optimization) |
| $S$ | Set of chargers |
| *Parameter* | |
| $T_0$ | The start time of the ride-hailing service (i.e., 6:00 a.m.) |
| $T$ | The end time of ride-hailing service hours (i.e., 20:00) |
| $\overline{H}$ | Planning horizon, i.e., $\overline{H} = [T_0, T]$ |
| $n_{H^{\ell}}$ | The number of charging decision epochs in $\widehat{H}^{\ell}$ where the SoC of vehicles is below 80% |
| $u$ | Type of chargers, $u \in U = \{\text{fast}, \text{slow}\}$ |
| $\Delta_t$ | Vehicle dispatch (serving customers) time interval (1 minute) |
| $\Delta_{\ell}$ | Day-ahead charging schedule planning decision time interval (e.g., 30 minutes) |
| $p_h$ | Average energy price for charging decision epoch $h$ (dollar/kWh) |
| $\tilde{p}_s$ | Charger-type-specific charging cost (dollar/kWh) |
| $\bar{p}$ | Flat fare for overnight charging on level-2 public charging stations |
| $\phi(q)$ | Energy price function during service hours; q denotes a charging session |
| $C$ | Average travel cost for recharge per charging operation (dollar) |
| $\varphi_s$ | Charging power of charger $s$ (kW) |
| $\mu$ | Energy consumption rate per kilometer travelled (kWh/km) |
| $\gamma$ | Average profit per minute travelled (dollar/minute) |
| $\alpha_s$ | Minimum charging time per charging operation for the charger type for the charger s (minute) |
| $\beta_0, \beta_1$ | $\beta_0$ is the base rate, and $\beta_1$ is a distance-based operating fee (dollar/km) for a service trip |
| $B_v$ | Battery capacity of vehicle $v$ (kWh) |
| $E_v^{min}$ | Minimum (reserve) SoC of vehicle $v$ (kWh) |
| $E_v^{max}$ | Maximum threshold of SoC of vehicle $v$ with theoretical maximum charging speed (kWh) (e.g., 80% of vehicle's battery capacity) |
| $E_v^0$ | Initial SoC of vehicle $v$ for the charging schedule planning (KWh) |
| $\bar{E}_v, \bar{E}_v(t)$ | Planned battery level after recharge for vehicle $v$ at the end of charging decision epoch $h$ ($h$ is dropped) or at decision time $t$ |
| $Y_s^{max}$ | The maximum amount of energy that can be charged on charger $s$ during one charging decision epoch (kWh) |
| $\delta_h$ | Average energy consumption of a vehicle for charging decision epoch $h, \forall h \in H^{\ell}$ (kWh) |
| $\bar{t}_{vs}, \bar{t}_{vo_r}$ | Travel time from vehicle $v$'s current location to the location of charger $s$ or customer $r$'s pickup location(minute) |
| $\bar{d}_{ij}$ | Shortest path distance between location $i$ and $j$ (km) |
| $W_{rt}$ | Experienced wait time of request $r$ by the time $t$ (minute) |
| $W^{max}$ | Maximum waiting time of customers (same for each customer) (minute) |
| $\bar{W}_{vs}^{t(i)}$ | Waiting time that vehicle $v$ would experience if it departs immediately to go to charger $s$ at decision time $t(i)$(minute) |
| $\bar{W}_{hs}$ | Expected waiting time of vehicles when arriving at a charger $s$ at epoch $h, \forall \in H^{\ell}$ |
| $\pi$ | Cost per kilometer travelled of vehicles (dollar) |
| $\theta$ | SoC threshold under which vehicles are added to the pool of go-charge vehicles (% of vehicles' battery capacity) |
| $o_r, d_r$ | Pickup/drop-off location of request $r$ |
| $\tau^{max}$ | Maximum charging duration per charging session (minute) |
| *Decision/auxiliary variables* | |
| $x_{hs}^v$ | Indicator: 1 if vehicle $v$ is assigned to charger $s$ in epoch $h$, and 0 otherwise |
| $y_{hs}^v$ | Amount of charged energy for vehicle $v$ at charger $s$ during epoch $h$ |
| $\hat{x}_{vs}$ | Indicator: 1 if a vehicle is assigned to charger $s$, and 0 otherwise |
| $\hat{y}_{vs}$ | Amount of charged energy for vehicle $v$ at charger $s$ |
| $m_{rv}$ | Indicator: 1 if request $r$ is assigned to vehicle $v$, and 0 otherwise |
| $e_{vh}, e_{vt(i)}, e_v(t)$ | Battery level of vehicle $v$ at the beginning time of epoch $h \in H^{\ell}$ (charging schedule planning), at time $t(i)$, $i \in H^{\mathscr{b}}$ or at time $t$ (kWh) |

3.1. Problem description

Consider a ride-hailing service operated by a transport network company (TNC, referred to as operator hereafter) with a fleet of homogeneous EVs. Vehicles are equipped with dedicated communication devices for real-time communication to the operator's control center of vehicle states (e.g., location, vehicle's activity, e.g., idle, charging, or serving customers, SoCs of vehicles, etc.). Vehicles' dispatch and charging operations are controlled by the operator. Customers arrive randomly and send their ride requests via a smartphone app by indicating their pickup location and desired pickup time. The operator applies a batch assignment method for vehicle dispatching. Customers have limited patience, which is quantified by a maximum waiting time threshold (assume identical for each customer). The operator receives revenue from customer rides, and the service fee depends on the distance travelled. We assume customers will use the service if their waiting time is below a predefined maximum threshold. Ride-hailing drivers own the vehicles. Due to working hours regulations, ride-hailing drivers can drive around 10 hours (or 10 hours of passenger time) per day[1]. We divide the 24-hour period into two time periods: service hours from 6:00 a.m. to 8:00 p.m. and off-service hours from 8:00 p.m. to 6:00 a.m. the following day. Due to the limited number of fast public charging stations, we assume the operator uses a mixed strategy to charge their fleet. The drivers charge their vehicles on operator-owned stations during service hours and on abundant public level-2 chargers (typical charging power from 6.4 kW to 7.2 kW or 19.2 kW) overnight. This strategy balances charging costs and efficiency to avoid long waiting time at public fast charging stations during service hours[2]. We assume that the operator has a limited number of charging hubs/stations dedicated to the fleet charging during service hours. The operator aims to optimize vehicle charging operations under stochastic customer demand to maximize the fleet's total profit over a 24-hour period. The fleet size is assumed to be much larger than the number of operator-owned chargers. Using fast chargers incurs higher charging costs due to charger-type-specific fees. The assumptions for the initial state of charge (SoC) of the vehicles, energy consumption, and charging operations are detailed below.

– Vehicles are the same type and are fully charged at the start of each day's service. Vehicles must keep their State of Charge (SoC) above a predefined minimum level as possible. During service hours, charging can only occur at operator-owned charging stations. This assumption reflects emerging operational models in which transportation network companies contract dedicated charging infrastructure to ensure service reliability and mitigate congestion risks. Moreover, this allows explicit modeling of charging congestion dynamics under controlled charging infrastructure.
– During off-service hours, drivers fully recharge their vehicles at other public level-2 charging stations for overnight charging. We assume that there are enough public level-2 chargers[3] and that the cost of overnight charging is based on the standard pricing scheme for public level-2 charging stations. This assumption reflects operational practice and allows the focus on charging coordination during service hours.
– The operator has full access to the occupancy state of each charger, as well as the status of vehicles at any time. The waiting time for each charger and the length of the charging queue are available to the operator at any given moment. This allows the operator to control and coordinate vehicle charging operations over time to enhance the operational efficiency of the system. Although vehicles may be driver-owned, we assume that the operator coordinates dispatching and charging decisions during service hours, and drivers fully comply with the operator's orders for vehicle dispatch and charging during service hours. This full driver compliance is imposed as a behavioral simplification to evaluate the potential of coordinated fleet management and charging strategies. Recent research has highlighted the influence of heterogeneous driver responses on ride-hailing dispatch outcomes (Luo et al., 2026). Relaxing this assumption by accounting for heterogeneous driver compliance behavior can be a future extension of this study.
– The cost of charging depends on the time of day and the type of charger. The operator pays for electricity costs of recharging from exogenous time-of-use (ToU) energy prices during service hours. The ToU energy prices are deterministic and known day-ahead. The net profit is calculated as the customer service revenue minus vehicle travel costs and daytime and overnight charging costs.

---

[1] Fatigued Driving Prevention - Frequently Asked Questions. https://www.nyc.gov/site/tlc/about/fatigued-driving-prevention-frequently-asked-questions.page

[2] https://www.blackcarnews.com/article/gravity-installs-super-fast-ev-chargers-in-nyc?

[3] According to ChargeHub, NYC has 2049 public charging station ports (level 2 and level 3) within 15km (97% of the ports are level 2). https://chargehub.com/en/countries/united-states/new-york/new-york.html?utm_source=chatgpt.com

- Charging facilities are heterogeneous, and the number of operator-owned chargers is fixed. When chargers are occupied, vehicles must wait. There is no overlap allowed for each charger; only one vehicle at a time. This enables realistic charging congestion modeling during service hours.
- Travel speed and energy consumption per kilometer travelled are based on averaged values and are assumed to be identical and deterministic for all vehicles. Vehicles' energy consumption and charging waiting times are calculated based on their realized trips and vehicles' arrival times at chargers and charging station occupation state dynamics in the simulation. As stochastic events (e.g., weather, road congestion, etc.) may deviate travel speed and energy consumption from their expected values, the reported results should be interpreted as performance under expected operational conditions. In reality, these stochastic events may change vehicle arrival times and dispatch feasibility, which require further online adjustment to meet real-time operational constraints.
- Customer demand (i.e., origins, destinations, and request times) is stochastic and revealed dynamically during service hours. We model stochastic demand as a non-homogeneous stochastic process based on the empirical distribution of historical ride data in the study area. Customer requests are generated according to this empirical temporal distribution, capturing realistic demand fluctuations over the day.

As the problem is a sequential decision-making problem involving large state-action space and stochastic demand, we propose a **two-stage sequential optimization approach** to solve the problem and compare its performance with a list of benchmark approaches. In the **first stage (day-ahead planning)**, we solve a deterministic optimization problem to generate an approximate fleet charging plan using expected system performance derived from historical data. This stage acts as a *(charger) resource allocation coordination* that anticipates charging congestion and distributes charging demand over time. In the **second stage (online control)**, decisions (vehicle dispatch to serve customers and recharging) are made dynamically at a fine temporal resolution (e.g., every minute) based on the realized system state and stochastic customer arrivals. The day-ahead charging planning is a static model for the planning horizon (service hours, e.g., 06:00-20:00), but it integrates the overnight charging costs for the full recharge so the SoC of vehicles is 100% (or maximum capacity in practice) at the beginning of the next day. Operators leverage historical data (including energy consumption of vehicles, average waiting time at each charger, average earnings per planning time interval, and energy prices) to derive a day-ahead charging plan that minimizes total charging costs. The model devises static vehicle-specific charging schedules to guide vehicles' charging times and target SoCs after recharging for each charging time slot (e.g., 30 minutes) for the day (Model P1, described later).

The **online dispatch and charging management (stage 2)** is a dynamic model with a homogeneous decision time interval of $\Delta_t$ minutes (i.e., 1 minute) during service hours. This dynamic decision process is modeled based on the discrete-event simulation framework (detailed later) with two online optimization models for vehicle dispatch and vehicle charging management, i.e., the optimisation models are deterministic approximations with stochastic customer arrival realizations. During the day, the online planner/controller dispatches idle vehicles to serve customers or assigns vehicles to recharge to maximize total profit, sequentially, until the end of service hours. As shown in the system framework (Figure 1), customers arrive randomly during the day. The operator dispatches idle vehicles to serve customers based on a batch dispatching model (Model P2) to maximize the profit of vehicle-customer matching, given the customer's maximum waiting time and vehicles' SoC constraint at each decision time $t(i), \forall i \in \{1,2,\dots,T-T_0\}$. The day-ahead charging schedule serves as input for the online planner and is updated at each decision time $t(i)$ based on the current system state and anticipated energy needs to the end of service hours. An online vehicle-charger assignment model (Model P3) determines where to charge to minimize the total charging operational time per assignment epoch. At the end of the service hours, drivers go to some public level-2 charging stations to charge to full. The overnight charging costs of the fleet are then deducted to calculate the net profit. We neglect the locations and access travel costs of overnight public charging stations. This aspect can be easily integrated into the practical application of the model. Considering travel time or energy consumption uncertainty is a possible extension of this work in the future.

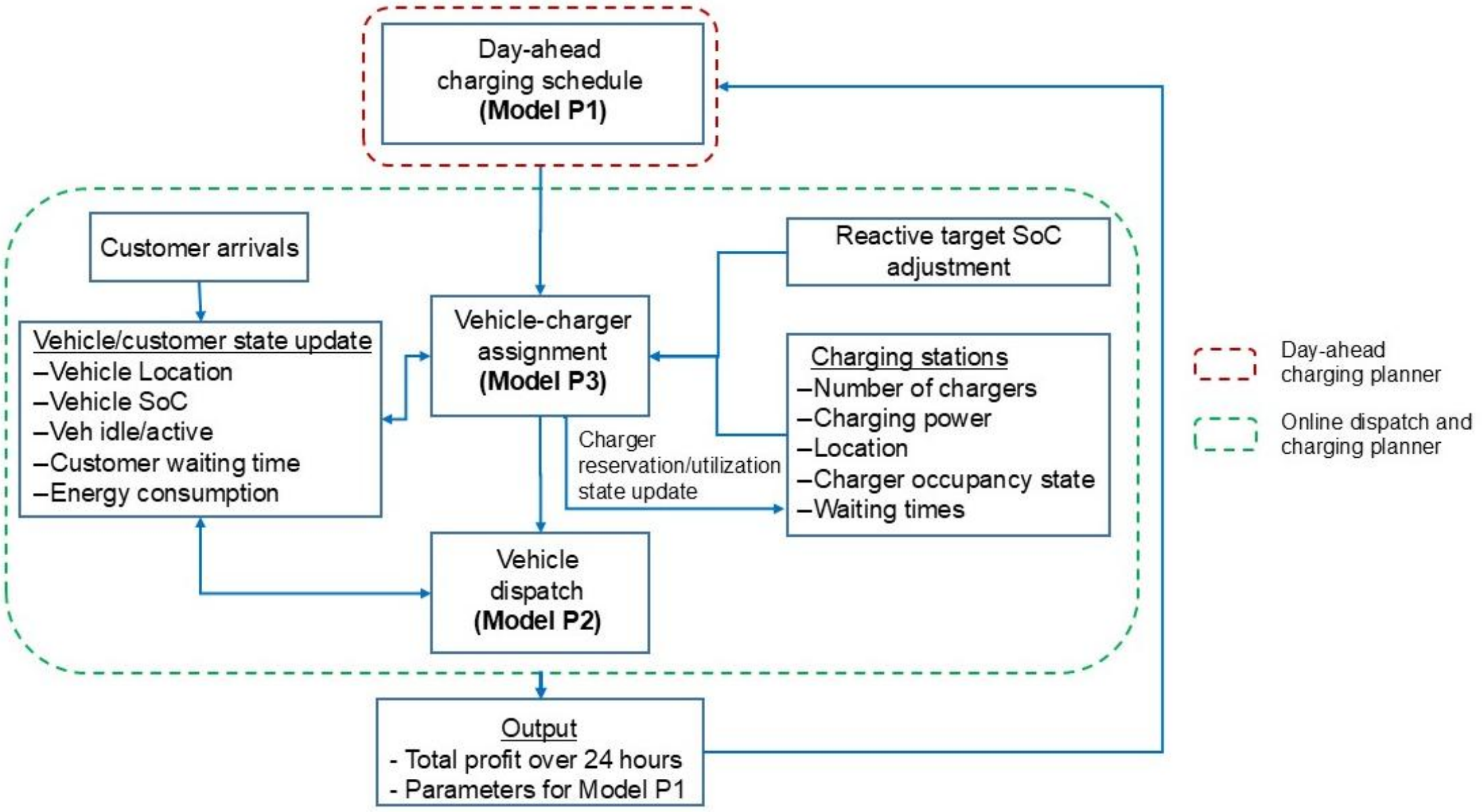

Figure 1. System Framework. P1 generates a static day-ahead charging plan for the fleet (stage 1), while P2 and P3 determine realized vehicle dispatch and vehicle-charger assignment under dynamic system state evolution (stage 2).

3.2. Sequential decision making for dynamic ride-hailing systems

We model the dynamic ride-hailing system as a sequential decision process. This modelling framework is based on the discrete-event simulation framework (Horváth and Tamási, 2025). The system comprises a fleet of vehicles $V$, a set of operator-owned heterogeneous chargers $S$ located at various charging hubs/stations, and stochastic customer requests (exogenous). The operator observes exogenous electricity prices that are known one day in advance. We model the system dynamics with the following key components: system states, decisions, state transition functions, and an objective function.

**System states**: The system state refers to the physical and informational state of the system components over time. These components include vehicles, chargers, and customers. Let $\vec{x}_v(t), \forall v \in V$ denote the state of vehicle $v$ at time $t$, i.e., $\vec{x}_v(t) = (l_v(t), e_v(t), \chi_v(t), \psi_v(t))$ where $l_v(t)$ denotes the location of the vehicle, $e_v(t)$ the SoC of the vehicle, $\chi_v(t)$ the status of the vehicle, and $\psi_v(t)$ the assigned task to the vehicle at time $t$. A vehicle's status can be one of the following statuses: {*idle, serving, charging, waiting, en-route, arrive_chg*}. Note '*arrive_chg'* denotes the vehicle's status when arriving at its assigned charging station/ charger. Assigned tasks include {*nothing, serving, charging*}. Information related to customers' pickup and drop-off locations or assigned chargers to be used and their locations are associated with them. A customer (request) $r$'s state at time $t$ is denoted as $\vec{x}_r(t) = (o_r, d_r, A_r, W_r(t), \chi_r(t))$ where $o_r$ and $d_r$ are the corresponding pickup and drop-off locations. $A_r$ and $W_r(t)$ are customer's service request time and waiting times experienced from the request time until $t$, respectively. $\chi_r(t)$ represents the status of the customers in the system, which can be one of the following: {waiting, serving, completed, rejected}. Let $\vec{x}_s(t) = (l_s, b_s(t))$ denote the state of a charger $s$ at time $t$, where $l_s$ are the x- and y- coordinates of charger $s$, and $b_s(t)$ is the reservation/utilization state of $s$ (i.e., the union of scheduled charging sessions up to $t$). Customer arrivals are exogenous, described by the exogenous information status $\vec{x}_I(t)$, which reveals new customer arrivals up to time $t$. The overall system state at time $t$ is then denoted as $\vec{x}(t) = \left(\vec{x}_v(t), \vec{x}_r(t), \vec{x}_s(t), \vec{x}_I(t)\right), \forall v \in V, s \in S$. Sequential decisions, such as vehicle dispatch and charging, are triggered at different times and the system states are updated based on the transition functions.

**Decisions**: The operator makes the following three possible decisions at periodic decision points $t(i), \forall i \in \{1,2, \dots, T - T_0\}$.

- *Vehicle dispatch*: The decision to dispatch idle vehicles to serve customers.
- *Vehicle-charger assignment*: The decision to determine where and how much energy to charge for vehicles.
- *Vehicle-charger reassignment*: If the waiting time for a vehicle at a charger exceeds a predefined maximum, the vehicle is reassigned to another charger. This heuristic rule is only applied to the benchmark policies (described later) to improve the performance.

**State transition equations**: The state transition equations model the corresponding system state changes given a system state-event pair. The general form of transition functions is defined as $f: \boldsymbol{X} \times \boldsymbol{\mathcal{E}} \to \boldsymbol{X}$, where $\boldsymbol{X}$ is the set of all feasible system states $\{\vec{x}(t)\}$, and $\boldsymbol{\mathcal{E}}$ is the set of all events. Each event is associated with a time. Events are triggered and managed in an event queue according to their occurrence times. Events are processed in chronological order. The simulation clock jumps from one event to another and the system states are updated accordingly by the state transition equations. For a more detailed description, readers are referred to Appendix B.

**Objective function:** We calculate the total profit of the fleet over a 24-hour period to evaluate the system's performance. Revenue is generated by serving customers, while the costs include vehicle travel costs and charging expenses. The revenue received when serving a customer $r$ is based on the corresponding service fare $\beta_0 + \beta_1 \bar{d}_{o_r d_r}$ where $\bar{d}_{o_r d_r}$ is the trip distance of the request $r$. The service fare is composed of a base rate $\beta_0$ and a distance-based operating fee $\beta_1 \bar{d}_{o_r d_r}$. Given a trip (serving customers or traveling to charging stations) from $i$ to $j$, the vehicle travel cost is proportional to the travel distance, defined as $\pi \bar{d}_{ij}$, where $\pi$ is the cost per kilometer traveled of vehicles (dollar). Note that the trip set $\mathcal{T}$ includes all vehicle movements during service hours. As the overnight public charging stations are not represented in our model, the costs of the last trips to them are neglected. Given a charging session $q(\Delta e, u, t_{start}, t_{end})$ with $\Delta e$ kWh charged on charger $s$ of type $u_s$ from $t_{start}$ to $t_{end}$, the charging cost can be easily computed using a charging cost function $\phi(q)$. Considering the service hour period, let $\mathcal{R}$ denote the set of all served requests, $\mathcal{T}$ the set of all trips, and $\mathcal{Q}$ the set of all charging sessions that occurred during service hours. The **total daily profit** for a 24-hour period from $T_0$ to $T_0 + 1440$ (minutes) can then be calculated as the total revenue minus vehicle travel costs and charging costs by Eq. (1).

$$Z = \sum_{r \in \mathcal{R}} \left(\beta_0 + \beta_1 \bar{d}_{o_r d_r}\right) - \sum_{(i,j) \in \mathcal{T}} \pi \bar{d}_{ij} - \sum_{q \in \mathcal{Q}} \phi(q) - \sum_{v \in V} \bar{p}(B_v - e_{vT}) \tag{1}$$

The last term is the overnight charging costs, where $e_{vT}$ is the SoC of a vehicle $v$ at the end of the service time $T$.

Table 1 in the next section provides the sequence of P1-P3 model calls and vehicle and charger state updates of the above discrete-event sequential decision-making process where periodic decisions for vehicle dispatch and charging operations are made at each decision epoch and the system states are updated accordingly.

### 3.3. Day-ahead charging schedule planning (stage 1)

To better utilize limited operator-owned chargers, the day-ahead charging planning aims to minimize the daily charging operational costs by considering charging access costs, energy costs (including overnight charging costs), and opportunity costs for charging operations (including charging time and expected waiting time for charging). The day-ahead charging plan serves as a charging coordination guidance (when to charge and target SoC level after recharge) rather than a rigid plan to avoid charging congestion during peak hours. Let $V$ denote the set of vehicles, $\boldsymbol{T_0}$ **and** $\boldsymbol{T}$ **the start and end times of service hours**, respectively. The service hours from $T_0$ to $T$ is divided into a set of discrete charging decision epochs $h \in \widehat{H}^{\ell} = \{1,2,\dots,\lceil\frac{T-T_0}{\Delta_\ell}\rceil\}$ with a homogeneous interval $\Delta_\ell$, where $\Delta_\ell$ denotes the charging schedule planning decision time interval (e.g., 30 minutes). **The day-ahead charging planning aims to determine when vehicles charge and their target SoCs for each epoch in** $\boldsymbol{\widehat{H}^{\ell}}$ **given vehicles' driving (energy) needs and charging infrastructure constraints for the day.** As charging speed decreases significantly when the vehicle's SoC is above around 80% of its battery capacity (Froger et al., 2019), vehicles are not to charge above this threshold to save charging times during service hours. Based on historical driving data of vehicles, the average consumption of vehicles for each charging time slot can be calculated. Given the fact that vehicles are fully recharged at $T_0$, we can identify the first epoch at which the SoC of vehicles is below 80% (say $h'$). The precedent charging decision epochs of $h'$ $(h < h', \forall h \in \widehat{H}^{\ell})$ can be removed as the vehicles' charging decision at these epochs is irrelevant. In doing so, the feasible solution space can be reduced significantly. Let $H^{\ell}$ denote the set of the shifted decision epochs for charging schedule planning, i.e., $H^{\ell} = \{1,2,\dots,n_{H^{\ell}}+1\}$ where $n_{H^{\ell}}$ is the number of the charging decision epochs in $\widehat{H}^{\ell}$ where the SoC of vehicles is below 80%. The charging costs on fast chargers include a charger-type-specific fee $\tilde{p}_s$, where s is the index of chargers. The overnight charging cost in the objective function takes into account the charging costs for a full recharge at level-2 public charging stations with a cheaper flat fare. For a typical 62 kWh commercial EV battery[4], ride-hailing drivers can fully recharge their vehicles overnight at public level 2 charging stations with charging power ranging from 6 kW to 19.2 kW. The first epoch $h$ in $H^{\ell}$ corresponds to $h'$ in $\widehat{H}^{\ell}$. The problem is formulated as a MILP below in the space of $H^{\ell}$. The decision variables are: $x_{hs}^{v}$ vehicle $v$ is assigned to charger $s$ at epoch $h$, and $y_{hs}^{v}$ the amount of charged energy for a vehicle $v$ at charger $s$ at epoch $h$.

**P1: Day-ahead charging schedule planning**

$$\min Z_1 = \sum_{v\in V}\sum_{h\in H^{\ell}}\sum_{s\in S}\left((p_h + \tilde{p}_s + \frac{\gamma}{\varphi_s})y_{hs}^{v} + (C + \gamma\overline{W}_{hs})x_{hs}^{v}\right) + \sum_{v\in V}\bar{p}(B_v - e_{v,n_{H^{\ell}}+1}) \quad (2)$$

Subject to

$$\sum_{s\in S} x_{hs}^{v} \le 1, \forall v \in V, h \in H^{\ell} \quad (3)$$

$$\sum_{v\in V} x_{hs}^{v} \le 1, \forall s \in S, h \in H^{\ell} \quad (4)$$

$$e_{v,h+1} = e_{vh} - \delta_h\left(1 - \sum_{s\in S} x_{hs}^{v}\right) + \sum_{s\in S} y_{hs}^{v}, \forall v \in V, h \in H^{\ell} \quad (5)$$

$$x_{hs}^{v} = 1 \Rightarrow \left(\frac{y_{hs}^{v}}{\varphi_s}\right) \ge \alpha_s, \forall v \in V, h \in H^{\ell}, s \in S \quad (6)$$

$$E_v^{min} \le e_{vh} \le E_v^{max}, \forall v \in V, h \in H^{\ell} \cup \{n_{H^{\ell}}+1\} \quad (7)$$

$$e_{v1} = E_v^{0}, \forall v \in V \quad (8)$$

$$y_{hs}^{v} \le M_1 x_{hs}^{v}, \forall v \in V, h \in H^{\ell}, s \in S \quad (9)$$

$$0 \le y_{hs}^{v} \le Y_s^{max}, \forall v \in V_h, h \in H^{\ell}, s \in S \quad (10)$$

$$x_{hs}^{v} \in \{0,1\}, \forall v \in V, h \in H^{\ell}, s \in S \quad (11)$$

The objective function (2) minimizes total charging costs for a 24-hour period with additional overnight charging costs for a full charge. The first term in Eq. (2) is related to charging costs for $y_{hs}^{v}$ where $p_h$ denotes the average energy price on $h$. $\tilde{p}_s$ is charger-type-specific charging costs. $\varphi_s$ is the charging power of charger

[4] https://www.uber.com/nz/en/drive/services/electric/

s. $\gamma$ is the average profit per vehicle-minute traveled based on the realized customer service/profit on the previous days. The second term is related to charging access distance costs $C$ and expected waiting times $\bar{W}_{hs}$ when arriving at the charger $s$ at epoch $h$. The third term is the overnight charging costs from the SOC at the end of service hour to a full recharge (i.e., $B_v - e_{v,n_{H^\ell}+1}$). $\bar{p}$ is the flat fee (\$/kWh) charged for using public level-2 chargers. Eqs. (3) and (4) state that each vehicle can be assigned to at most one charger, and each charger can be assigned to at most one vehicle for each $h$, respectively. Eq. (5) states vehicles' SoC changes from $h$ to $h+1$ with the charged amount of energy when recharging and with average energy consumption $\delta_h$ of vehicles otherwise. $\delta_h$ represents the expected energy consumption of a vehicle during charging decision epoch $h$, estimated from historical energy consumption data of the fleet as the averaged energy consumption of a vehicle during epoch $h$. Eq. (6) states that a minimum charging time $\alpha_s$ (e.g. 10 minutes) is implied for each charging operation. Eq. (7) and (8) set the range of $e_{vh}$ and the initial battery level $E_0$ at $h = 1$, respectively. Eq. (9) binds $x_{hs}^v$ and $y_{hs}^v$. Eq. (10) ensures the maximum amount of energy can be recharged from charger $s$ during one charging decision epoch. $M_1$ is a positive number with $M_1 = max_{s\in S}\{Y_s^{max}\}$ where $Y_s^{max}$ is the maximum amount of energy that can be charged on charger $s$ during one charging decision epoch. The model parameters $\delta_h$, $\bar{W}_{hs}$, $\gamma$, and $C$ can be estimated based on historical vehicle driving and charging data. Note that in Eq. (6) (and other occurrence in this study), we use the symbol "⇒" to denote the if-then condition, which is handled as an indicator constraint in commercial solvers[5]. It means that if the expression on the left side of ⇒ is true, then the associated constraint on its right side is activated. An alternative formulation uses big-M tricks and is less efficient than this indicator-based approach when inappropriate big-M values are used.

**Based on the outputs of P1 (i.e., the solution of $e_{vh}$, $x_{hs}^v$ and $y_{hs}^v$), the operator obtains an approximate schedule to assign vehicles to charge and their target SoCs after recharge over $\hat{H}^\ell$. This schedule is then adapted based on an online reactive model to determine the amount of energy to be charged and where to charge for vehicles during the day.** This reactive model maintains a pool of go-charge-vehicles based on the day-ahead charging plan (excluding vehicles currently serving customers), which is adapted with additional vehicles to recharge; either with low SoC (i.e. less than 20% of their battery level) or previously delayed vehicles for charging or due to the number of go-charge-vehicles exceeds the number of chargers. An online vehicle-to-charger assignment model below is applied to minimize total charging operational times for vehicle-to-charger assignment. A more detailed description is presented in the simulation framework (Table 1). The online vehicle-to-charger assignment problem is formulated as a MILP as follows.

### 3.4. Vehicle dispatch and vehicle-charger assignment (stage 2)

We adopt a batch dispatch optimization approach to match unserved requests with idle vehicles. Customer arrivals are grouped into batches at each decision epoch $i \in H^{\mathscr{b}} = \{1,2, \dots T - T_0\}$ with one-minute time interval. We denote $t(i)$ the start time of decision epoch $i \in H^{\mathscr{b}}$. A batch dispatching optimization is executed at each batch decision epoch to maximize the profit of serving customers. As the stochastic demand is revealed in real-time, the operator's profit maximization problem cannot be optimized for the entire day, but instead must be optimized in a sequential way based on customer arrivals every minute. Customer patience is represented by a maximum acceptable waiting time, which is integrated into the vehicle dispatching model to maximize the total profit of vehicles' dispatches. If a customer's waiting time exceeds this patience threshold, the customer is considered unserved and leaves the system. Let $R_{t(i)}$ denote the pool of unserved requests at time $t(i)$, and $V_{t(i)}$ the set of idle vehicles at $t(i)$. The batch dispatching problem at epoch $i$ (i.e., $i$ minutes from the start of the planning horizon as a one-minute time interval is used) is formulated as a MILP problem. The decision variable is $m_{rv}$, determining the vehicle-to-request assignment. $e_{vt(i)}$ and $W_{rt(i)}$ are the SoC of vehicle $v$ and the cumulative waiting time of customer $r$ at time $t(i)$, $\forall i \in$

[5] https://docs.gurobi.com/projects/optimizer/en/current/concepts/modeling/constraints.html

$H^{\mathcal{b}}$, respectively. The problem is solved every minute, and the total revenue is aggregated over the service hours.

**P2: Batch dispatch**

$$\max Z_2^i = \sum_{r \in R_{t(i)}} \sum_{v \in V_{t(i)}} \left(\beta_0 + \beta_1 \bar{d}_{o_r d_r} - \pi(\bar{d}_{o_r d_r} + \bar{d}_{v o_r})\right) m_{rv} \tag{12}$$

Subject to

$$\sum_{v \in V_{t(i)}} m_{rv} \leq 1, \forall r \in R_{t(i)} \tag{13}$$

$$\sum_{r \in R_{t(i)}} m_{rv} \leq 1, \forall v \in V_{t(i)} \tag{14}$$

$$m_{rv} = 1 \Rightarrow E_v^{min} \leq e_{vt(i)} - \mu\left(\bar{d}_{o_r d_r} + \bar{d}_{v o_r}\right), \forall r \in R_{t(i)}, v \in V_{t(i)} \tag{15}$$

$$W_{rt(i)} + \bar{t}_{v o_r} m_{rv} \leq W^{max}, \forall r \in R_{t(i)}, v \in V_{t(i)} \tag{16}$$

$$m_{rv} \in \{0,1\}, \forall r \in R_{t(i)}, v \in V_{t(i)} \tag{17}$$

The objective function (12) maximizes the profit of vehicle dispatch at time $t$ where the net profit of a customer-vehicle match $(r, v)$ is calculated as the difference of the service fare $\beta_0 + \beta_1 \bar{d}_{o_r d_r}$ and vehicle's travel cost $\pi\left(\bar{d}_{o_r d_r} + \bar{d}_{v o_r}\right)$, where $\bar{d}_{o_r d_r}$ is the trip distance of request $r$, and $\bar{d}_{v o_r}$ is the distance from the vehicle's current location to pick up customer $r$ at $o_r$. Constraints (13)-(14) ensure that one vehicle can serve at most one customer and vice versa. Constraint (15) ensures that a matched vehicle needs to have sufficient energy to reach the pickup location of the assigned customer and serve that trip (i.e., the vehicle's SoC needs to be no less than a minimum level after subtracting the energy consumption of serving that ride from its current SoC). Energy consumption is assumed proportional to the travel distance with a constant energy consumption rate $\mu$. Constraint (16) ensures customers' waiting time cannot exceed a maximum threshold $W^{max}$ (e.g., 7-10 minutes). Note that $W_{rt(i)}$ is customer $r$'s cumulative waiting time up to $t$. Note that unassigned customers wait in the system to be assigned until the next decision time $t(i+1)$. If a customer's maximum waiting time is reached, they leave the system and are not served by the service. For unmatched vehicles (remaining idle), their SoCs remain unchanged until $t(i+1)$.

**P3: Online vehicle-charger assignment**

We utilize several different sets to manage online vehicle-to-charger assignment. Let $\Omega_{t(i)}$ be the set of need-to-charge vehicles at decision time $t(i)$, $\forall i \in H^{\mathcal{b}}$. It includes the set of vehicles that are scheduled for recharging based on the day-ahead charging plan whose charging requirements have not been implemented up to $t(i)$ (e.g., due to serving customers), denoted as $\Omega_{t(i)}^0$, and the set of idle vehicles**,** not in $\Omega_{t(i)}^0$, whose actual SoCs are below the charging threshold, denoted as $N_{t(i)}$. Let $\widetilde{\Omega}_{t(i)}$ be the subset of $\Omega_{t(i)}$ that are relevant to be assigned to chargers for recharging at time $t(i)$, i.e. a subset of vehicles in $\Omega_{t(i)}$ that **are idle** at time $t(i)$ and need to charge to their target SoCs satisfying a minimum amount of charged energy requirement, calculated as the product of the fast-charger power and the corresponding minimum charging duration per charging operation. We denote $M_{t(i)}$ as the set of vehicles in $\widetilde{\Omega}_{t(i)}$ whose charger assignments are successfully implemented at $t(i)$ (described below). Therefore, $\Omega_{t(i+1)} = \{\Omega_{t(i)} \backslash M_{t(i)}\} \cup N_{t(i+1)} \cup \widetilde{V}_{h|h=t(i+1)/\Delta h}$ where $\widetilde{V}_{h|h=t(i+1)/\Delta h}$ is the newly added need-to-charge vehicles based on the day-ahead plan at time $t(i+1)$, if any. At each decision point $i$ for vehicle (batch) dispatch, the operator needs to decide whether currently idled vehicles in $\widetilde{\Omega}_{t(i)}$ need to recharge. If not, the remaining idled vehicles are available to serve customers. Given the current decision time $t(i)$, the current location of vehicles and the utilization state of chargers, the objective function (18) aims to minimize the total charging access time ($\bar{t}_{vs}$) and waiting time ($\bar{W}_{vs}^{t(i)}$), and charging time ($y_{vs}/\varphi_s$) of the assignment of vehicles to chargers. Note that while the day ahead planning of P1 optimizes charging costs, here P3 seeks to minimize charging time, prioritizing vehicle availability to avoid the cost of missed rides. The decision variables are: $\hat{x}_{vs}$ whether a vehicle $v$ is assigned to a charger $s$, and $\hat{y}_{vs}$ the amount of energy vehicle $v$ charges at a charger. $e_{vt(i)}$ is the SoC of vehicle $v$ at decision time $t(i)$. Given the current utilization state of chargers, $\bar{W}_{vs}^{t(i)}$ is the waiting time that vehicle $v$ would experience if it departs immediately to go to charger $s$ at decision time $t(i)$. The operator looks across all current charger queues and obtains the waiting time for every charger. This information can be obtained online from the operator's charging network management system.

$$\min Z_3^i = \sum_{v\in\widetilde{\Omega}_{t(i)}} \sum_{s\in S} ((\bar{t}_{vs}+\bar{W}_{vs}^{t(i)})\hat{x}_{vs} + \hat{y}_{vs}/\varphi_s) \tag{18}$$

$Subject\ to$

$$\sum_{s\in S} \hat{x}_{vs} = 1, \forall v \in \widetilde{\Omega}_{t(i)} \tag{19}$$

$$\sum_{v\in\widetilde{\Omega}_{t(i)}} \hat{x}_{vs} \leq 1, \forall s \in S \tag{20}$$

$$\hat{x}_{vs} = 1 \Longrightarrow e_{vt(i)} - \mu\bar{d}_{vs}\hat{x}_{vs} \geq 0,\ \forall v \in \widetilde{\Omega}_{t(i)}, s \in S \tag{21}$$

$$\hat{x}_{vs} = 1 \Longrightarrow e_{vt(i)} - \mu\bar{d}_{vs}\hat{x}_{vs} + \hat{y}_{vs} \geq \bar{E}_v,\ \forall v \in \widetilde{\Omega}_{t(i)}, s \in S \tag{22}$$

$$\hat{x}_{vs} = 1 \Rightarrow \left(\frac{\hat{y}_{vs}}{\varphi_s}\right) \geq \alpha_s, \forall v \in \widetilde{\Omega}_{t(i)}, s \in S \tag{23}$$

$$\hat{y}_{vs}/\varphi_s \leq \tau^{max}\hat{x}_{vs},\ \forall v \in \widetilde{\Omega}_{t(i)}, s \in S \tag{24}$$

$$0 \leq \hat{y}_{vs} \leq (B - e_{vt(i)} + \mu\bar{d}_{vs})\hat{x}_{vs},\ \forall v \in \widetilde{\Omega}_{t(i)}, s \in S \tag{25}$$

$$\hat{x}_{vs} \in \{0,1\}, \forall v \in \widetilde{\Omega}_{t(i)}, s \in S \tag{26}$$

When the number of to-recharge vehicles does not exceed the number of chargers, Eqs. (19) and (20) ensure that each vehicle is assigned to exactly one charger, and each charger can be connected to at most one vehicle. In the other case, these two equations are replaced by Eqs. (27)-(28).

$$\sum_{s\in S} \hat{x}_{vs} \leq 1, \forall v \in \widetilde{\Omega}_{t(i)} \tag{27}$$

$$\sum_{v\in\widetilde{\Omega}_{t(i)}} \hat{x}_{vs} = 1, \forall s \in S \tag{28}$$

Eq. (21) ensures that the vehicle's SoC is always non-negative when arriving at the charger's location. Eq. (22) states that vehicle $v$ needs to be recharged to the target SoC $\bar{E}_v$. Note that the values of $\bar{E}_v$ for vehicles in $\Omega_{t(i)}^0$ are based on the day-ahead charging plan. When a vehicle is delayed in being added to the pool due to serving customers, its target SoC remains the planned one based on the output of P1. However, for vehicles in $N_{t(i)}$, their target SoCs are adapted to reflect anticipated future energy needs. These vehicles' respective target SoCs are set based on the average energy consumption $\delta_h$ (estimated from historical data) by anticipating their energy consumption to the end of service hours. To maximize vehicles' availability, vehicles' SoC at the end of service hours should be as close as possible to $E_v^{min}$. For this purpose, we apply the following rule to determine vehicles' target SoCs when their SoCs are below $\theta B$. Let $SoC_v(t(i))$ denote vehicle $v$'s SoC at time $t(i)$, and $h(t(i))$ be the corresponding $h$ index of $t(i)$. The target SoC of vehicle $v$ when adding it to the to-charge vehicle pool at $t(i)$ is defined as

$$SoC_v^{target}(t(i)) = min\,(\tilde{E}_v, SoC_v(t(i)) + \vec{E}_v(t(i))) \quad (29)$$

where $\vec{E}_v(t(i))$ is the energy needed to the end of service hours $T$, calculated as $\vec{E}_v(t(i)) = \sum_{h=h(t(i))}^{n_{H^\ell}} \delta_h - m\bar{\delta}_{tail}$ with $m$ being the approximated number of epochs whose total energy consumption is around $E_v^{min}$(i.e. $m\bar{\delta}_{tail} \leq E_v^{min}$), and $\bar{\delta}_{tail}$ is the average $\delta_h$ over the last $m$ epoch until the end of service hours. In doing so, Vehicles' SoCs at the end of service hours would be a little more above $E_v^{min}$, given vehicles' current SoC is around $\theta B$ where $B$ denotes the vehicle's battery capacity. To further reduce vehicles' waiting time for charging, we uniformly randomly draw $\tilde{E}_v$ from the interval $[0.5B, E_v^{max}]$. This policy is more flexible than using $E_v^{max}$ as vehicles can recharge again sometime later as far as their SoCs are below $\theta B$. Eq. (23) ensures a minimum charging time $\alpha_s$ when vehicles go charging. $\varphi_s$ is the charging rate of charger $s$. Eq. (24) binds the variables $\hat{x}$ and $\hat{y}$, and limits the charging duration of each charging session to a reasonable upper bound, represented by a predefined maximum charging duration (e.g., 60 minutes). This parameter can be tuned according to the interplay between charging-infrastructure supply, charging demand, and operational performance. Eq. (25) ensures that the amount of energy charged does not exceed the remaining vehicle battery capacity upon arrival at the assigned charger.

Once the P3 has been solved, if a vehicle $v$ is assigned to $s$ for recharge and its expected waiting time $\overline{W}_{vs}^{t(i)}$ exceeds a maximum waiting time limit (e.g., 30 minutes), the charging assignment is not implemented, and the vehicle remains idle to serve customers. Otherwise, the charging request is accepted, and the charger reservation scheduler updates the utilization/reservation state of the corresponding charger immediately by adding the scheduled charging sessions. Consequently, a charging assignment accepted at the current decision epoch is incorporated into the charger state used at subsequent decision epochs, preventing overlapping charging sessions and allowing the waiting time estimation to account for previously accepted assignments. P3 nevertheless remains myopic: it assigns only the vehicles currently in the pool and does not anticipate vehicles that will need a charger shortly afterwards, so a charger reserved now may be unavailable to a vehicle with a more urgent need a few minutes later. Re-optimising every minute limits this effect; a rolling-horizon version of P3 that includes forecast charging requests would be a possible extension.

Note that if no feasible solutions can be found for $\tilde{\Omega}_{t(i)}$, vehicles with the highest SoCs are removed from $\tilde{\Omega}_{t(i)}$, and then the problem is solved again until the optimal solution is found. P3 can be infeasible if a vehicle has insufficient energy to reach an eligible charger, cannot charge to its planned target energy level within the allowable maximum charging time per session, or there are too few eligible chargers. The maximum charging waiting time restriction is checked after solving P3 and therefore does not cause MILP infeasibility. Instead, vehicle charging assignments of P3 that exceed this maximum waiting time limit are deferred to avoid excessive waiting at charging stations.

If the removed vehicles are not dispatched for serving customers (remain idle) at $t(i)$, they are retained in the pool $\Omega_{t(i+1)}$ for charging at the next decision time $t(i+1)$. Note that for each vehicle-to-charger assignment epoch $i$, $\Omega_{t(i)}$ is filtered by retaining a subset of vehicles in $\Omega_{t(i)}$ where the calculated amount of energy to be charged (depending on vehicles' current SoC and their target $\bar{E}_v$) needs to be at least equivalent to the amount of charging $\alpha$ minutes (minimum charging duration) on a fast charger. If the amount of energy to be charged is below this minimum amount of energy, vehicles remain idle for serving customers. In doing so, multiple short-duration recharging operations with short charging times can be avoided, significantly reducing the operator's total charging access costs. Note the sets of $\tilde{\Omega}_t$, $\Omega_t$ and $V_t$ depend on system states or decisions and are updated over time (see details in the simulation pseudocode in Table 1).

**Table 1 shows the pseudocode for implementing the discrete-event-based, sequential decision-making framework described in Section 3.2**. This framework governs how the system state evolves over time at each decision epoch/point. Specifically, the system evolves when new customers arrive and when periodic decisions are made at each discrete decision time (every minute) for vehicle dispatch and charging operations. The **interdependence of Models P1-P3** is as follows: (i) day-ahead charging planning (P1)

provides initial charging plans and target SoC after recharge of vehicles to guide vehicle charging under stochastic demand (lines 3 and 7-9); (ii) vehicle batch dispatch (P2) assigns idle vehicles to service unserved requests (lines 15-16); (iii) online vehicle-charger assignment coordinates charging operations in response to current system state (lines 10-14).

At each discrete decision epoch $i$, decisions are made based on the updated system state up to $t(i)$. The resulting decisions evoke sequential state transitions according to the sequential decision process described in Section 3.2. Table 1 summarizes the different steps in the simulation. Step 1 reads the input data. Step 2 estimates the parameters used for the day-ahead charging planning model (described in more detail in Sect. 4). Step 3 solves P1 to obtain the day-ahead charging plan of vehicles. Step 4 sets up the initial conditions for the simulation. Step 5 involves the simulator clock loop for batch dispatching where $t(i)$ corresponds to $i$ minutes after the start time (6:00), and $T$ corresponds to the end of service hours (20:00). Note that we use continuous time to track the system state in the simulator. Steps 6 to 8 add vehicles to the pool based on the day-ahead charging plan. Step 10 adds additional vehicles to the pool if their SoCs are below the threshold $\theta B$. Steps 10-14 filter idle vehicles in the pool and assign them to chargers for recharge by solving P3. Steps 15-16 update the list of idle vehicles and unserved customers and dispatch vehicles to serve customers based on P2. Then update the system state to the beginning of next decision point $t(i+1)$. Step 18 computes the overnight charging costs from the vehicles' SoCs at the end of service hours to full capacity, based on the flat fare $\bar{p}$. Step 19 outputs the total profit of the fleet for the 24-hour period. Note that Steps 7-16 represent charging-first-and-service-second decision priorities. To justify this decision ordering setting, an online charging-first/service-first decision-order sensitivity analysis is conducted in Appendix I.

Table 1. Simulation-based implementation of the sequential decision-making framework for dynamic charging planning, vehicle-to-charger assignment, and batch vehicle dispatch.

1. Input: Time-dependent energy prices, customer demand, a fleet of vehicles, and charging facilities.
2. Compute the average energy consumption ($\delta_h$) and average charging waiting times ($\overline{W}_{hs}$) of vehicles at chargers per charging decision epoch up to the previous day.
3. Solve the day-ahead charging planning problem **P1** to get initial charging plans and target SoC after recharge of vehicles ($\tilde{V}_h, \bar{E}_v, v \in \tilde{V}_h$) for the planning horizon (i.e. 6:00-20:00).
4. Initialization: Initialize SoCs and locations of vehicles. Initialize the set of idle vehicles $V = \emptyset$. Set the pool of need-to-charge vehicles as empty, $\Omega = \emptyset$.
5. **For** the vehicle dispatching decision points, $i \in \{1,2,\dots,T-T_0\}$
6.   Update the system state to $t$
7.   **If** $t(i)$ mod $\Delta h = 0$
8.     Add the planned need-to-charge vehicles $\tilde{V}_{h|h=t(i)/\Delta h}$ to $\Omega_{t(i)}$
9.   **End**
10.   Find the subset $N_{t(i)}$ of idle vehicles in $V_{t(i)}$ with SoCs lower than $\theta B$ and add them to $\Omega_{t(i)}$.
11.   Find the subset of vehicles $\tilde{\Omega}_{t(i)}$ in $\Omega_{t(i)}$ that are idle and must charge to their target SoCs satisfying a minimum amount of charged energy requirement.
12.   **If** $\tilde{\Omega}_{t(i)}$ is not empty
13.     **Solve P3** for $\tilde{\Omega}_{t(i)}$ and obtain vehicle-charger assignments. If a vehicle’s waiting time exceeds the maximum waiting time limit, its charging assignment is not implemented, and the vehicle remains idle to serve customers. Otherwise, the charging request is accepted and assign vehicles to their assigned chargers, execute their charging operations, and update vehicle state and the utilization/reservation state of the corresponding chargers. When a vehicle completes its charging operation at time $t'$, it becomes idle immediately and is added to the set of idle vehicles. If there are no feasible solutions, remove the vehicle with the highest SoC from $\widetilde{\Omega}_{t(i)}$ and solve P3 again. Continue until a feasible solution can be found.

14. **End**
15. Update the lists of unserved customers $R_{t(i)}$. Remove vehicles going for charging $(M_{t(i)})$ from the list of idle vehicles $V_{t(i)}$ and $\Omega_{t(i)}$.
16. Solve the batch dispatching problem **P2**, dispatch vehicles to serve customers, and update the state of vehicles (SoCs and locations). Update $V_{t(i)}$ by removing these vehicles at time $t(i)$ and update the corresponding sets of idle vehicles when dropping off customers.
17. **End**
18. Compute the overnight charging costs $\sum_{v \in V} \bar{p}(B_v - e_{vT})$
19. Output the profit of the fleet for a 24-hour period

**Illustrative example.** To make the decision sequence and the dispatch–charging interaction concrete, consider a single vehicle $v$ and two fast chargers $s_1$ and $s_2$ at the 14:00 decision epoch $t(i)$. Suppose $v$ is idle with SoC $e_v(t)$ = 12 kWh, just below the $\theta B$ = 20% trigger for a 62 kWh battery, and that the day-ahead plan (P1) has assigned $v$ a target SoC $\bar{E}_v$= 35 kWh for this epoch. First, the vehicle is added to the need-to-charge pool, and the vehicle-charger assignment model P3 chooses between $s_1$ and $s_2$ by minimizing the total charging operational time (access time + expected queueing time + charging time). Suppose $s_1$ is closer but already holds two reserved sessions, so the reservation scheduler returns an expected wait of 25 minutes, while $s_2$ is slightly farther but has no queue. P3 then assigns $v$ to $s_2$; the scheduler reserves the session $\zeta(v, t_{start}, t_{end})$ and updates the charger's reservation state so that no later request can overlap it. If instead both chargers returned waits above the 30-minute limit, $v$ would remain idle and available to serve customers. The online controller then solves the batch-dispatch model P2: if an unserved request $r$ can be reached within the customer's maximum waiting time and $v$ has enough energy to complete the trip while staying above $E_v^{min}$(i.e., 0.1$B$), then $v$ is matched to $r$ and its charging is deferred to be reconsidered for charging at the next epoch. This vignette highlights the central interaction: at every epoch, dispatching and charging compete for the same vehicle, and the day-ahead target SoC $\bar{E}_v$ guides rather than forces when and how much $v$ recharges.

### 3.5. Benchmark charging policies

To validate the proposed methodology, five benchmark charging policies selected from the literature are compared. The benchmark charging policies can be classified into rule-based policies (benchmark policies from a. to d. below) and optimization-based policy (benchmark policy e. below). The rule-based policies assume that vehicles go to recharge when their SoC is lower than a pre-defined threshold. Furthermore, these policies apply different methods to determine when, where, and how much energy to charge. Existing studies assume vehicles are willing to wait without limits at a charger (Jamshidi et al., 2021). We consider more realistic queuing scenarios at chargers for the benchmark policies by assuming a maximum waiting time limit in a queue (i.e., **30 minutes, identical for all vehicles and for all charging policies**). Vehicles are assigned to chargers based on the applied benchmark charging policy. When arriving at the assigned charging stations, if the waiting time exceeds the maximum threshold, vehicles move away to another charger with the least waiting time when their SoCs are sufficient to reach there. In case the vehicle's SoC is too low to reach the targeted charger, vehicles go to the nearest one (if the vehicle's SoC is feasible) or remain at the same charger (if the vehicle's SoC is not feasible). This allows not to have unrealistically too long queue at charging stations. **The benchmark policies are aware of the queueing information, same as the CongestionAware charging policy**. The benchmark charging policies are described as follows.

a. Nearest charging policy (**Nearest**) (Bischoff and Maciejewski, 2014): Vehicles go to the nearest charger to recharge to $E_v^{max}$ (i.e., 80% of their battery capacity) when their SoCs are below the threshold $\theta B$ (i.e. 20% of their battery capacity).
b. Fastest charging policy (**Fastest**): Vehicles go to the fastest charger to recharge to $E_v^{max}$ when their SoCs are below the threshold $\theta B$. In case there is more than one fastest charger, a randomly selected one is used. Note that if charging demand on fast chargers is not high (no congestion), one can assign vehicles to the nearest one. On the contrary, when the number of fast chargers is scarce, using this random-assignment policy allows not to assign too many vehicles to a geographically well-situated charging station (e.g., the one located at the middle of our study area), resulting in over-saturated utilization of certain fast chargers/charging stations.
c. Charging operational time minimization approach (**MinChgOpT**) (Ma and Xie, 2021): When vehicles' SoCs are below the threshold, vehicles are assigned to the charger with a minimum charging operational time to charge to $E_v^{max}$, including access time, waiting time when arriving at chargers and charging time.
d. Dynamic charging threshold policy (**DynaThreshold**) (Ahadi et al., 2023): Different from the above benchmark policies, DynaThreshold activates vehicles' charging operations earlier when their SoCs are still much higher than the threshold in order to avoid charging during peak charging demand periods in the afternoon. In doing so, vehicles can save (charging) waiting time and the charging facility can be utilized more effectively. We adopt the hourly dynamic charging thresholds used in Ahadi et al. (2023), where the average hourly charging threshold in the morning is around 50% while it becomes around 30% in the afternoon (see Table 2). When going charging, vehicles are assigned to the fastest charger to charge to $E_v^{max}$ so we expect that vehicles charges smaller amount of energy compared to the situation in the afternoon. In case there is more than one fastest charger, a randomly selected one is used.
e. Optimal charging scheduling and assignment policy (**OptChg**) (Ma, 2021): This policy is a similar two-stage optimization approach developed in our previous study for electric ride-hailing services. Unlike our CongestionAware policy, this policy obtains a day-ahead charging plan based on vehicle-specific energy consumption on each charging decision epoch without considering charging station capacity constraints and heterogeneous charging stations. We adapt this model and formulate an alternative MILP model to obtain the day-ahead charging plan (the corresponding MILP model based on Ma (2021) can be found in Appendix F). For the online vehicle-charger assignment model, we use the same model P3 and vehicles' energy need anticipation approach as the CongestionAware policy as it significantly outperforms the online vehicle–charger assignment approach used in Ma (2021), which relies on a full recharge to $E_v^{max}$ policy and without considering the vehicle's minimum charging duration requirement.

Note that we assume that the operator has real-time access to occupancy, queuing, and waiting time information for each charger. **All fast and slow chargers are available for benchmarking using the same experimental setup as the CongestionAware charging policy**. The vehicle-customer dispatch model (P2) is used for all policies. The rule-based benchmark policies do not use the day-ahead plan, whereas OptChg employs an alternative day-ahead charging planning model adapted from Ma (2021). OptChg uses the same vehicle-charger assignment model P3 and energy-need anticipation approach as the CongestionAware policy. Therefore, the comparison between CongestionAware and OptChg allows us to isolate the effect of the day-ahead charging planning component. To ensure a fair comparison, the profit calculation for all benchmark policies includes the cost to recharge vehicles to 100% SoC overnight.

Table 2. Dynamic hourly SoC threshold for activating charging operations (Ahadi et al., 2023).

| Hour (morning) | 1 | 2 | 3 | 4 | 5 | 6 | 7 | 8 | 9 | 10 | 11 | 12 |
|---|---|---|---|---|---|---|---|---|---|---|---|---|
| SoC threshold ($\theta$)* | 0.45 | 0.6 | 0.65 | 0.62 | 0.58 | 0.55 | 0.52 | 0.5 | 0.4 | 0.4 | 0.4 | 0.4 |
| Hour (afternoon) | 13 | 14 | 15 | 16 | 17 | 18 | 19 | 20 | 21 | 22 | 23 | 24 |
| SoC threshold ($\theta$) | 0.38 | 0.35 | 0.32 | 0.25 | 0.25 | 0.2 | 0.2 | 0.25 | 0.27 | 0.35 | 0.35 | 0.4 |

*: % of vehicles' battery capacity.

## 4. Computational study

In this section, we first describe the test instance and parameter setting based on Manhattan yellow taxi data. Then we present the computational results for different demand scenarios. A sensitivity analysis is conducted to evaluate the impact of different model parameters.

### 4.1. Test instance generation

We test the proposed approach on a Manhattan-like 4×20 km² area divided into 80 one-kilometre zones, with a fleet of 100 homogeneous EVs. Service hours are 06:00–20:00; between 20:00 and 06:00 the next day, drivers recharge to full at public level-2 charging stations. We consider two demand scenarios: 3000 customers/day (low) and 4000 customers/day (high), both over 06:00–20:00. For reference, the average yellow taxi made 20.5 trips/day in July 2019 according to the NYC Taxi and Ride-hailing Usage dashboard[6]. Demand is generated from the New York yellow-taxi data[7] for a typical weekday in July 2019. From this dataset we take only the customer pickup times, which preserves a realistic temporal demand profile (notably the evening peak); origins and destinations are generated synthetically within the study area, since the method's focus is charging coordination rather than spatial demand modelling. Generating the spatial demand lets us evaluate the policies under controlled, repeatable conditions instead of replaying a single historical day: the demand level is varied systematically (3000 and 4000 customers/day, extended to 1000–6000 in Appendix D), and independent realizations are drawn so that performance differences are attributable to the charging logic rather than to one day's trip pattern. For each scenario we generate 15 independent demand datasets: 10 for estimating the P1 parameters ($\delta_h$, $\overline{W}_{hs}$ and $\gamma$) and 5 for the test of the performance of the proposed methodology. Origins lie within the 4×20 km² area while destinations may fall outside it, and we impose a minimum trip length of 5 km, which is somewhat longer than typical Manhattan trips.

Figure 2 illustrates the charging station locations and initial locations of vehicles in the study area. For the operator-owned charging infrastructure, we assume there are 6 fast and 6 slow chargers located at 4 different charging stations (i.e. each fast (slow) charging station has three fast (slow) chargers, see Figure 2). We locate one fast charging station around the bottom center and another around the center right based on the current public charging station locations in Manhattan[8]. The charging power of the operator-owned dedicated charging stations is assumed 50kW (fast) and 11kW (slow), respectively. As aforementioned (Section 3.1), the vehicles are charged on operator-owned dedicated charging stations during the service hours and on public level-2 chargers for overnight charging.

The number of customer arrivals per 10 minutes for low and high-demand scenarios is shown in Figure 3. We can observe that both demand profiles are very irregular with higher peaks during certain time slots (e.g. 18:00-20:00). The peak demand is around 80 and 50 customer arrivals per 10 minutes for the high and low demand scenario, respectively. We assume that Nissan Leaf e+ is the used vehicle with a battery capacity of 62 kWh and an energy consumption rate of 0.25 kWh/kilometer traveled. The impact of battery capacity on the performance of the proposed methodology will be analyzed in the sensitivity analysis. For energy prices, we assume that ToU energy prices change every 15 minutes based on the time-varying day-ahead electricity prices, adapted to the real charging fee of NYC public charging stations. Charging costs on the operator-owned chargers depend on the type of charger used, i.e., there is an additional fee for using fast chargers. We assume the cost of $0.2/kWh for using fast chargers[9].

[6] https://toddwschneider.com/dashboards/nyc-taxi-ridehailing-uber-lyft-data/
[7] https://www.nyc.gov/site/tlc/about/tlc-trip-record-data.page
[8] https://www.nyc.gov/html/dot/html/motorist/electric-vehicles.shtml#/find/nearest?location=NYC

For a level-2 charger in NYC, the charging cost is $1 per hour for overnight charging[10]. We estimate the cost of overnight charging as $1/6.24kWh=$0.16/kWh. The typical charging power is 6.4kW to 7.2kW or higher 19.2 kW[11] in NYC, and there are an estimated 2000 level-2 chargers[12]. We assume drivers can fully charge their vehicles overnight from 8 p.m. to 6 a.m. the next day.

Figure 4 shows the ToU electricity prices per kWh charged during the service hours. The charging costs are in the range of $0.3-$0.7/kWh for using the fast chargers and around $0.1- $0.5/kWh for slow chargers. Note that other practical ToU energy prices can be applied in practice. The detailed parameter setting for the computational study is shown in Table 3.

Table 3. Simulation parameter settings of the case study.

| System parameter | Value |
|---|---|
| Fleet size | 100 |
| Number of charging stations | 4 |
| Number of chargers[1] | 6 DC fast (50kW) and 6 slow chargers (11 kW) |
| Number of customers per day | 3000 and 4000 |
| Vehicle speed | 30 km/hour |
| Taxi fare[2] | $\beta_0 = \$8$ (base rate), $\beta_1 = 3.1$ ($/km operating fee) |
| Maximum waiting time of customers | 10 minutes |
| Minimum charging duration ($\alpha_s$) | 10 minutes |
| Battery capacity[3] ($B_v$) | 62 (kWh) |
| Travel distance cost[4] ($\pi$) | 0.53 (dollars/km) |
| Energy consumption rate[3] ($\mu$) | 0.25 (kWh/km) |
| $E_{max}(E_{min})$ | 80% (10%) of $B_v$ |
| $E_{init}$ | 62 (kWh) |
| Time-dependent energy price[5] | Day-ahead electricity Prices (15-minute resolution) |
| Service time | 6:00-20:00 |
| Charge decision time interval ($\Delta h$) | 30 minutes for the c3000 scenario and 20 minutes for the c4000 scenario[6] |
| Maximum charging duration per session ($\tau^{max}$) | 60 minutes[7] |
| Batch dispatching time interval ($\Delta_t$) | 1 minute |
| $\tilde{p}_{fast_charger}$ | 0.2 dollars/kWh |
| $\bar{p}$ (flat fee for overnight charging) | 0.16 dollars/kWh |

Remark: 1. Based on the two main charger types used in New York City (Taxi & Limousine Commission, 2022). 2. Approximated based on the current yellow taxi fare in NYC. The base fee includes the surcharges (https://www.nyc.gov/site/tlc/passengers/taxi-fare.page). 3. Adapted from the characteristics of Nissan Leaf e+ (https://evadept.com/calc/tesla-supercharger-charging-cost-calculator). 4. Based on Ahadi et al. (2023). 5. Adapted from https://transparency.entsoe.eu/, bidding zones GE_LU, July 2019. The adapted electricity price variation range covers the real charge fee in New York City ($0.39 per kWh consumed; see https://www.nyc.gov/html/dot/html/motorist/electric-vehicles.shtml#/dc). 6. The charging decision interval depends on factors such as EV battery capacity and fast-charger power. Reducing the decision epoch from 30 to 20 minutes increases the number of charging decision epochs per hour by 50%, thereby enabling P1 to obtain feasible solutions for the c4000 scenario. This parameter can therefore be calibrated to achieve an appropriate balance between feasibility, computational time, and system profit. The impact of different decision interval lengths for the c3000 and c4000 scenarios is provided in Appendix I. 7. This value is based on a sensitivity analysis considering values ranging from 30 minutes to beyond 60 minutes. The results indicate that imposing a 60-minute limit has

[9] https://www.blackcarnews.com/article/gravity-installs-super-fast-ev-chargers-in-nyc?utm_source=chatgpt.com

[10] Curbside Level 2 Charging Project FAQ. https://www.nyc.gov/html/dot/downloads/pdf/curbside-level-2-charging-pilot-faq.pdf

[11] https://www.flo.com/new-york-city/

[12] https://chargehub.com/en/countries/united-states/new-york/new-york.html?utm_source=chatgpt.com

a negligible impact on profitability compared with allowing unlimited charging durations, resulting in profit reductions of only 0.1% and 0.2% for the c3000 and c4000 scenarios, respectively.

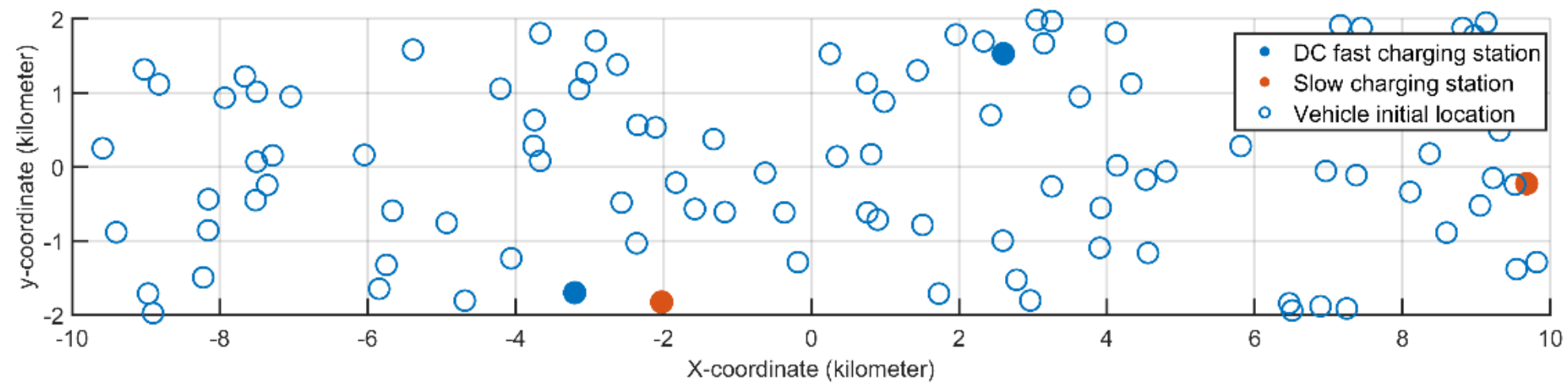

Figure 2. Charging station distribution in the study area. There are 4 charging stations, each with three chargers of the same type. A total of 6 DC fast chargers (50kW) and 6 slow chargers (11kW).

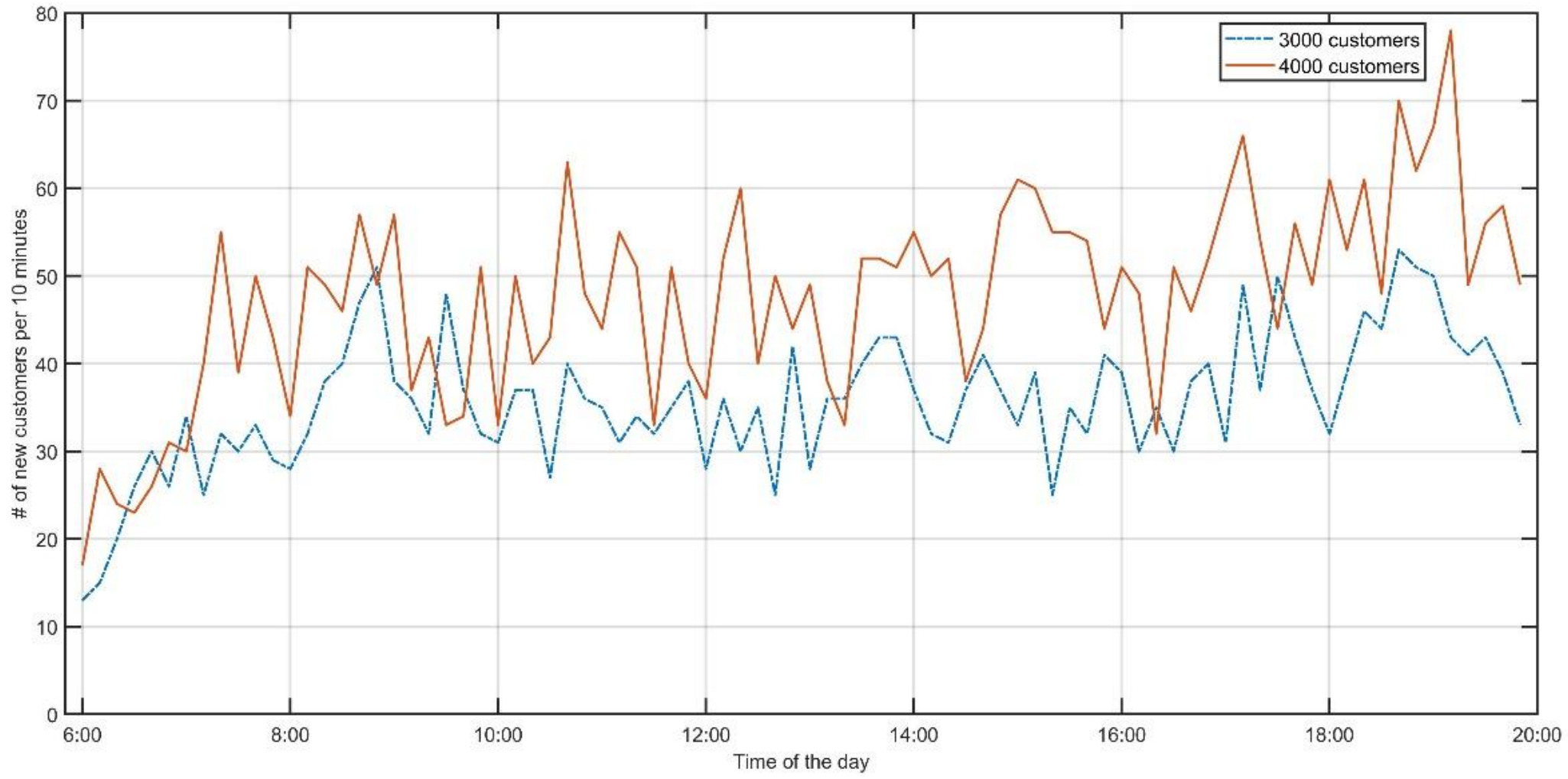

Figure 3. An example of customer arrival intensity per 10 minutes for two different demand scenarios.

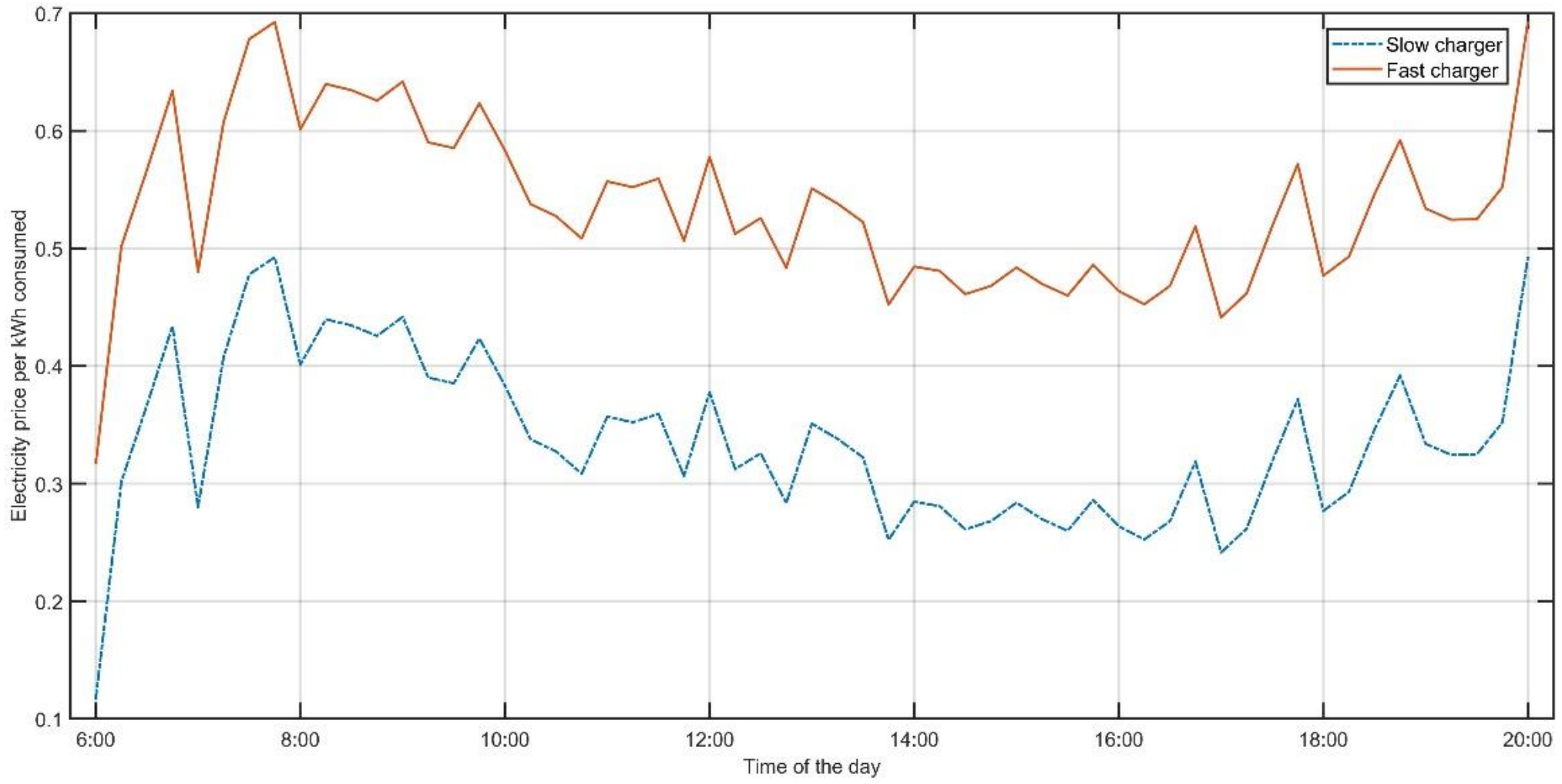

Figure 4. Electricity prices per kWh charged on different types of chargers during service hours. The difference between fast and slow chargers is the specific fee for using fast chargers.

## 4.2. Results

We test the performance of the CongestionAware charging policy on the test instances and compare it with the benchmark approaches. The implementation is based on Julia on a laptop with Intel(R) Core(TM) i7-11800H and 64 GB memory. The MILP models P1 to P3 are solved using Gurobi v12.0.3. We define a set of key performance indicators (KPIs) to evaluate the performance of different charging policies, including total profit (PF), revenue (TR), travel costs (TC), total charging costs (CC), charging costs during service hours (CC_D) and overnight (CC_N), total energy charged (ENG), energy charged during service hours (ENG_D) and overnight (ENG_N), customer service rate (SR), vehicle kilometer traveled (KMT), charging waiting time during service hours (WT), and charging time during service hours (CT). Total profit is calculated by subtracting travel and charging costs from total revenue. The charging costs include charging costs during service hours and for charging to full overnight. The acronyms and used measurement units are shown in Table 4.

Table 4. Acronyms used in the tables of computational studies.

| Acronym | Meaning | Unit |
|---|---|---|
| PF | Profit | 1000 USD |
| TR | Revenue | 1000 USD |
| TC | Travel costs | 1000 USD |
| CC | Total charging costs (CC_D+CC_N) | 1000 USD |
| CC_D | Charging costs during service hours | 1000 USD |
| CC_N | Charging costs overnight | 1000 USD |
| ENG | Total Energy charged (ENG_D+ENG_N) | kWh |
| ENG_D | Energy charged during service hours | kWh |
| ENG_N | Energy charged overnight | kWh |
| SR | Service rate | % |
| KMT | Vehicle-kilometer traveled | 1000 kilometers |
| WT | Charging waiting time during service hours | hour |
| CT | Charging time during service hours | hour |

### a. Model parameter estimation and base results

The model parameters for the day-ahead charging planning P1 are $C$, $\delta_h$, $\gamma$, and $\overline{W}_{hs}$, which need to be estimated based on historical driving and charging operations data. We estimate the average charging access cost $C$ as approximately \$2.7 based on an approximate average distance cost to the charging stations. $\delta_h$ is estimated by conducting a simulation using the Fastest charging policy. For $\overline{W}_{hs}$ and $\gamma$, we simulate the system for the two demand scenarios using the Fastest charging policy and obtain their respective averaged values. As mentioned previously, 10 independent datasets are used for the simulation to obtain these parameters for each demand scenario. Then we test the performance of the proposed method on the 5 test instances (days) and report the average results. Note that to further evaluate the robustness of the performance comparison, an instance-level paired difference comparison and statistical analysis are conducted using 30 independent test instances. The reader is referred to Appendix J for a more detailed description.

The results for two demand scenarios with 3000 and 4000 customers a day are shown in Table 5. Compared with the benchmark, the CongestionAware policy has the highest profit for both scenarios. Total profit increases 3.56%-14.90% for the c3000 scenario and 3.40%-21.68% for the c4000 scenario. Regarding service rate, the CongestionAware policy (86.5%) outperforms the benchmark by increasing 4.34%-18.33% and 3.87%-22.22% for the c3000 and c4000 scenarios, respectively. For the c3000 scenario, the rule-based benchmark charging policies result in much higher total charging waiting times (doubled or tripled) during service hours than the CongestionAware policy. As a result, the service rates of the benchmark policies are much lower due to lower vehicle availability. Comparing with the OptChg policy, using the CongestionAware policy results in similar waiting time but higher charging time (+7.85%) and energy charged (+5.01%). Regarding total charged energy, the CongestionAware policy charges 5.01%-17.47% more energy compared to the benchmark due to higher energy consumption to serve more customers. In terms of total charging costs, it is $1.77k, ranked as the second highest compared with the benchmark. Table 6 shows the charged energy and costs during service hours and for overnight over different charging policies. We can observe that only 25.63% of charged energy are realized during the service hours, while the rule-based benchmark charges more during the service hours (28.98%-37.35%). The OptChg policy exhibits the lowest percentage (21.73%) as it applies the same driving-need based approach as the CongestionAware policy for vehicle recharge. The corresponding charging costs for the CongestionAware policy is 51.21% of the total due to cheaper overnight charging costs. The charging cost ratio during the service hours for the benchmark ranges from 44.90%-63.61%.

For the c4000 high-demand scenario, similar results can be observed. The CongestionAware policy has much lower charging waiting time than the rule-based benchmark policies, resulting in a higher service rate (67.1%) and profit ($79.02k). Similar trend can be observed when comparing with the OptChg policy. As the service rate is higher for the CongestionAware policy, which means (generally) more KMT to serve more customers with more energy use and higher TC. The energy used per served customer is 2.80 (2.97) kWh for the c3000 (c4000) scenario for the CongestionAware policy, which is similar to the benchmark policies (2.78-2.89 kWh for c3000 scenario and 2.91-3.02 kWh for the c4000 scenario). Note that among the benchmark policies, the OptChg policy has the best performance as it applies the same online optimization-based approach to maximize vehicle availability. The second-best benchmark policy is the DynaThres because vehicles go to recharge earlier and avoid charging congestions in the peak charging demand period (to be explained below). We observe that the total energy charged is increased by +9.87% when customer demand is increased from 3000 to 4000. The energy charged increases 36.52% during the day to meet higher energy needs during the service hours. The energy charged during the night remains at a similar level (+0.68%). For the rule-based benchmark policy, the total charged energy increases less significantly (+3.35% to +4.13%) compared to the CongestionAware policy. The energy charged during the service hours remains similar compared to that of the C3000 scenario for the rule-based benchmark policies.

Table 5. Comparison of the KPIs for different charging policies.

| Scenario | Charging policy | PF | TR | TC | CC | ENG | SR | KMT | WT | CT |
|---|---|---|---|---|---|---|---|---|---|---|
| c3000 (3000 requests) | Nearest | 63.27 | 77.91 | 13.11 | 1.53 | 6186 | 73.1% | 24.7 | 82.9 | 58.8 |
| | Fastest | 63.83 | 78.66 | 13.25 | 1.58 | 6250 | 73.8% | 25.0 | 79.1 | 58.1 |
| | MinChgOpT | 64.63 | 79.68 | 13.42 | 1.63 | 6330 | 74.8% | 25.3 | 64.4 | 59.3 |
| | DynaThreshold | 66.79 | 82.98 | 14.32 | 1.86 | 6755 | 78.0% | 27.0 | 96.2 | 77.9 |
| | OptChg | 70.20 | 86.44 | 14.67 | 1.57 | 6920 | 82.9% | 27.7 | 26.8 | 38.2 |
| | **CongestionAware** | **72.70** | **89.88** | **15.41** | **1.77** | **7267** | **86.5%** | **29.1** | **26.8** | **41.2** |
| | Nearest | 64.94 | 80.05 | 13.55 | 1.55 | 6393 | 54.9% | 25.6 | 113.8 | 62.9 |
| | Fastest | 65.97 | 81.35 | 13.78 | 1.59 | 6502 | 55.6% | 26.0 | 104.3 | 61.8 |
| c4000 (4000 requests) | MinChgOpT | 66.64 | 82.17 | 13.91 | 1.63 | 6560 | 56.2% | 26.2 | 96.3 | 62.9 |
| | DynaThreshold | 70.01 | 86.79 | 14.91 | 1.87 | 7034 | 58.3% | 28.1 | 97.4 | 80.4 |
| | OptChg | 76.42 | 94.63 | 16.28 | 1.93 | 7679 | 64.6% | 30.7 | 45.6 | 57.4 |
| | **CongestionAware** | **79.02** | **98.06** | **16.93** | **2.11** | **7984** | **67.1%** | **31.9** | **43.0** | **62.8** |

Table 6. Charging costs and the amount of charged energy during the service hours and overnight.

| Scenario | Charging policy | ENG | ENG_D | | ENG_N | | CC_D | | CC_N | |
|---|---|---|---|---|---|---|---|---|---|---|
| | | KWh | KWh | % | KWh | % | 1000 USD | % | 1000 USD | % |
| c3000 (3000 requests) | Nearest | 6186 | 1793 | 28.98 | 4393 | 71.02 | 0.83 | 54.18 | 0.70 | 45.82 |
| | Fastest | 6250 | 1874 | 29.98 | 4376 | 70.02 | 0.88 | 55.58 | 0.70 | 44.42 |
| | MinChgOpT | 6330 | 1982 | 31.32 | 4347 | 68.68 | 0.93 | 57.28 | 0.70 | 42.72 |
| | DynaThreshold | 6755 | 2523 | 37.35 | 4232 | 62.65 | 1.18 | 63.61 | 0.68 | 36.39 |
| | OptChg | 6920 | 1504 | 21.73 | 5416 | 78.27 | 0.71 | 44.90 | 0.87 | 55.10 |
| | **CongestionAware** | **7267** | **1862** | **25.63** | **5405** | **74.37** | **0.91** | **51.21** | **0.86** | **48.79** |
| c4000 (4000 requests) | Nearest | 6393 | 1816 | 28.41 | 4577 | 71.59 | 0.82 | 52.84 | 0.73 | 47.16 |
| | Fastest | 6502 | 1880 | 28.91 | 4622 | 71.09 | 0.86 | 53.63 | 0.74 | 46.37 |
| | MinChgOpT | 6560 | 1947 | 29.68 | 4613 | 70.32 | 0.89 | 54.62 | 0.74 | 45.38 |
| | DynaThreshold | 7034 | 2461 | 34.99 | 4572 | 65.01 | 1.14 | 60.82 | 0.73 | 39.18 |
| | OptChg | 7679 | 2202 | 28.68 | 5477 | 71.32 | 1.05 | 54.56 | 0.88 | 45.44 |
| | **CongestionAware** | **7984** | **2542** | **31.84** | **5442** | **68.16** | **1.24** | **58.78** | **0.87** | **41.22** |

Figures 5 and 6 analyze the number of vehicles charging (subfigures (b)) and of waiting (subfigures (c)) for the day using different charging policies. Subfigure (a) shows the average waiting time at fast chargers that vehicles would have experienced using the Fastest charging policy. For the c3000 scenario, Subfigure (a) in Figure 5 shows the charging waiting times at fast chargers increase significantly from the $8^{th}$ hour from the beginning of service (14:00) to the $12^{th}$ hour (18:00) when using the Fastest charging policy. As a result, the day-ahead charging plan model devises a charging plan to avoid vehicle charging and waiting during peak charging demand hours. Subfigure (b) shows the CongestionAware charging policy activates vehicle charging operations earlier to reduce charging congestion compared to the benchmark. During peak charging demand hours (from $4^{th}$ to $8^{th}$ hours), the number of vehicles charging is around 15 per 30 minutes, lower than the benchmark with around 1 or 2 more vehicles. As this policy anticipates vehicles' energy (driving) need to the end of service hours, less energy is charged in the late evening, resulting in a higher service rate and fewer vehicles charging and waiting. Subfigure (c) shows that using the first three benchmark policies results in high charging congestion (more vehicles waiting for charging) during the evening peak (from the $10^{th}$ hour to the $12^{th}$ hour) due to the increasing number of vehicles' SoCs falling below the threshold of 20% battery capacity. The charging demand decreases after the $12^{th}$ hour in the evening. For the DynaThreshold policy, the number of vehicles waiting increases gradually from the $6^{th}$ hour until the $11^{th}$ hours (17:00). This is because the hourly charging thresholds for the second half of the day are lower (around 50% in the morning and 30% in the afternoon on average), resulting in more vehicles charging and waiting at the chargers due to charging longer time in the afternoon and evening, and insufficient number of fast chargers in the system. However, this DynaThreshold policy performs better than the other three benchmark policies as it applies a smarter partial recharge policy. For the OptChg policy, its charging demand (number of vehicles charging) starts from around the $6^{th}$ hour based on its day-ahead charging plan (See Appendix F for the corresponding MILP formulation and detailed analysis).

For the c4000 scenario, the charging congestion starts earlier (mainly located between the $7^{th}$ and $10^{th}$ hours instead of between the $8^{th}$ and $12^{th}$ hours; see subfigures (a) of Figure 6 and Figure 5 due to higher customer arrival intensity for the c4000 scenario. Consequently, the obtained day-ahead charging plan reacts in response to the charging waiting time signals, resulting in significantly higher charging operations in the early hours of the day. As the CongestionAware policy applies a smarter partial recharge strategy (Eq. (29)) to delay vehicles' charging operations when charging waiting time on a queuing charger exceeds a maximum threshold (e.g. 30 minutes) and minimize the charging operational times for online vehicle-to-charger assignment, it allows significantly increasing vehicles' availability to serve customers.

More detailed results related to the standard deviations of the KPIs, histograms of the SoC of vehicles at the end of service hours, and a representative vehicle's energy profile and its state changes using different charging policies during the service hours can be found in Appendix G.

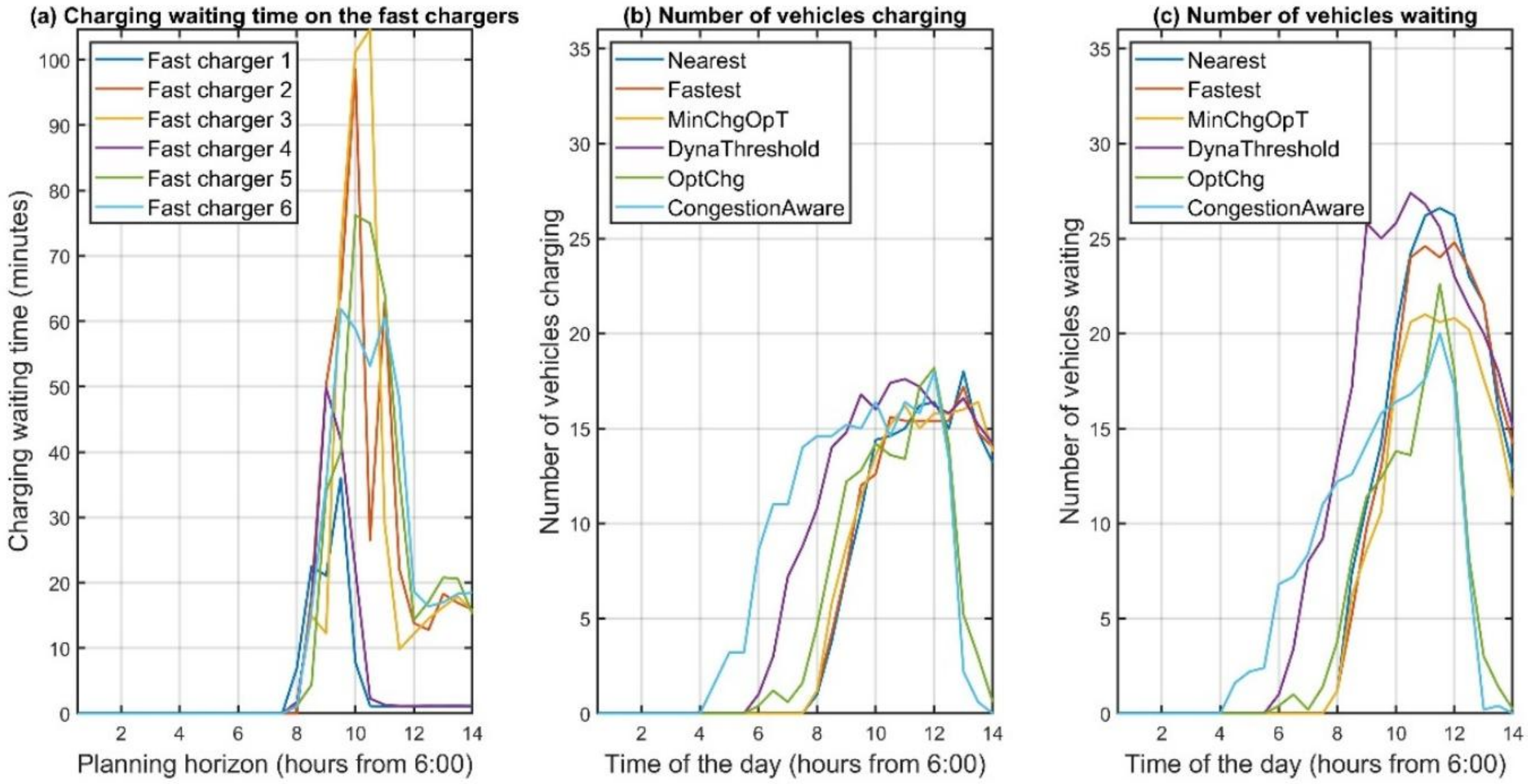

Figure 5. Comparison of charging waiting times at fast chargers, number of vehicles charging, and number of vehicles waiting during the day (service hours) (c3000 scenario).

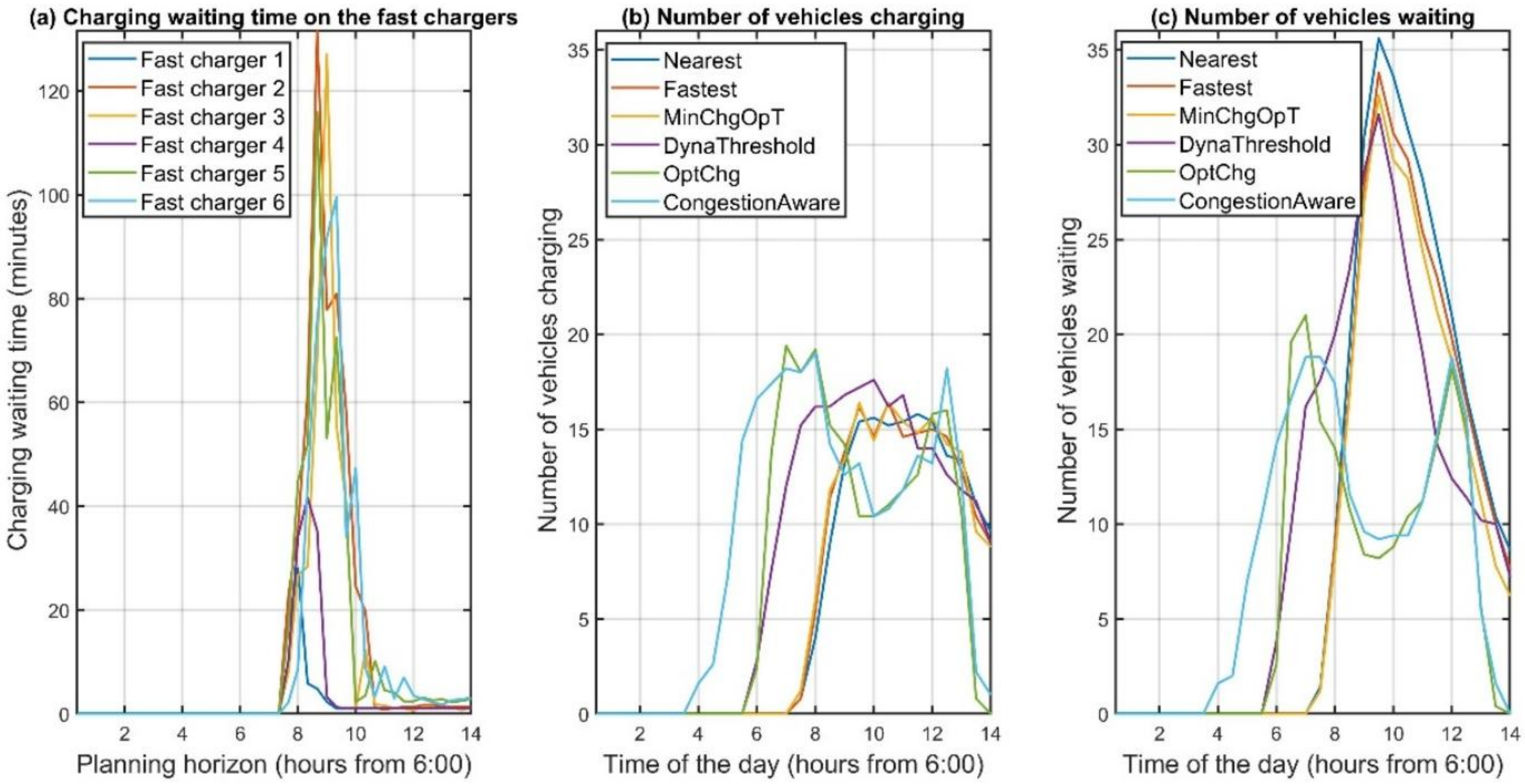

Figure 6. Comparison of charging waiting times at fast chargers, number of vehicles charging, and number of vehicles waiting during the day (service hours) (c4000 scenario).

b. Analysis of realized charging sessions and occupancy state of chargers

Figure 7 and Figure 8 compare the distributions of charging session durations and the amount of charged energy per charging session using different charging policies where Subfigure (a) is related to the c3000 scenario, while Subfigure(b) to the c4000 scenario. The CongestionAware policy has a significantly lower charging duration and charged energy per charging session compared with the benchmark. For the c3000 scenario, the median charging session duration and median charged energy for the CongestionAware policy is 14.9 minutes (75$^{th}$ percentile = 19.1 minutes) and 12.1 kWh (75$^{th}$ percentile = 15.0 kWh), respectively. For the c4000 scenario, the median duration of charging sessions and median charged amount of energy is 15.4 minutes (75$^{th}$ percentile = 27.9 minutes) and 11.9 kWh (75$^{th}$ percentile = 17.5 kWh), respectively. Compared with the CongestionAware policy, the median charging durations and median charged energy per session of the rule-based benchmark policies are much higher (around 50 minutes and 35-42kWh for both scenarios). The OptChg policy has slightly higher median charging durations and charged energy per session than the CongestionAware policy, but also exhibits a wider distribution for both metrics. As we require a minimum charging duration of 10 minutes for vehicle recharging, we can observe the minimum charging times are above it for different charging policies (Figure 7). For the charged energy, we observe two or three vehicles charged around 2 or 3 kWh only for some benchmark policies. This is because these vehicles started charging near the end of the day on slow chargers for about 10–15 minutes.

Figure 9 compares the distribution of charging waiting times per charging session using different charging policies. A vehicle's charging waiting times for a charging session account for its waiting times on the queue of the used charger and waiting times for previously unsuccessful attempts on other chargers (charger-chasing B, see Appendix C). For the c3000 scenario, the waiting time median of the CongestionAware policy is 9.5 minutes, while it ranges from 24.7 to 40.3 minutes for the rule-based benchmark. However, the 75th percentile of waiting times for the CongestionAware policy is much lower (24 minutes) compared to the rule-based benchmark (54-101 minutes) due to longer charging times for more charged energy. The OptChg policy achieves the lowest waiting times among all benchmark policies, owing to its partial recharging strategy based on the vehicle's anticipated driving need until the end of the service hours.

Figure 10 shows the number of vehicles on fast and slow chargers over different charging policies, where Subfigures (a) (fast chargers) and (b)(slow chargers) are related to using the benchmark policies, while Subfigures (c) (fast chargers) and (d) (slow chargers) are related to using the CongestionAware policy. For the CongestionAware policy, all fast chargers (6 in total) are almost fully occupied from the 8th to the 12th hours. The CongestionAware policy has a higher occupancy rate on fast chargers compared to the benchmark, except the DynaThreshold policy, between the 6$^{th}$ and 9$^{th}$ hours. For slow chargers, using the rule-based benchmark policies results in a higher charger occupancy on slow chargers (11kW) after the 7$^{th}$ hour. This is because vehicles move away from fast chargers after waiting for a maximum waiting time (30 minutes) and go to another fast/slow charger with the least waiting time. Note that we might apply different waiting policies at chargers. We test the effect of applying different waiting policies (i.e., no limit waiting (naïve), charger-chasing A using fast chargers only or charger-chasing B using fast or slow chargers). The results show the above charger-chasing B policy has better performance compared to the other two waiting policies (see Appendix C). The charger occupancy on slow chargers is very low (0 most of the time or <3 between the 9$^{th}$ and 12$^{th}$ hours) when using the CongestionAware policy. This is because the online vehicle-to-charge assignment model P3 minimizes the total charging operation times (including charging access time, waiting time, and charging time) under the constraints that the SoCs of vehicles after recharge need to be no less than their target energy levels (Eq. (22)) and satisfy a minimum charged energy requirement. Given their lower charging power (11 kW), slow chargers generally require longer charging durations to deliver the same amount of energy as fast chargers. Moreover, the imposed minimum charged-energy requirement and the maximum charging duration per session may restrict the feasibility of assigning vehicles to slow chargers.

Comparing the CongestionAware policy with the OptChg benchmark reveals distinct charger occupancy patterns. The CongestionAware policy starts using fast chargers earlier and progressively increases their occupancy, whereas the OptChg policy exhibits a later and sharper increase in occupancy. For slow chargers, the OptChg policy results in substantially higher occupancy. These differences can be mainly attributed to their day-ahead charging planning models. Unlike CongestionAware, OptChg does not explicitly account for charger capacity constraints or distinguish between charger types in its day-ahead planning model, resulting in more concentrated charging demand during certain periods (see the comparison of day-ahead charging plans in Figure F1, Appendix F). In contrast, CongestionAware better distributes charging demand over time, enabling more efficient charging operations under limited charging capacity.

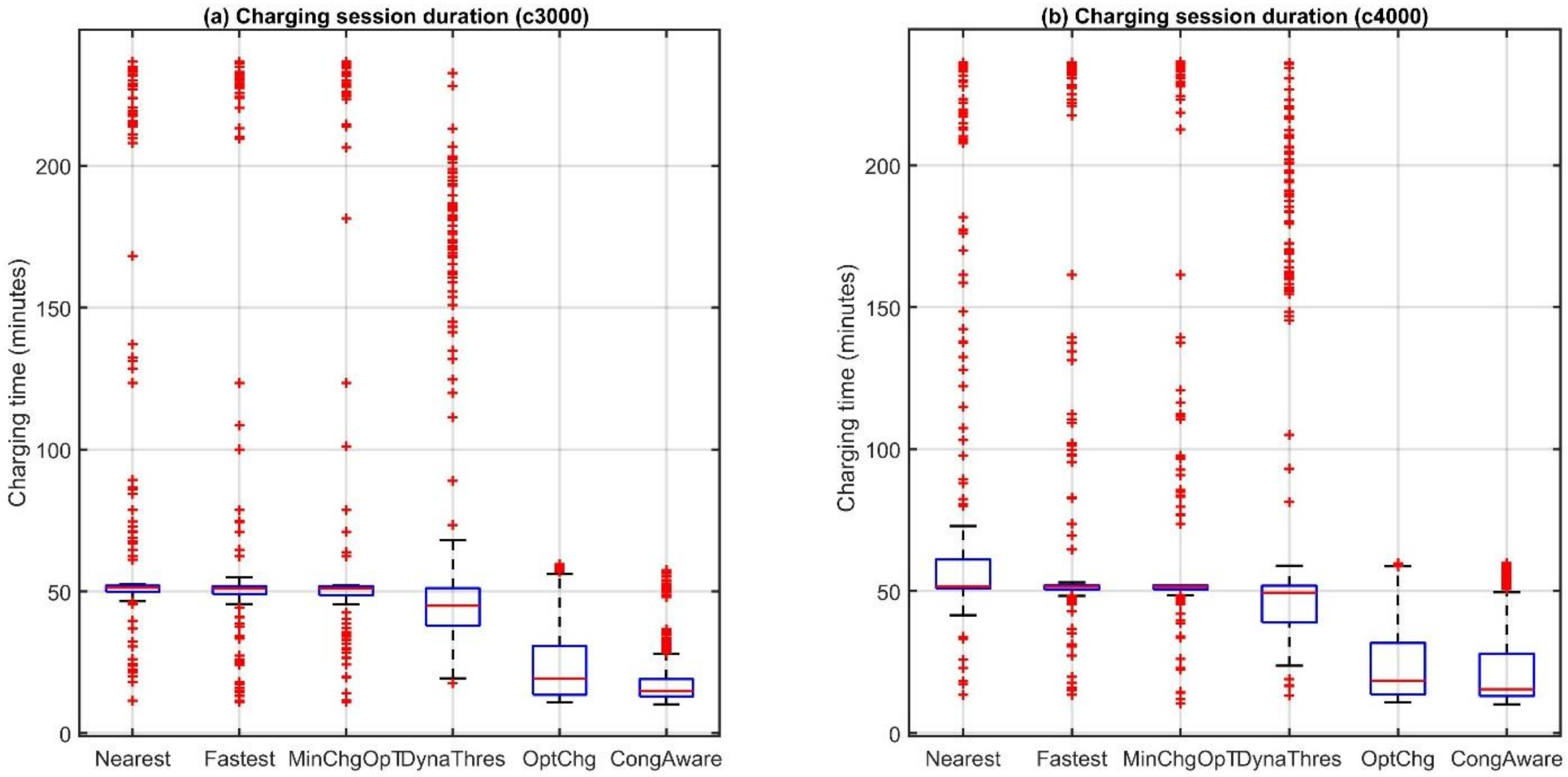

Figure 7. Distribution of the charging session durations for the c3000 scenario (on the left) and the c4000 scenario (on the right).

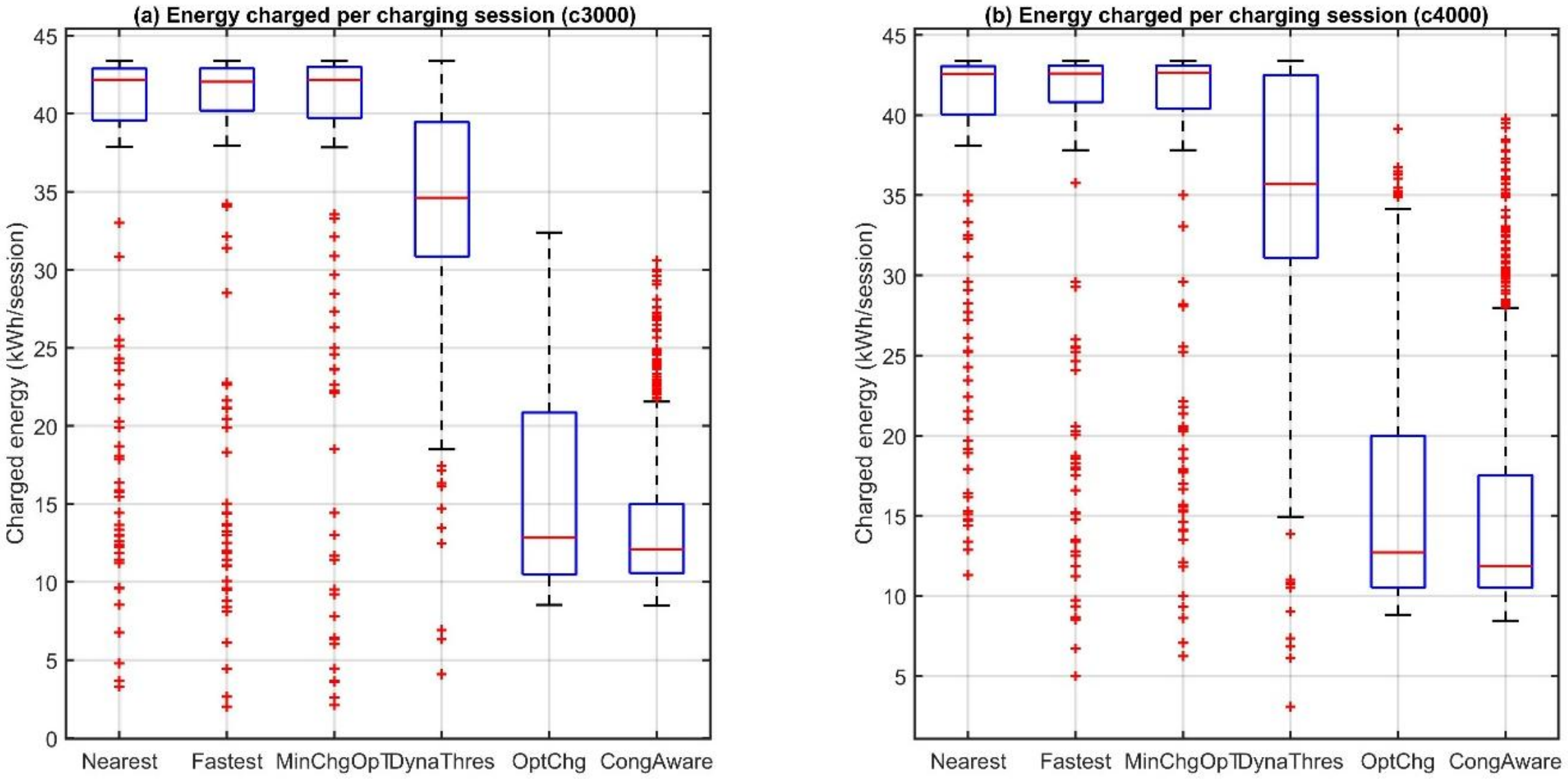

Figure 8. Distribution of the amount of charged energy for the c3000 scenario (on the left) and the c4000 scenario (on the right).

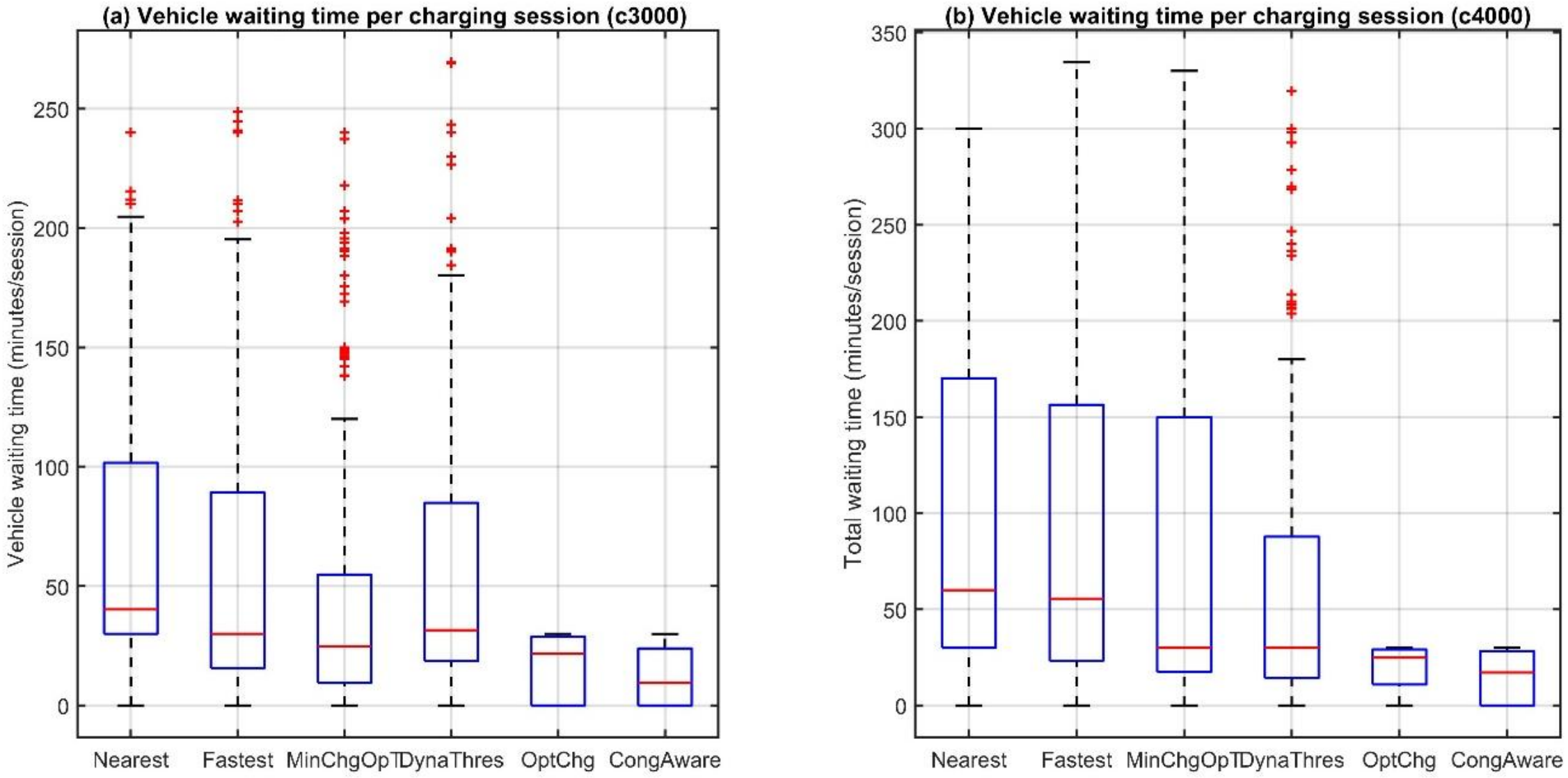

Figure 9. Distribution of charging waiting times per charging session for the c3000 scenario (on the left) and the c4000 scenario (on the right).

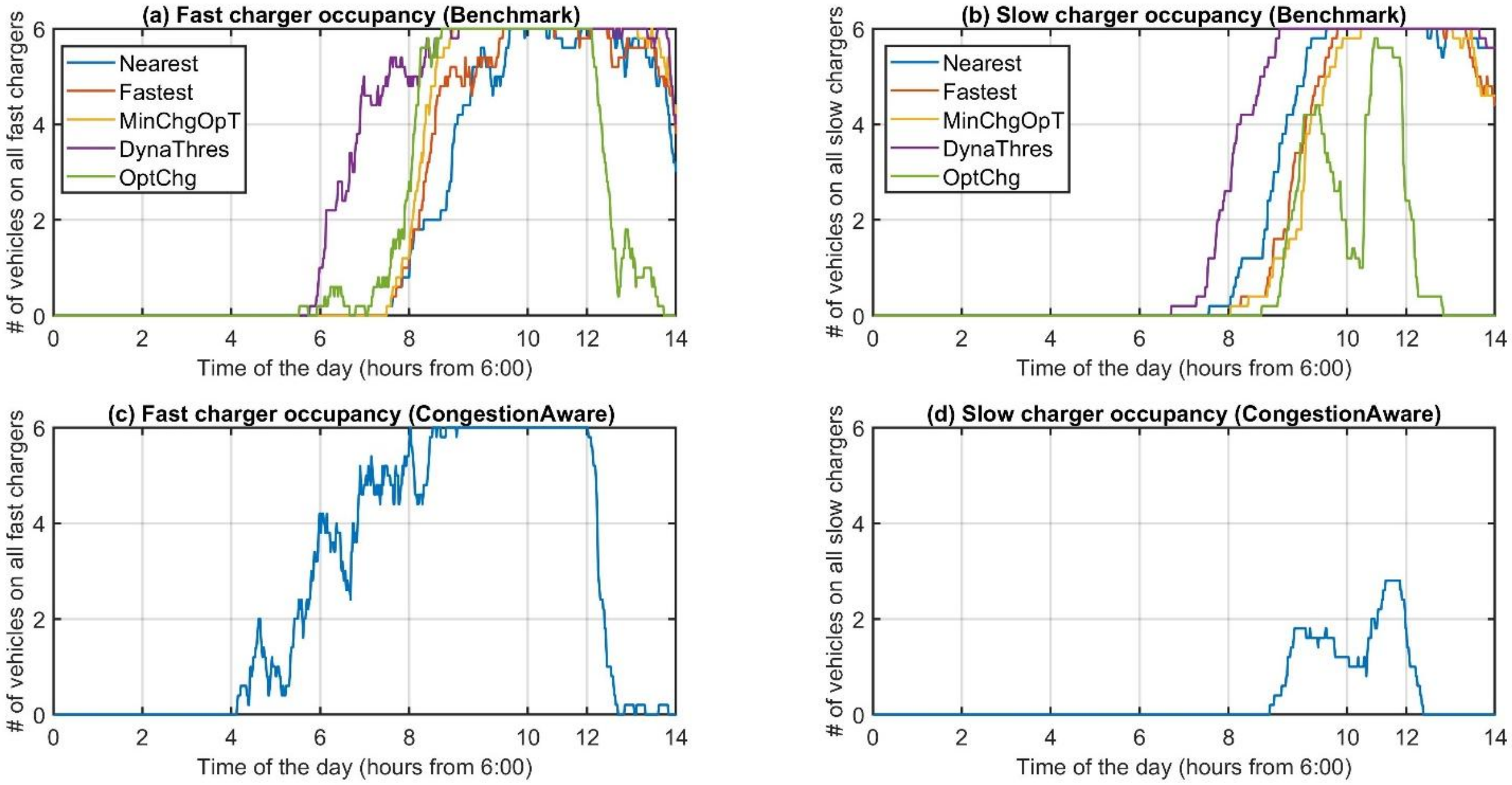

Figure 10. Comparison of the average number of vehicles on fast chargers (on the left) and slow chargers (on the right) using different charging policies for the c3000 scenario.

c. Effects of the day-ahead charging plan and vehicles' energy need anticipation

For the CongestionAware policy, we further investigate the benefits of the day-ahead charging plan by comparing system performance with and without it under the c3000 and c4000 scenarios. Without P1, charging is triggered by the SoC threshold, $\theta$, while the congestion-aware assignment (P3) and partial-recharge anticipation are retained. The results are presented in Table 7. We observe that incorporating the day-ahead charging plan increases the service rate by 5.0% and by 9.8% for the c3000 and c4000 scenarios, respectively. Total profit follows a similar trend, increasing by 4.1% for c3000 and 9.5% for c4000. The results demonstrate the effectiveness of integrating the day-ahead charging plan in improving the system performance, in particular for higher customer demand. Similarly, for the energy-need anticipation strategy, which determines the target energy level (Eq. (29)), we compare the system performance with and without using it, i.e., vehicles charge to $E_v^{max}$ (80%) from their current battery levels. The results are shown in Table 8. We can observe that using this strategy allows for increasing customers' service rate (+4.5% to +8.1%) and profits (+3.7% to +7.8%). When not applying this anticipative strategy, the CongestionAware policy still significantly outperforms the rule-based benchmark in terms of total profit and service rate (see Table 5).

Table 7. Comparison of the CongestionAware policy with and without day-ahead charging plan.

| Scenario | Day-ahead plan | PF | TR | TC | CC | ENG | SR | KMT | WT | CT |
|---|---|---|---|---|---|---|---|---|---|---|
| c3000 | **Yes** | **72.70** | **89.88** | **15.41** | **1.77** | **7267** | **86.5%** | **29.1** | **26.8** | **41.2** |
| | No | 69.87 | 85.92 | 14.53 | 1.52 | 6853 | 82.4% | 27.4 | 26.3 | 32.9 |
| c4000 | **Yes** | **79.02** | **98.06** | **16.93** | **2.11** | **7984** | **67.1%** | **31.9** | **43.0** | **62.8** |
| | No | 72.18 | 88.98 | 15.13 | 1.67 | 7135 | 61.1% | 28.5 | 28.1 | 37.4 |

Table 8. Comparison of the CongestionAware policy with and without anticipating energy needs.

| Scenario | Anticipate | PF | TR | TC | CC | ENG | SR | KMT | WT | CT |
|---|---|---|---|---|---|---|---|---|---|---|
| c3000 | **Yes** | **72.70** | **89.88** | **15.41** | **1.77** | **7267** | **86.5%** | **29.1** | **26.8** | **41.2** |
| | No | 70.10 | 86.72 | 14.72 | 1.89 | 6943 | 82.8% | 27.8 | 23.1 | 48.8 |
| c4000 | **Yes** | **79.02** | **98.06** | **16.93** | **2.11** | **7984** | **67.1%** | **31.9** | **43.0** | **62.8** |
| | No | 73.28 | 91.08 | 15.68 | 2.12 | 7397 | 62.1% | 29.6 | 37.0 | 60.1 |

### 4.3. Sensitivity analysis

In this section, we investigate the impact of different model parameters on the performance of the CongestionAware policy. These parameters include the maximum charging waiting times at a queuing charger, battery capacity, and number of fast and slow chargers in the system.

#### a. Impact of the maximum charging waiting time of vehicles at a queuing charger

The maximum charging waiting time at a queuing charger affects vehicles' availability to serve customers as the operator can delay vehicles' charging operations (idle for serving customers) when the estimated queuing time on the assigned charger (has minimum total charging operational time among all chargers) exceeds the maximum threshold. Table 9 compares the performance using 20, 30, 40 minutes, and no maximum waiting time limit for different demand scenarios. We can observe that when there is no maximum charging waiting time, the profit is the lowest compared with the other cases for high-demand scenarios (total charging waiting time is 146.1 hours for no-limit waiting compared with the other cases (less than 59.9 hours)). This is because the number of fast chargers is very insufficient in the system, and vehicles need to wait a long time at fast chargers if there are no maximum charging waiting time limits. On the other hand, when limiting vehicles' maximum waiting time at a queuing charger, vehicles can be idle to serve customers. Note that when increasing this maximum threshold, the total waiting time increases accordingly. However, the total customer service rate and profit increase if the maximum waiting time threshold is limited, as is the case for a threshold of no more than 40 minutes. In practice, the operator could learn/adjust this threshold based on its day-to-day realized vehicle queuing patterns and charging demands.

Table 9. Impact of maximum waiting time at chargers for the CongestionAware policy.

| Scenario | Max. waiting time* | PF | TR | TC | CC | ENG | SR | KMT | WT | CT |
|---|---|---|---|---|---|---|---|---|---|---|
| c3000 | 20 | 72.12 | 89.16 | 15.29 | 1.75 | 7210 | 85.9% | 28.8 | 17.8 | 37.1 |
| | 30 | 72.70 | 89.88 | 15.41 | 1.77 | 7267 | 86.5% | 29.1 | 26.8 | 41.2 |
| | **40** | **73.21** | **90.46** | **15.47** | **1.78** | **7299** | **87.1%** | **29.2** | **33.9** | **43.0** |
| | no limit | 73.27 | 90.56 | 15.50 | 1.79 | 7313 | 87.1% | 29.3 | 36.9 | 43.3 |
| c4000 | 20 | 78.24 | 97.02 | 16.71 | 2.07 | 7884 | 66.6% | 31.5 | 28.4 | 56.0 |
| | 30 | 79.02 | 98.06 | 16.93 | 2.11 | 7984 | 67.1% | 31.9 | 43.0 | 62.8 |
| | **40** | **79.14** | **98.30** | **17.01** | **2.16** | **8021** | **67.3%** | **32.1** | **59.9** | **68.5** |
| | no limit | 75.71 | 94.29 | 16.37 | 2.21 | 7721 | 64.2% | 30.9 | 146.1 | 79.3 |

*: in minutes.

#### b. Impact of battery capacity

To test the impact of battery capacity, we consider three battery sizes, i.e., 62, 72, and 82 kWh. We consider only the high demand scenario of 4000 customers with identical parameter settings (results for the c3000 scenario would draw similar conclusions). Table 10 compares the performance of different charging policies. As expected, increasing the battery size could significantly increase the customer service rate and the system's profit (i.e., the customer service rate increases around 10% for different charging policies if the battery size is increased from 62 kWh to 82 kWh). Compared with the benchmark policies, the CongestionAware policy has higher profits and service rates for all battery sizes. When comparing with the OptChg policy, the profit and customer service rate increase ranging from 3.40% to 7.73% and from 3.87% to 8.72%, respectively. When comparing with the best rule-based benchmark, the profit and customer service rate increase from 12.87% to 14.67% and from 14.29% to 16.61%, respectively. The results demonstrate that by increasing the battery size (given the same charging station capacity limit), the CongestionAware policy further outperforms the benchmark. However, the relative improvement over the OptChg policy peaks at 72 kWh and is not monotonic. This is because when battery size (capacity) increases, vehicles' charging times and waiting times become longer, which could increasingly harm vehicles' availability when charging to 80% for each charging operation. However, the CongestionAware policy devises a smarter charging plan that anticipates vehicles' charging waiting times and energy needs and adapts it (reactive model) to minimize total system costs, resulting in more effective utilization of congested charging facilities.

Table 10. Impact of battery capacity on the KPI using different charging policies (c4000 scenario).

| Battery | Charging policy | PF | TR | TC | CC | ENG | SR | KMT | WT | CT |
|---|---|---|---|---|---|---|---|---|---|---|
| 62 kWh | Nearest | 64.94 | 80.05 | 13.55 | 1.55 | 6393 | 54.9% | 25.6 | 113.8 | 62.9 |
| | Fastest | 65.97 | 81.35 | 13.78 | 1.59 | 6502 | 55.6% | 26.0 | 104.3 | 61.8 |
| | MinChgOpT | 66.64 | 82.17 | 13.91 | 1.63 | 6560 | 56.2% | 26.2 | 96.3 | 62.9 |
| | DynaThreshold | 70.01 | 86.79 | 14.91 | 1.87 | 7034 | 58.3% | 28.1 | 97.4 | 80.4 |
| | OptChg | 76.42 | 94.63 | 16.28 | 1.93 | 7679 | 64.6% | 30.7 | 45.6 | 57.4 |
| | **CongestionAware** | **79.02** | **98.06** | **16.93** | **2.11** | **7984** | **67.1%** | **31.9** | **43.0** | **62.8** |
| 72 kWh | Nearest | 69.67 | 85.83 | 14.56 | 1.60 | 6866 | 59.9% | 27.5 | 106.0 | 53.7 |
| | Fastest | 70.22 | 86.56 | 14.69 | 1.65 | 6930 | 60.3% | 27.7 | 101.4 | 54.6 |
| | MinChgOpT | 70.70 | 87.14 | 14.77 | 1.68 | 6967 | 60.6% | 27.9 | 90.8 | 55.4 |
| | DynaThreshold | 73.41 | 90.86 | 15.56 | 1.88 | 7339 | 62.0% | 29.4 | 87.6 | 73.0 |
| | OptChg | 78.10 | 96.26 | 16.42 | 1.74 | 7746 | 66.5% | 31.0 | 27.2 | 38.6 |
| | **CongestionAware** | **84.14** | **104.33** | **18.05** | **2.13** | **8515** | **72.3%** | **34.1** | **34.2** | **54.7** |
| 82 kWh | Nearest | 75.97 | 93.40 | 15.85 | 1.58 | 7478 | 65.9% | 29.9 | 79.1 | 39.9 |
| | Fastest | 76.42 | 94.05 | 15.99 | 1.64 | 7542 | 66.2% | 30.2 | 73.2 | 41.3 |
| | MinChgOpT | 76.83 | 94.55 | 16.05 | 1.67 | 7573 | 66.5% | 30.3 | 70.0 | 42.2 |
| | DynaThreshold | 76.70 | 94.67 | 16.19 | 1.78 | 7638 | 65.9% | 30.6 | 86.8 | 54.1 |
| | OptChg | 84.59 | 104.35 | 17.91 | 1.85 | 8448 | 72.5% | 33.8 | 30.5 | 42.7 |
| | **CongestionAware** | **88.10** | **108.89** | **18.76** | **2.04** | **8847** | **76.0%** | **35.4** | **22.7** | **37.5** |

c. Impact of the number of fast and slow chargers

We further analyze the sensitivity of increasing the number of fast and slow chargers in the system for the c3000 (Table 11) and c4000 (Table 12) scenarios. The number of fast and slow chargers is increased from 12 chargers (6 fast and 6 slow) to 20 (10 fast and 10 slow) at the same charging stations (i.e. from 3 fast/slow chargers to 5 fast/slow chargers per charging station). The battery size is 62 kWh, and the other parameters are identical as described in Section 4.1. For the **c3000 scenario**, increasing the number of chargers in the system does reduce the charging congestion for the benchmark policies. The profit, service rate, and charged amount of energy increase along with more fast and slow chargers in the system. The service rate increases around 2.1% to 7.4% if the number of chargers is increased from 12 to 20. For the CongestionAware policy, the benefit is less significant in terms of customer service rate and total system profit, although the total charging waiting time decreases accordingly when the number of chargers increases.

For the **c4000 scenario**, using the rule-based benchmark policies results in a significant reduction in total (charging) waiting time (see WT column in Table 12) when increasing the number of chargers from 12 to 16 (i.e., from -17.4% to -41.3%). Adding more chargers from 12 to 20 chargers can further effectively reduce total waiting time (i.e., from -34.7% to -56.3%). The profit (+3.1% to +5.9%) and service rate (+2.6% to +6.3%) increase accordingly due to the total waiting time reduction when the number of chargers increases from 12 to 20. Increasing the number of chargers improves the performance of the CongestionAware policy, although its relative improvement in profit is smaller than that of the OptChg policy. Consequently, the relative advantage of CongestionAware generally diminishes as charging infrastructure becomes less constrained. The profit and service rate are increased by 6.5% (from 79.02k to 84.16k) and 8.0% (from 67.1% to 72.5%), respectively, if the number of chargers is changed from 12 to 20. Interestingly, we observe that total charging waiting time is reduced slightly from 43.0 hours (12 chargers) to 40.9 hours (20 chargers), but the total charged amount of energy is increased significantly from 7984 kWh to 8576 kWh (+7.4%). In practice, the operator can further invest in their charging infrastructure to increase the system's profitability when charging capacity is scarce.

Table 11. System performance for different numbers of fast and slow chargers (c3000 scenario).

| # of chargers | Charging policy | PF | TR | TC | CC | ENG | SR | KMT | WT | CT |
|---|---|---|---|---|---|---|---|---|---|---|
| 12 | Nearest | 63.27 | 77.91 | 13.11 | 1.53 | 6186 | 73.1% | 24.7 | 82.9 | 58.8 |
| | Fastest | 63.83 | 78.66 | 13.25 | 1.58 | 6250 | 73.8% | 25.0 | 79.1 | 58.1 |
| | MinChgOpT | 64.63 | 79.68 | 13.42 | 1.63 | 6330 | 74.8% | 25.3 | 64.4 | 59.3 |
| | DynaThreshold | 66.79 | 82.98 | 14.32 | 1.86 | 6755 | 78.0% | 27.0 | 96.2 | 77.9 |
| | OptChg | **70.20** | **86.44** | **14.67** | **1.57** | **6920** | **82.9%** | **27.7** | **26.8** | **38.2** |
| | **CongestionAware** | **72.70** | **89.88** | **15.41** | **1.77** | **7267** | **86.5%** | **29.1** | **26.8** | **41.2** |
| 16 | Nearest | 64.25 | 79.18 | 13.32 | 1.62 | 6285 | 74.2% | 25.1 | 62.8 | 70.2 |
| | Fastest | 65.08 | 80.29 | 13.52 | 1.69 | 6376 | 75.4% | 25.5 | 45.6 | 69.6 |
| | MinChgOpT | 66.41 | 81.98 | 13.81 | 1.76 | 6512 | 77.0% | 26.0 | 36.2 | 69.0 |
| | DynaThreshold | 69.68 | 86.57 | 14.86 | 2.03 | 7010 | 81.9% | 28.0 | 68.1 | 91.5 |
| | OptChg | **72.47** | **89.41** | **15.23** | **1.71** | **7185** | **85.9%** | **28.7** | **31.0** | **38.3** |
| | **CongestionAware** | **73.35** | **90.54** | **15.43** | **1.77** | **7277** | **87.4%** | **29.1** | **16.8** | **37.0** |
| 20 | Nearest | 64.54 | 79.57 | 13.38 | 1.66 | 6310 | 74.6% | 25.2 | 54.8 | 76.2 |
| | Fastest | 65.85 | 81.25 | 13.68 | 1.72 | 6454 | 76.4% | 25.8 | 35.4 | 74.5 |
| | MinChgOpT | 67.95 | 83.89 | 14.12 | 1.82 | 6660 | 79.0% | 26.6 | 25.1 | 65.8 |
| | DynaThreshold | 71.07 | 88.21 | 15.10 | 2.05 | 7121 | 83.8% | 28.5 | 35.1 | 93.4 |
| | OptChg | 73.21 | 90.31 | 15.37 | 1.73 | 7251 | 87.2% | 29.0 | 19.0 | 36.5 |
| | **CongestionAware** | **73.30** | **90.44** | **15.39** | **1.75** | **7258** | **87.5%** | **29.0** | **3.7** | **36.5** |

Table 12. System performance for different numbers of fast and slow chargers (c4000 scenario).

| # of chargers | Charging policy | PF | TR | TC | CC | ENG | SR | KMT | WT | CT |
|---|---|---|---|---|---|---|---|---|---|---|
| 12 | Nearest | 64.94 | 80.05 | 13.55 | 1.55 | 6393 | 54.9% | 25.6 | 113.8 | 62.9 |
| | Fastest | 65.97 | 81.35 | 13.78 | 1.59 | 6502 | 55.6% | 26.0 | 104.3 | 61.8 |
| | MinChgOpT | 66.64 | 82.17 | 13.91 | 1.63 | 6560 | 56.2% | 26.2 | 96.3 | 62.9 |
| | DynaThreshold | 70.01 | 86.79 | 14.91 | 1.87 | 7034 | 58.3% | 28.1 | 97.4 | 80.4 |
| | OptChg | 76.42 | 94.63 | 16.28 | 1.93 | 7679 | 64.6% | 30.7 | 45.6 | 57.4 |
| | **CongestionAware** | **79.02** | **98.06** | **16.93** | **2.11** | **7984** | **67.1%** | **31.9** | **43.0** | **62.8** |
| 16 | Nearest | 65.91 | 81.23 | 13.74 | 1.58 | 6481 | 55.6% | 25.9 | 94.0 | 71.8 |
| | Fastest | 67.70 | 83.49 | 14.15 | 1.64 | 6673 | 56.9% | 26.7 | 71.1 | 70.9 |
| | MinChgOpT | 68.66 | 84.67 | 14.32 | 1.69 | 6754 | 57.7% | 27.0 | 64.0 | 71.3 |
| | DynaThreshold | 72.83 | 90.24 | 15.45 | 1.95 | 7288 | 60.8% | 29.2 | 57.2 | 85.9 |
| | OptChg | 81.28 | 100.95 | 17.45 | 2.22 | 8230 | 69.5% | 32.9 | 55.2 | 65.4 |
| | **CongestionAware** | **81.85** | **101.71** | **17.59** | **2.27** | **8296** | **70.1%** | **33.2** | **38.7** | **64.9** |
| 20 | Nearest | 66.96 | 82.52 | 13.96 | 1.59 | 6586 | 56.3% | 26.3 | 74.3 | 77.3 |
| | Fastest | 68.55 | 84.53 | 14.31 | 1.67 | 6751 | 57.6% | 27.0 | 49.3 | 75.5 |
| | MinChgOpT | 69.95 | 86.27 | 14.59 | 1.72 | 6881 | 58.8% | 27.5 | 42.1 | 76.4 |
| | DynaThreshold | 74.12 | 91.88 | 15.75 | 2.01 | 7428 | 62.0% | 29.7 | 48.3 | 95.3 |
| | OptChg | 82.82 | 103.08 | 17.89 | 2.37 | 8439 | 71.2% | 33.8 | 56.6 | 65.2 |
| | **CongestionAware** | **84.16** | **104.78** | **18.18** | **2.44** | **8576** | **72.5%** | **34.3** | **40.9** | **64.2** |

Note that the day-ahead charging plan model P1 can be solved efficiently using a commercial solver to obtain good approximate solutions given a reasonable computational time limit when the problem size is not too large. For our computational study with 12 chargers and a 62 kWh battery size, we use a one-hour computational time limit for the c3000 scenario (see Appendix E for more details). The obtained solutions have gaps to the lower bound of around 14%. For the c4000 scenario, however, the problem becomes harder to solve. With a 4-hour computational time limit, we can obtain feasible solutions with an optimality gap of around 4%. When the problem size (in terms of the number of customers, vehicles, and fast/slow chargers) increases, we can reduce the computational time by decomposing it into smaller problem-size blocks with proportional customer demand, number of vehicles, and number of chargers in the system to obtain a charging plan for the vehicles of the block and replicate it for the vehicles of the other blocks. Another option is developing efficient heuristics to get good solutions, which remains a future research avenue of this study. We further analyse the impact of unbalanced fast and slow charger numbers and demand variation. The results are reported in Appendix D. The impact of using re-estimated $\delta_h$, based on post-realized vehicle energy consumption, on the performance of the CongestionAware policy is presented in Appendix H.

## 5. Conclusions and discussion

In this study, we develop a dynamic charging approach that coordinates vehicle dispatch and charging operations for electric ride-hailing systems under stochastic demand, variable energy prices, and congested charging stations. We focus on maximizing the total system profit by anticipating vehicles' energy needs and waiting time for charging on different chargers during the day to reduce vehicle unavailability and increase the service rate for customers. The proposed sequential MILP approach first determines charging time and target SoCs of vehicles for a long planning horizon, based on which an online reactive model optimizes vehicle-to-charger assignment for a short planning horizon to minimize total charging operational costs given the current system state. This reactive model adjusts vehicles' charging time and target SoCs after recharge based on vehicles' energy needs and waiting time on chargers to use congested chargers more efficiently. Five benchmark charging policies are used to compare the performance of the proposed method: Nearest, Fastest, Minimum charging operational time, dynamic hourly charging thresholds and optimization-based charging scheduling. We also model charger queueing more realistically, with no overlapping sessions on a charger and a maximum waiting time per queue, and we impose a minimum charging duration per session, which to our knowledge prior studies omit.

We conducted a simulation case study using the synthetic data generated from NYC yellow taxi data in a Manhattan-like area with two demand scenarios using a fleet of 100 EVs and limited fast and slow charging stations. The computational results show that the developed methodology outperforms the benchmark approaches in terms of higher profit and customer service rates under different scenarios. Overall, compared with the benchmark, the proposed approach increases total profit by 3.56%-14.90% for the scenario of 3000 customers per day and 3.40%-21.68% for that of 4000 customers per day. Similarly, the customer service rate is increased by 4.34%-18.33% and 3.87%-22.22% for the c3000 and c4000 scenarios, respectively. When increasing the battery size of vehicles and the number of chargers in the system, the service rate could be further improved from 67.1% (12 chargers and 62kWh battery of vehicles) to 76.0% (12 chargers and 82kWh battery) and 72.5% (20 chargers and 62kWh battery) for the CongestionAware policy with 4000 customers/day. Compared with the benchmark, the CongestionAware policy increases the service rate systematically (up to +4.83%-15.33%, for the c4000 scenario with 12 chargers and 82kWh battery). Moreover, the total charging waiting time is significantly reduced compared with the benchmark. The approach can help ride-hailing operators manage charging more efficiently when charging capacity is limited and demand is uncertain.

In our test instance generation, customer origin-destination (OD) pairs are generated synthetically within the study area rather than drawn directly from the observed taxi-zone trip records, which limits the spatial realism of the demand. This simplification may affect vehicle access time to charging stations, spatial distribution of available vehicles, and vehicle queueing at charging stations, in particular, clustered demand or clustered OD patterns. The magnitude of the impact depends on the interplays between the spatio-temporal relationship between customer demand, fleet size and charging infrastructure supply.

Several limitations should be noted. Travel times and energy consumption are modeled as deterministic averages, so the present analysis isolates uncertainty in customer demand and in charging-station queueing rather than travel-time or traffic-congestion uncertainty. We model travel distance and time based on shortest path distances rather than a realistic road network, which could be relaxed in future research by incorporating a realistic road network. The estimated waiting time and energy-consumption inputs for the day-ahead planning model P1 are based on the Fastest policy, so the performance of the CongestionAware policy is conditional on that policy. Using the Fastest policy may overestimate waiting time on fast chargers and thus affect the derived charging plan. Although the online charger assignment model uses updated charger availability and waiting times information, it does not eliminate the dependence from its day-ahead plan**.** In practice, these inputs can be better estimated with more accurate predictive models based on historical data (Soufiani et al., 2025). Finally, when scaling up the problem size to larger fleet and higher customer demand, the day-ahead charging plan may not be solved to optimality within reasonable computational time. Partitioning the fleet and charger capacities into smaller blocks is a potential direction for future research to obtain a sufficiently good feasible charging plan within a limited computational time budget.

Future extensions include generating realistic demand using observed NYC taxi zone data or origin-destination demand patterns and analysing the impact of spatial demand heterogeneity on the operational performance. Alternative operating policies can be designed; for example, allowing vehicles waiting at a charger to remain available for dispatch. This would combine early positioning for charging with continued availability for service, at the cost of additional empty travel and deferred charging. Our experiments provide a bound for the potential additional availability: CongestionAware shows 27 (c3000) and 43 (c4000) vehicle-hours per day waiting at chargers (about 2–3% of fleet service time).

Other research avenues include developing efficient solution approaches for scaling up the system, applying this approach for optimizing vehicle battery configuration, fleet size, and charging infrastructure planning, etc. Other interesting research avenues include integrating smart charging strategies to mitigate the impact of charging operations on the power grid during peak hours, extending this approach for different systems (e.g., shared mobility systems or regular bus services), or considering dynamic pricing to mitigate charging congestion. Incorporating different sources of uncertainty (e.g., stochastic travel times or traffic congestion) could be another interesting research avenue for a more realistic system performance evaluation.

**Acknowledgments**

The work was supported by the Luxembourg National Research Fund (C20/SC/14703944).

**Disclosure of interest**
The authors report there are no competing interests to declare.
The authors report ChatGPT was used to support the coding for the revision of the manuscript.

## Appendix A. Methodological comparison of existing representative literature and this study on charging scheduling for electric ride-hailing systems.

| Study | Decision architecture | Day-ahead charging plan | Online dispatch/charger assignment | Charging oper. modeling | HC/EP | MCD | ORC |
|---|---|---|---|---|---|---|---|
| **Ma (2021)** | Two-stage charging scheduling: day-ahead plan + online charger assignment | Yes. Approximate day-ahead charging plan (when to charge and each individual vehicle's target energy level after recharge); historical energy consumption, expected charging waiting time, and charging cost are used to derive charging schedules. | Use a MILP model to minimize the total weighted time of charging operations. Vehicle-customer dispatch is based on a lowest-marginal-system-cost approximation function. | Charger availability/queueing time is represented and minimized for vehicle–charger assignment over a rolling short horizon. | Partially. Heterogeneous charging stations but constant electricity prices | No | No |
| **Jamshidi et al. (2021)** | Three-stage charging scheduling: charging-duration planning, zone-based charger assignment and vehicle relocation, online vehicle dispatch, and charging. | Yes. Aggregated daily charging planning (how many vehicles should start charging and the planned charging duration for each 30-minute decision interval) using historical information. | Use a hierarchical procedure for charging planning: aggregate zone-based vehicle relocation, followed by vehicle–charger assignment subject to charging-station capacity constraints; charging queues are not explicitly modeled. | A SoC-threshold-driven approach is used to determine which vehicles to recharge; waiting queues are not explicitly modeled, and queues are prevented through charging-station capacity constraints. | Yes. Consider slow and fast chargers with time-dependent electricity prices. | No | No |
| **Zalesak & Samaranayake (2021)** | Online batch optimization for high-capacity electric ride-pooling with coordinated charging | Yes. A full-day charging-time plan. Charging-station assignment is determined online based on a short-horizon discretization. | Extend the ride-pooling assignment of Alonso-Mora et al. (2017) with a two-step iterative charging-schedule planning procedure; charging queues are not explicitly modeled. | Uses hard charging-station capacity constraints to avoid charging congestion; charging waiting/queueing times are not minimized. | No. Use a more realistic charging function, homogeneous charging stations, and constant electricity prices. | No | No |
| **Representative RL-based literature (Kullman et al., 2021; Yan et al., 2023)** | Markov Decision Process + approximate dynamic programming or reinforcement learning | No. Value functions are learned online. | Joint vehicle dispatch/charging decisions using reinforcement learning; charging waiting/queueing time is not considered. | Uses a charging-station-capacity-constraint-based approach; charging waiting/queueing times are neither represented nor minimized. | No. Homogeneous charging stations; energy costs are not considered. | No | No |
| **This study** | Two-stage charging scheduling: day-ahead plan + online charger assignment | Yes. Extend Ma (2021) with additional minimum charging duration constraints and incorporate overnight replenishment costs in the objective function for a 24-hour planning horizon. | Extend Ma (2021) in the following aspects: (a) incorporates a minimum charging-duration requirement in the online vehicle–charger assignment MILP model; (b) uses an adaptive target SoC based on anticipated energy needs until the end of the service period for each individual vehicle; (c) uses a reservation-based, queue-aware charger-assignment policy; and (d) uses a vehicle–customer matching MILP model to maximize profit for each batch dispatch interval. | Extends Ma (2021) with explicit charger-level waiting/congestion modeling and a maximum waiting-time restriction. Total charging operation time (including waiting/queueing time) is minimized. | Yes. Heterogeneous chargers with charger-type-specific fees and time-dependent electricity prices. | Yes | Yes. |

Remark: **HC/EP**: Heterogeneous chargers / energy prices; **MCD**: Minimum charging duration / maximum waiting time; **ORC**: Explicit overnight replenishment cost.

**Appendix B. Details of state transition equations.**

We use $x \overset{\mathcal{E}}{\rightarrow} y$ to describe the relevant state changes from $x$ to $y$ due to the implied event $\mathcal{E}$. The symbol $\leftarrow$ assigns the values from the right-hand side to a variable. The reader is referred to Horváth and Tamási (2025) for the discrete-event based decision process modelling and simulation. We describe how the system state evolves over time based on vehicle dispatch and charging decisions as follows.

– *Vehicle dispatch*

Consider a vehicle $v$, which has the state $\vec{x}_v(t) = (l_v(t), e_v(t), \text{idle}, \text{nothing})$ and is dispatched at decision time $t$ to serve customer $r$ (i.e., decision event $\mathcal{E}_t^{\text{dispatch}}(v,r)$). The customer's state is $\vec{x}_r(t) = (o_r, d_r, A_r, W_r(t), \text{waiting})$. $\mathcal{E}_t^{\text{dispatch}}(v,r)$ indicates that the dispatch event occurs at time $t$ and involves a vehicle $v$ and customer $r$. The sequential system state changes from $t$ to the completion of this request are described by Eq. (1).

$$\vec{x}_v(t) \xrightarrow{\mathcal{E}_t^{\text{dispatch}}(v,r)} \begin{cases} \vec{x}_v(t) \leftarrow (l_v(t), e_v(t), \text{serving}, \text{serving}) \\ \vec{x}_r(t) \leftarrow (o_r, d_r, A_r, W_r(t), \text{serving}) \end{cases} \xrightarrow{\mathcal{E}_{t'}^{\text{completed}}(v,r)} \begin{cases} \vec{x}_v(t') \leftarrow (d_r, e_v(t'), \text{idle}, \text{nothing}) \\ \vec{x}_r(t') \leftarrow (o_r, d_r, A_r, W_r(t'), \text{completed}) \end{cases} \quad \text{(B1)}$$

The first part on the right corresponds to changes in the system state after the decision is made (occurring at time $t$). The second part corresponds to the state changes when the vehicle arrives at the customer $r$'s destination at time $t' = t + \bar{t}_{l_v(t)o_r} + \bar{t}_{o_r d_r}$. $l_v(t')$ updates the vehicle's location to the customer $r$'s drop-off location, $d_r$, and the vehicle's SoC is updated as $e_v(t') = e_v(t) - \mu\left(\bar{d}_{l_v(t)o_r} + \bar{d}_{o_r d_r}\right)$, where $\mu$ is the energy consumption rate per kilometer traveled. $\bar{t}_{ij}$ and $\bar{d}_{ij}$ are the travel time and travel distance between $i$ and $j$, respectively. Travel time is calculated by dividing the total travel distance by the average vehicle speed. Customers' waiting times are updated as $W_r(t') = W_r(t) + \bar{t}_{l_v(t)o_r}$. Vehicles are not allowed to change their assigned tasks after the dispatch. If a customer $r$ cannot be served at time $t$, their waiting time at the next decision time $t + \Delta_t$ is updated as $W_r(t + \Delta_t) = W_r(t) + \Delta_t$. When $W_r(t + \Delta_t)$ exceeds the maximum predefined customer waiting times, the customer's status becomes rejected, and they leave the system.

– *Vehicle-charger assignment and* charging requests scheduling

Consider a vehicle with state $\vec{x}_v(t) = (l_v(t), e_v(t), \text{idle}, \text{nothing})$ that is assigned to recharge at a charger $s$ at time $t$ (i.e., decision event $\mathcal{E}_t^{\text{charging}}(v,s)$). The current state of a charger $s$ is $\vec{x}_s(t) = (l_s, b_s(t))$. If the vehicle departs at time $t$, its arrival time at $s$ can be calculated as $t_{arr}(v,s) = t + \bar{t}_{l_v(t)l_s}$. Given the current SoC $e_v(t)$ and the target battery level after recharge $\bar{E}_v(t)$, the amount of energy to be charged at $s$ can be calculated as $\Delta e_v(t) = \bar{E}_v(t) - e_v(t) + \mu\bar{d}_{l_v(t)l_s}$. The charging time can then be deduced as $\tau_{charging}(v,s) = \Delta e_v(t)/\varphi_s$ where $\varphi_s$ is the charging power of $s$. A charging request $g(v,s,t) = (v, t_{arr}(v,s), \Delta e_v)$ is then generated at request time $t$ and to be scheduled at charger $s$. $t_{arr}(v,s)$ is the arrival time, and $\Delta e_v$ is the requested amount of energy to charge.

We define a **charger reservation scheduler** $\omega(g(v,s,t), b_s(t), s)$ that schedules a charging request $g(v,s,t)$ at a charger given the current reservation/utilization state of the charger. Figure B1 illustrates how a charging request is scheduled at a charger based on its current reservation and utilization rates. Consider a charger where two charging sessions have already been scheduled. The red dashed box indicates a charging request, showing the arrival time of the vehicle and the end time of charging if the vehicle charges immediately upon arrival. Charging request (1) has no conflicts with other scheduled charging operations. Its waiting time is zero, while the earliest time to insert the charging request (2) is $t_3$. The waiting time of the request (2) is calculated $t_3 - t_1$.

When a charging request is accepted, the reservation/utilization state of the charger $s$ is updated accordingly at the request time $t$ and remains unchanged. The scheduler's output is the schedule for the charging request, denoted as $\zeta(v, t_{start}(v,s), t_{end}(v,s))$ with the start time $t_{start}(v,s)$ and end time $t_{end}(v,s)$ of the scheduled charging session. Vehicle's waiting time at the charger is then defined as $\tau_{wait}(v,s) = t_{start}(v,s) - t_{arr}(v,s)$. The vehicle's SoC after recharging at $s$ can be obtained as $e_v(t_{end}(v,s)) = e_v(t) - \mu \bar{d}_{l_v(t) l_s} + \Delta e_v(t)$. The charger reservation scheduler ensures correct modeling of charger capacity and vehicle queuing (no charging time overlaps will occur at any charger). Other scheduling rules (e.g., given priority to urgent requests, etc.) can be applied, if needed. The sequential system state changes triggered by a vehicle-charger assignment decision are described below.

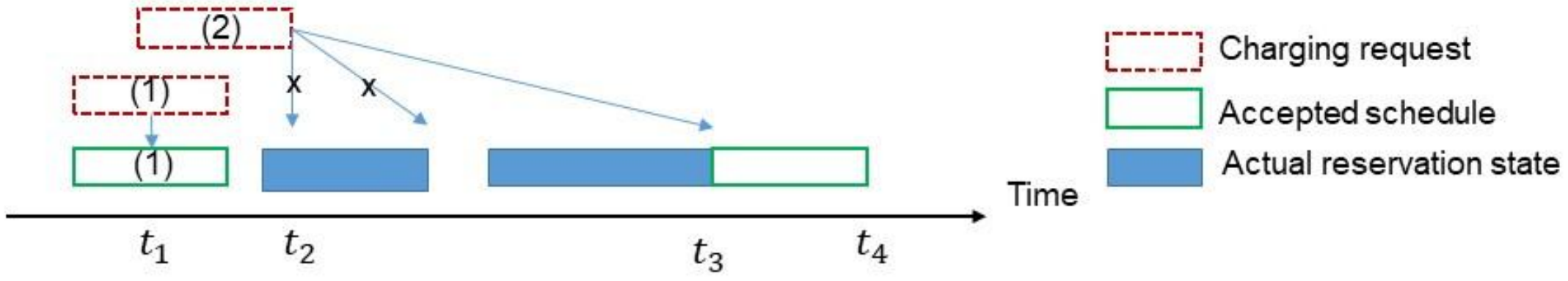

Figure B1. Charging request scheduling on a charger.

a) Charging request scheduling and vehicle arriving at the location of charger $s$

$$\vec{x}_v(t) \xrightarrow{\mathcal{E}_t^{\text{charging}}(v,s)} \begin{Bmatrix} \{\vec{x}_v(t) \leftarrow (l_v(t), e_v(t), \text{en_route}, \text{charging}) \\ \vec{x}_s(t) \leftarrow \left(l_s, b_s^{'}(t)\right) \end{Bmatrix} \quad \text{(B2)}$$
$$\xrightarrow{\mathcal{E}_{t_{arr}(v,s)}^{\text{arrive_chg}}(v,s)} \{\vec{x}_v\left(t_{arr}(v,s)\right) \leftarrow (l_s, e_v(t) - \mu \bar{d}_{l_v(t) l_s}, \text{arrive_chg}, \text{charging}\}$$

Given a vehicle-charger assignment decision $\mathcal{E}_t^{\text{charging}}(v,s)$ at time $t$ (output of Model P3), the vehicle's status and assigned task are updated to "en_route" and "charging", respectively. The reservation status of $s$ is updated immediately as $b'_s(t) = b_s(t) \cup \{\zeta(v, t_{start}(v,s), t_{end}(v,s))\}$. The second state change occurs when the vehicle arrives at $s$, and the vehicle's SoC is updated by subtracting the energy consumed during that trip. The vehicle's location and status are updated as $l_s$ and "arrive_chg", respectively.

b) If there is no need to wait, the vehicle charges immediately. The system state updates as follows.

$$\vec{x}_v\left(t_{arr}(v,s)\right)$$
$$\xrightarrow{\mathcal{E}_{t_{arr}(v,s)}^{\text{charging}}(v,s)} \{\vec{x}_v\left(t_{arr}(v,s)\right) \leftarrow \left(l_s, e_v\left(t_{arr}(v,s)\right), \text{charging}, \text{charging}\right)\} \quad \text{(B3)}$$
$$\xrightarrow{\mathcal{E}_{t_{end}}^{\text{idle}}(v,s)} \{\vec{x}_v\left(t_{end}(v,s)\right) \leftarrow \left(l_s, e_v\left(t_{end}(v,s)\right), \text{idle}, \text{nothing}\right)\}$$

The first state change indicates that the vehicle's status changes to "charging" upon arrival. The second state change indicates that the vehicle's status and assigned task change to "idle" and "nothing," respectively, at the end of the charging time. The vehicle's SoC is updated accordingly.

c) Otherwise, the vehicle needs to wait. If the waiting time $\tau_{wait}(v,s)$ is no greater than $\varpi$, the vehicle starts charging after waiting $\tau_{wait}(v,s)$ minutes. The system state is then updated as follows.

$$\vec{x}_v\left(t_{arr}(v,s)\right) \quad \text{(B4)}$$

$$\xrightarrow{\mathcal{E}^{\text{waiting}}_{t_{arr}(v,s)}(v,s)} \{\vec{x}_v(t) \leftarrow (l_s, e_v(t_{arr}(v,s)), \text{waiting}, \text{charging})\}$$
$$\xrightarrow{\mathcal{E}^{\text{charging}}_{t_{start}(v,s)}(v,s)} \{\vec{x}_v(t_{start}(v,s)) \leftarrow (l_s, e_v(t_{arr}(v,s)), \text{charging}, \text{charging})\}$$
$$\xrightarrow{\mathcal{E}^{\text{idle}}_{t_{end}(v,s)}(v,s)} \{\vec{x}_v(t_{end}(v,s)) \leftarrow (l_s, e_v(t_{end}(v,s)), \text{idle}, \text{nothing})\}$$

The first state change updates the vehicle's status from "arrived_chg" to "waiting" upon arrival at the charger. The second state change indicates that the vehicle begins charging at time $t_{start}(v,s,t)$, after which the vehicle's status becomes "charging". The third state change updates the vehicle's SoC after charging is complete, as well as its status to idle with no assigned tasks. Otherwise, the system state remains constant between events.

d) When $\tau_{wait}(v,s) > \varpi$, the following *vehicle-charger reassignment decision is evoked at time* $t_{arr}(v,s) + \varpi$.

– *Vehicle-charger reassignment*

Consider the current decision time $t_{reassign} = t_{arr}(v,s) + \varpi$, a newly reassigned charger $s_{new}$ and the vehicle's current state $\vec{x}_v(t_{reassign}) = (l_s, e_v(t_{reassign}), \text{waiting}, \text{charging})$. Similar to the *vehicle-charger assignment,* calculate the vehicle's arrival time from its current location at $s$ to that of $s_{new}$ and the amount of energy to charge. If the vehicle's SoC is sufficient to reach $s_{new}$, cancel the scheduled charging session for charger $s$ and update the reservation state of the charger. Then generate a new charging request for $s_{new}$ and get it scheduled by the scheduler. The system state is updated based on the same rules (Eqs. (B2)-(B4)). Otherwise, the vehicle waits at the current charger $s$ until its turn, after which the vehicle's state is updated accordingly. Note that the reassignment for $s_{new}$ depends on the applied policy. Here, we reassign vehicle $v$ to the charger with the least waiting time among all fast and slow chargers (see Appendix C).

**Appendix C. Comparison of different modeling approaches for vehicle queuing at chargers for the benchmark charging policies.**

In the literature, existing studies assume a simplified vehicle queuing behavior modeling at chargers/charging stations, i.e., vehicles wait in a queued charger/charging station for charging without time limits. This might not be realistic when the vehicle's queuing times are very long (e.g., several hours). In this case, drivers (vehicles) might prefer to move to other nearby chargers or the least-waiting-time chargers to recharge and return to serve customers earlier. To investigate the impact of queuing behavior at chargers, three modeling approaches are tested as follows. Note that vehicles' charger assignments are determined by the applied charging policies.

a. **Naïve queuing**: Vehicles wait in a queue for their target (assigned) charger without a time limit.
b. **Charger-chasing A**: Vehicles wait at a charger (current charger) with a maximum waiting time of 30 minutes, then move away to a **fast charger** (next charger) with the least waiting time when arriving at chargers' locations. If vehicles cannot reach the next charger due to insufficient SoCs, vehicles go to the closest charger to recharge. If unable to reach the closest one, vehicles wait at the current charger until its turn.
c. **Charger-chasing B**: Vehicles wait at a charger (current charger) with a maximum waiting time of 30 minutes, then move away to a **fast or slow** charger with the least waiting time when moving away from the queue on current chargers. If vehicles cannot reach the next charger due to insufficient SoCs, vehicles go to the closest charger to recharge or wait at the current charger until its turn if unable to reach the closest one.

Table C1 reports the performance of using different vehicle queuing modeling approaches. The upper block in Table C1 is related to the c3000 scenario, while the lower block is related to the c4000 scenario. The results show that using charger-chasing B has the highest service rates and profits for both demand scenarios. The total waiting time at chargers using charger-chasing B is systematically lower than that of charger-chasing A. However, using the naïve queue approach does not necessarily make it worse in terms of total waiting time compared to the charger-chasing approach. It may depend on the applied charging policy and uncertain queuing situations at chargers.

We further look into the details of the realized charging sessions for the c3000 scenario (The c4000 scenario has similar results; we neglect it here). Table C2 shows the results of using different queuing modeling approaches for the benchmark charging policies for the c3000 scenario. When using the naïve queuing approach, the number of realized charging operations is significantly fewer than the charger-chasing approaches due to the long waiting time on queued chargers. The charger-chasing B has the lowest average queuing time at chargers per charging session (13.2 minutes on average compared with the charger-chasing A (15.2 minutes on average) and the naïve queuing (86.8 minutes). The average amount of energy charged is similar. Still, the average charging times of charger-chasing A are lower than those of charger-chasing B as the latter considers both fast and slow chargers when vehicles go away from the current charger, resulting in more utilization of slow chargers to charge from vehicles' current SoC to $E_v^{max}$. Figure C1 reports the boxplots for vehicle queuing time at chargers for realized charging sessions for the naïve queue and the charger-chasing B. It shows that using the naïve queue approach might result in an unrealistically long queuing time at a charger, while using the charger-chasing B would not have this issue. Figure C2 reports vehicle queuing time at chargers (i.e., average queuing (waiting) times per charging session (2) in Table C2), taking into account the vehicle's waiting times for previously unsuccessful attempts on other chargers. The median of the average waiting time (2) of Charger-chasing B ranges from 24.7 to 40.3 minutes, which is much lower compared with that of Naïve, ranging from 64.3 to 91.4 minutes. The results confirm the charger-chasing B results in reducing unnecessary queuing times for vehicle recharging.

Table C1. Effect of different modeling approaches for vehicle queuing at chargers.

| Scenario | Queuing approach | Charging policy | PF | TR | TC | CC | ENG | SR | KMT | WT | CT |
|---|---|---|---|---|---|---|---|---|---|---|---|
| c3000 | Naïve queuing | Nearest | 60.42 | 73.87 | 12.31 | 1.15 | 5806 | 69.7% | 23.2 | 32.7 | 23.0 |
| | | Fastest | 63.95 | 78.74 | 13.22 | 1.57 | 6236 | 73.9% | 24.9 | 60.2 | 40.5 |
| | | MinChgOpT | 64.59 | 79.56 | 13.36 | 1.61 | 6302 | 74.8% | 25.2 | 55.7 | 45.3 |
| | | DynaThreshold | 66.93 | 82.78 | 14.05 | 1.81 | 6625 | 77.9% | 26.5 | 105.2 | 51.1 |
| | Charger-chasing A | Nearest | 63.39 | 78.06 | 13.13 | 1.54 | 6194 | 73.2% | 24.8 | 89.8 | 50.0 |
| | | Fastest | 63.88 | 78.72 | 13.26 | 1.57 | 6257 | 73.9% | 25.0 | 88.8 | 44.6 |
| | | MinChgOpT | 64.55 | 79.56 | 13.40 | 1.61 | 6319 | 74.7% | 25.3 | 76.1 | 46.9 |
| | | DynaThreshold | 66.07 | 82.13 | 14.23 | 1.83 | 6713 | 77.1% | 26.9 | 127.7 | 51.5 |
| | **Charger-chasing B** | **Nearest** | **63.27** | **77.91** | **13.11** | **1.53** | **6186** | **73.1%** | **24.7** | **82.9** | **58.8** |
| | | **Fastest** | **63.83** | **78.66** | **13.25** | **1.58** | **6250** | **73.8%** | **25.0** | **79.1** | **58.1** |
| | | **MinChgOpT** | **64.63** | **79.68** | **13.42** | **1.63** | **6330** | **74.8%** | **25.3** | **64.4** | **59.3** |
| | | **DynaThreshold** | **66.79** | **82.98** | **14.32** | **1.86** | **6755** | **78.0%** | **27.0** | **96.2** | **77.9** |
| c4000 | Naïve queuing | Nearest | 60.43 | 74.08 | 12.48 | 1.17 | 5887 | 51.6% | 23.5 | 45.7 | 25.3 |
| | | Fastest | 65.42 | 80.66 | 13.64 | 1.59 | 6434 | 55.2% | 25.7 | 80.3 | 43.5 |

| | | | | | | | | | | |
|---|---|---|---|---|---|---|---|---|---|---|
| | MinChgOpT | 66.56 | 82.08 | 13.87 | 1.65 | 6544 | 56.1% | 26.2 | 69.5 | 51.8 |
| | DynaThreshold | 69.46 | 86.00 | 14.65 | 1.89 | 6912 | 57.7% | 27.6 | 117.6 | 55.6 |
| Charger-chasing A | Nearest | 64.95 | 80.05 | 13.55 | 1.55 | 6392 | 54.8% | 25.6 | 122.2 | 52.6 |
| | Fastest | 65.82 | 81.21 | 13.78 | 1.62 | 6498 | 55.6% | 26.0 | 116.8 | 51.5 |
| | MinChgOpT | 66.37 | 81.86 | 13.85 | 1.64 | 6533 | 55.9% | 26.1 | 109.2 | 51.9 |
| | DynaThreshold | 68.81 | 85.56 | 14.82 | 1.93 | 6991 | 57.5% | 28.0 | 135.3 | 61.9 |
| **Charger-chasing B** | **Nearest** | **64.94** | **80.05** | **13.55** | **1.55** | **6393** | **54.9%** | **25.6** | **113.8** | **62.9** |
| | **Fastest** | **65.97** | **81.35** | **13.78** | **1.59** | **6502** | **55.6%** | **26.0** | **104.3** | **61.8** |
| | **MinChgOpT** | **66.64** | **82.17** | **13.91** | **1.63** | **6560** | **56.2%** | **26.2** | **96.3** | **62.9** |
| | **DynaThreshold** | **70.01** | **86.79** | **14.91** | **1.87** | **7034** | **58.3%** | **28.1** | **97.4** | **80.4** |

Table C2. Statistics of realized charging sessions for different modeling approaches for vehicle queuing at chargers (c3000 scenario).

| Queuing approach | Charging policy | Average # of realized charging sessions | Average queuing (waiting) times per charging session (1)[1] | Average queuing (waiting) times per charging session (2)[2] | Average charging times per charging session | Average charged energy per charging session |
|---|---|---|---|---|---|---|
| Naïve queuing | Nearest | 19.0 | 103.2 | 103.2 | 72.7 | 37.8 |
| | Fastest | 45.2 | 79.9 | 79.9 | 53.8 | 38.3 |
| | MinChgOpT | 46.8 | 71.4 | 71.4 | 58.0 | 39.4 |
| | DynaThreshold | 68.0 | 92.8 | 92.8 | 45.1 | 33.0 |
| Charger-chasing A | Nearest | 44.6 | 13.1 | 65.8 | 67.3 | 38.9 |
| | Fastest | 45.8 | 15.2 | 61.5 | 58.5 | 38.6 |
| | MinChgOpT | 46.4 | 15.5 | 47.3 | 60.6 | 39.8 |
| | DynaThreshold | 66.8 | 17.1 | 65.0 | 46.3 | 34.0 |
| Charger-chasing B | Nearest | 47.4 | 11.4 | 65.9 | 74.4 | 37.8 |
| | Fastest | 49.2 | 11.9 | 57.3 | 70.8 | 38.1 |
| | MinChgOpT | 51.0 | 14.3 | 44.6 | 69.8 | 38.9 |
| | DynaThreshold | 73.2 | 15.0 | 55.4 | 63.9 | 34.5 |

Remark: 1. Average times that a vehicle waits in a queue for recharging for a realized charging session. 2: Average queuing (waiting) times per charging session (1) + a vehicle's waiting times for previously unsuccessful attempts on other chargers. Time is measured in minutes; energy is measured in KWh.

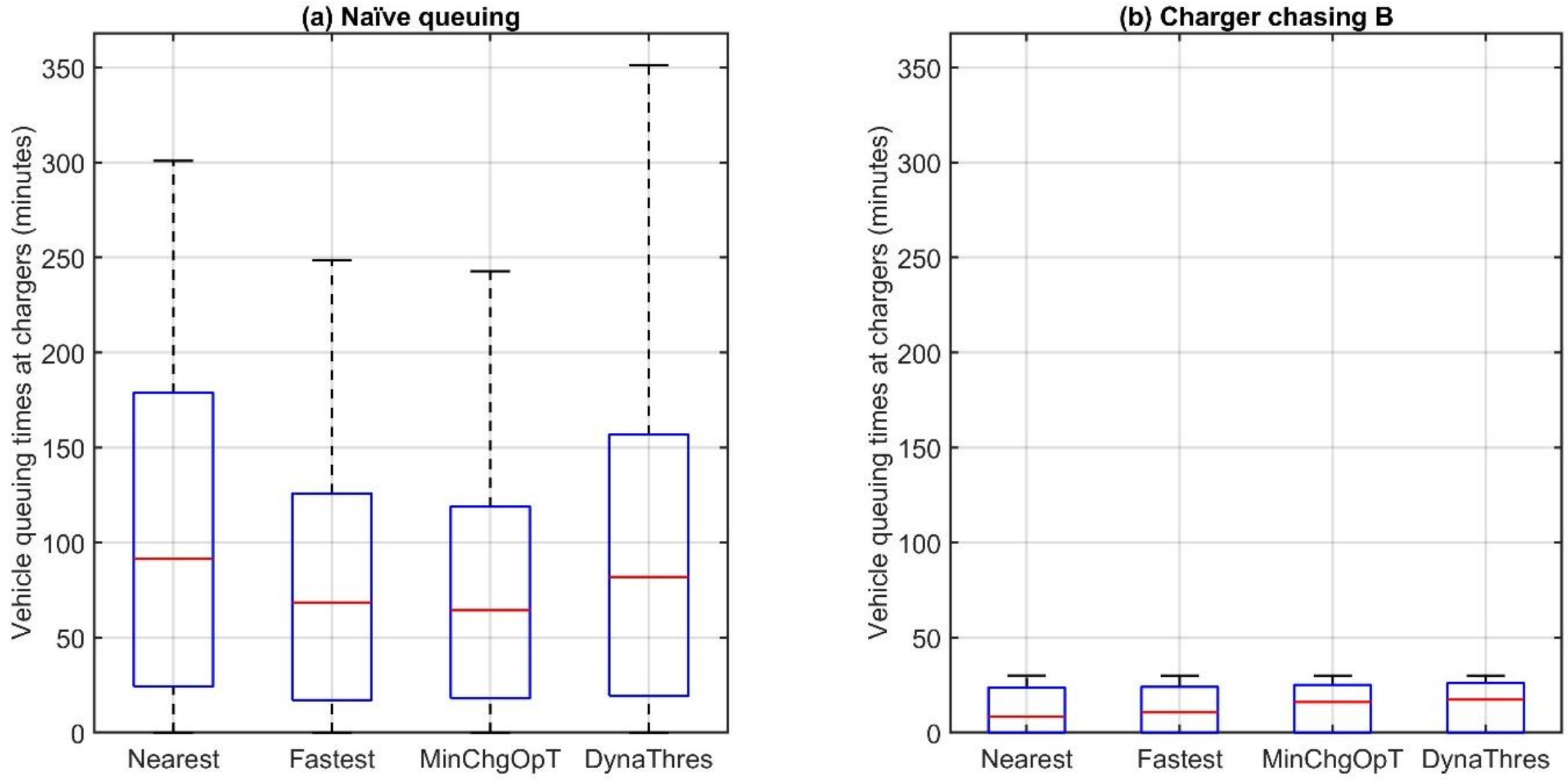

Figure C1. Boxplots for average queuing times that vehicles wait in a queue for recharging over the benchmark charging policies (i.e., average queuing (waiting) times per charging session (1)). Naïve queueing (on the left) and charger-chasing B (on the right).

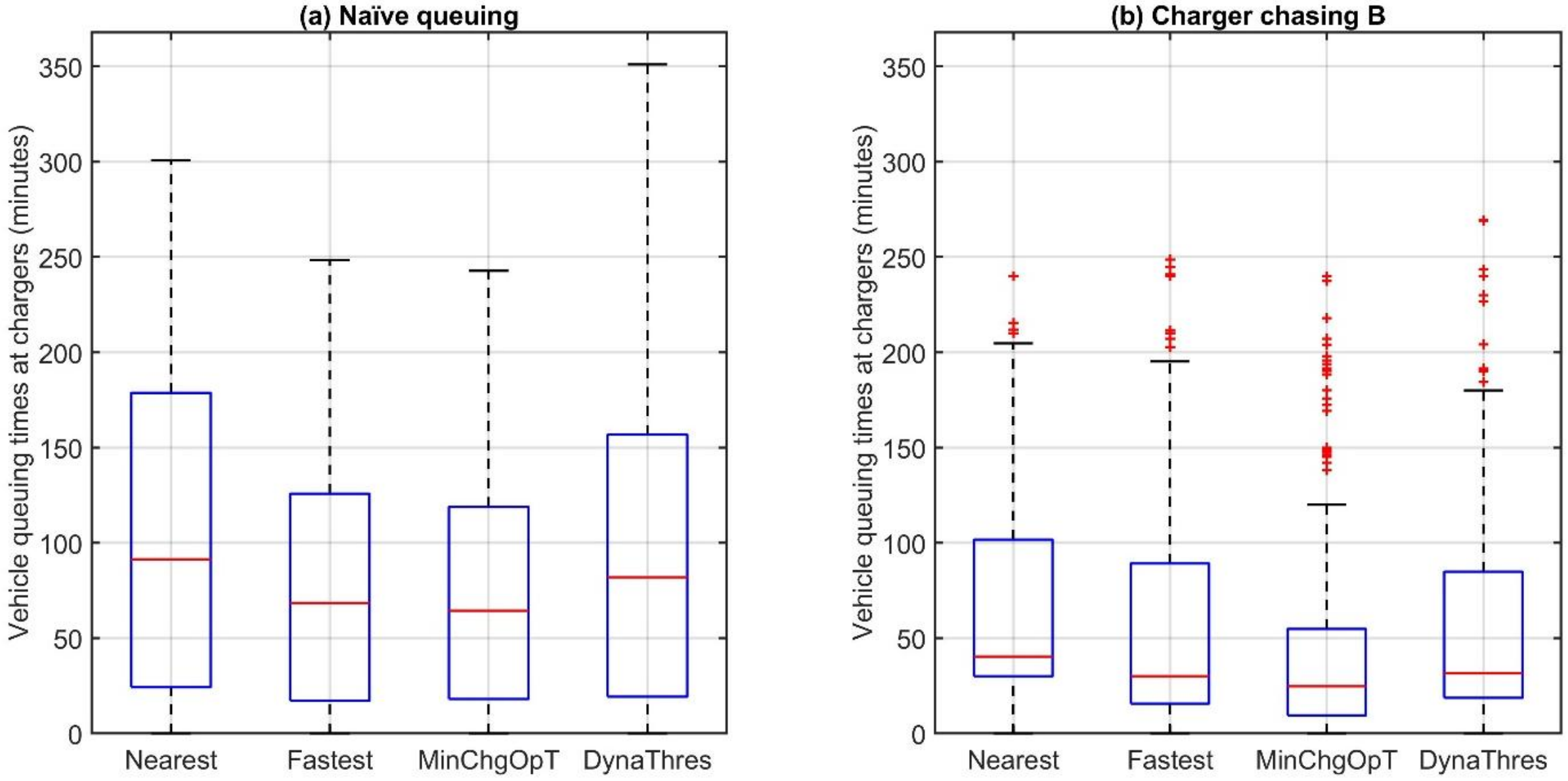

Figure C2. Boxplots for average queuing times that vehicles wait in a queue for recharging over the benchmark charging policies (i.e., average queuing (waiting) times per charging session (2)). Naïve queueing (on the left) and charger-chasing B (on the right).

**Appendix D. Impact of charging station availability and imbalanced number of fast and slow chargers**

We evaluate the effect of extreme charging station availability and demand variability on the performance of different charging policies, including: i) extremely low/high charging station availability; ii) imbalanced number of fast and slow chargers; iii) variation of customer demand. We design three scenarios as detailed in Table D1. The computational time details are reported in Appendix E.

Table D1. Scenarios for the system evaluation with different charging station availability and imbalanced numbers of fast and slow chargers and customer demand variation.

| Scenario | Description | Tested values | | |
|---|---|---|---|---|
| | | Fleet size | Demand (number of requests/day) | Number of fast and slow chargers (fast, slow) |
| 1 | Charging station availability | 100 | 3000 | (2,2),(20,20) |
| 2 | Imbalanced numbers of fast and slow chargers | 100 | 3000 | (2,10), (4,20) |
| 3 | Variation of customer demand | 100 | 1000,2000,3000,4000,5000,6000 | (6,6) |

Remark: Battery size = 62 kWh. The number and location of fast and slow charging stations remain the same (2 fast and 2 slow charging stations (see Figure 2) with balanced numbers of chargers). The other parameters are based on Table 3.

**Scenario 1**

We solve the day-ahead charging schedule planning problem with one-hour computational time to obtain a good solution. For the case of extremely low numbers of chargers (2 fast and 2 slow), we cannot obtain a feasible solution for the fleet size of 100 vehicles. Consequently, we reduce the number of vehicles for P1 gradually (reduce 10 vehicles at a time) to obtain feasible solutions for at most 40 vehicles. The remaining 60 vehicles go to recharge when idle and their SoCs are below the threshold $\theta$. The performance of the CongestionAware and other benchmark charging policies is shown in Table D2. As expected, the results demonstrate that the CongestionAware policy outperforms the rule-based benchmark when the availability of charging facilities is low (+6.87% total profit compared to the best rule-based benchmark). The total waiting time for the benchmark policies are extremely high (more than 159 hours in total) when vehicles' charging operations are not well coordinated to avoid peak congestions or myopically charging to $E_v^{max}$ (i.e. 80%); while the CongestionAware policy and OptChg avoid recharging vehicles during charging peak hours with a partial recharge policy, resulting in significantly lower charging waiting time (12.1-12.9 hours in total). On the other hand, when the fast and slow chargers are abundant, the two policies achieve comparable performance under abundant charging capacity. The charging waiting time is significantly reduced (0 hour) compared with the DynaThreshold (16.6 hours) and other rule-based policies.

Table D2. System performance of different charging policies with an extremely low/high number of fast and slow chargers.

| # of chargers[1] | Charging policy | PF | TR | TC | CC | ENG | SR | KMT | WT | CT |
|---|---|---|---|---|---|---|---|---|---|---|
| (2,2) | Nearest | 59.57 | 73.08 | 12.35 | 1.15 | 5826 | 68.8% | 23.3 | 163.7 | 23.2 |
| | Fastest | 59.42 | 72.88 | 12.32 | 1.15 | 5812 | 68.7% | 23.2 | 159.2 | 22.4 |
| | MinChgOpT | 59.49 | 73.00 | 12.36 | 1.15 | 5828 | 68.8% | 23.3 | 160.9 | 22.3 |
| | DynaThreshold | 57.29 | 71.32 | 12.78 | 1.25 | 6026 | 66.6% | 24.1 | 220.1 | 29.4 |
| | OptChg | 62.60 | 76.64 | 12.86 | 1.18 | 6066 | 72.7% | 24.3 | 12.1 | 17.3 |

| | CongestionAware | **63.66** | **77.96** | **13.08** | **1.22** | **6168** | **74.1%** | **24.7** | **12.9** | **17.5** |
|---|---|---|---|---|---|---|---|---|---|---|
| (20,20) | Nearest | 64.54 | 79.53 | 13.38 | 1.61 | 6312 | 74.6% | 25.2 | 42.3 | 85.8 |
| | Fastest | 67.17 | 82.87 | 13.96 | 1.74 | 6584 | 78.1% | 26.3 | 16.9 | 72.4 |
| | MinChgOpT | 69.74 | 86.08 | 14.48 | 1.86 | 6831 | 81.5% | 27.3 | 5.5 | 62.0 |
| | DynaThreshold | 72.78 | 90.38 | 15.47 | 2.13 | 7298 | 86.3% | 29.2 | 16.6 | 89.8 |
| | OptChg | 73.51 | 90.55 | 15.33 | 1.72 | 7229 | 87.6% | 28.9 | 0.3 | 35.8 |
| | **CongestionAware** | **73.17** | **90.24** | **15.33** | **1.74** | **7233** | **87.4%** | **28.9** | **0.0** | **36.3** |

Remark: (fast, slow)

**Scenario 2**

This scenario considers the case with relatively low availability of fast chargers. We solve P1 with a 1-hour computational time to obtain a good charging schedule for the fleet. The KPIs of different charging policies are shown in Table D3. The results show that the CongestionAware policy has the highest total profit and customer service rate. When doubling the number of both types of chargers, the customer service rate of the CongestionAware policy is improved from 78.2% to 84.9%, resulting in a profit increase of 7.15%. The total waiting time is reduced significantly for the rule-based benchmark approaches with doubled chargers, but remain much higher than the CongestionAware and OptChg policies. The profit of the CongestionAware policy is higher than that of the benchmark approaches.

Table D3. System performance for different charging policies with an imbalanced number of fast and slow chargers.

| # of chargers[1] | Charging policy | PF | TR | TC | CC | ENG | SR | KMT | WT | CT |
|---|---|---|---|---|---|---|---|---|---|---|
| (2,10) | Nearest | 60.30 | 73.97 | 12.45 | 1.21 | 5873 | 69.5% | 23.5 | 122.3 | 59.0 |
| | Fastest | 60.02 | 73.65 | 12.42 | 1.21 | 5859 | 69.3% | 23.4 | 117.9 | 57.9 |
| | MinChgOpT | 60.03 | 73.68 | 12.43 | 1.21 | 5865 | 69.4% | 23.5 | 121.4 | 57.1 |
| | DynaThreshold | 59.73 | 74.06 | 12.99 | 1.34 | 6125 | 69.0% | 24.5 | 173.5 | 78.9 |
| | OptChg | 63.97 | 78.38 | 13.18 | 1.23 | 6219 | 74.3% | 24.9 | 14.4 | 31.2 |
| | **CongestionAware** | **66.84** | **82.22** | **13.97** | **1.41** | **6590** | **78.2%** | **26.4** | **24.7** | **45.6** |
| (4,20) | Nearest | 61.82 | 76.04 | 12.79 | 1.44 | 6032 | 71.3% | 24.1 | 62.5 | 100.8 |
| | Fastest | 62.13 | 76.47 | 12.88 | 1.46 | 6074 | 71.8% | 24.3 | 65.9 | 89.7 |
| | MinChgOpT | 62.39 | 76.80 | 12.94 | 1.47 | 6106 | 72.1% | 24.4 | 65.5 | 85.9 |
| | DynaThreshold | 65.34 | 80.92 | 13.88 | 1.70 | 6547 | 75.8% | 26.2 | 70.7 | 129.0 |
| | OptChg | 68.52 | 84.21 | 14.26 | 1.44 | 6725 | 80.3% | 26.9 | 21.7 | 47.9 |
| | **CongestionAware** | **71.62** | **88.40** | **15.12** | **1.67** | **7133** | **84.9%** | **28.5** | **27.4** | **54.3** |

Remark: (fast, slow)

**Scenario 3**

In this scenario, we vary daily customer demand from 1000 to 6000 with an interval of 1000 customers. As for the previous scenarios, a fleet of 100 vehicles with a 62 kWh battery is considered. The charging facility has 6 fast and 6 slow chargers. For low-demand cases with 1000 and 2000 customers/day, there are no charging operations scheduled from P1 (CongestionAware). The total charging time during the day of the CongestionAware policy is a small fraction of the other benchmark approaches. At very low demand (1000–2000 customers/day), the CongestionAware policy does not outperform the benchmark; at 2000 customers/day its profit and service rate are in fact marginally lower than the best benchmark (see Table D4), because the day-ahead plan and partial-recharge anticipation provide little benefit when charging capacity is not a binding constraint. The advantage of the proposed coordination emerges only once charging congestion becomes significant, which is the operating regime this study targets. However, for higher customer demand cases (from 3000 to 6000 customers/day), the CongestionAware policy systematically outperforms the benchmark approaches with the least waiting times and the highest profit (see Table D4 and Figure D1). Note that the location of charging stations may affect the access costs of charging operations. The sensitivity analysis related to the impact of charging station locations remains for future extensions of this study.

Table D4. Comparison of the KPIs for different charging policies for low customer demand.

| Scenario | Charging policy | PF | TR | TC | CC | ENG | SR | KMT | WT | CT |
|---|---|---|---|---|---|---|---|---|---|---|
| c1000* | Nearest | 27.51 | 33.14 | 5.21 | 0.42 | 2459 | 96.8% | 9.8 | 0.2 | 2.6 |
| | Fastest | 27.50 | 33.14 | 5.22 | 0.43 | 2462 | 96.8% | 9.8 | 0.0 | 1.9 |
| | MinChgOpT | 27.50 | 33.14 | 5.21 | 0.43 | 2459 | 96.8% | 9.8 | 0.0 | 2.0 |
| | DynaThreshold | 27.48 | 33.14 | 5.23 | 0.44 | 2465 | 96.8% | 9.9 | 0.2 | 2.6 |
| | OptChg | 27.54 | 33.14 | 5.21 | 0.39 | 2457 | 96.8% | 9.8 | 0.0 | 0.0 |
| | **CongestionAware** | **27.54** | **33.14** | **5.21** | **0.39** | **2457** | **96.8%** | **9.8** | **0.0** | **0.0** |
| c2000 | Nearest | 52.20 | 63.96 | 10.60 | 1.15 | 5002 | 93.0% | 20.0 | 27.3 | 37.9 |
| | Fastest | 52.50 | 64.40 | 10.69 | 1.20 | 5043 | 93.8% | 20.2 | 16.3 | 34.5 |
| | MinChgOpT | 52.72 | 64.70 | 10.73 | 1.25 | 5061 | 94.4% | 20.2 | 12.6 | 34.6 |
| | DynaThreshold | 52.83 | 64.90 | 10.81 | 1.26 | 5099 | 94.7% | 20.4 | 14.1 | 36.5 |
| | OptChg | 51.93 | 63.11 | 10.34 | 0.85 | 4875 | 91.8% | 19.5 | 0.1 | 4.3 |
| | **CongestionAware** | **51.75** | **62.85** | **10.27** | **0.83** | **4846** | **91.4%** | **19.4** | **0.0** | **3.4** |
| c5000 | Nearest | 64.83 | 79.79 | 13.52 | 1.45 | 6376 | 42.4% | 25.5 | 102.5 | 56.2 |
| | Fastest | 66.08 | 81.35 | 13.77 | 1.50 | 6496 | 43.1% | 26.0 | 89.2 | 56.6 |
| | MinChgOpT | 66.28 | 81.59 | 13.82 | 1.49 | 6519 | 43.1% | 26.1 | 85.1 | 56.2 |
| | DynaThreshold | 70.46 | 87.27 | 14.94 | 1.87 | 7045 | 45.1% | 28.2 | 114.2 | 79.8 |
| | OptChg | 76.29 | 94.42 | 16.19 | 1.93 | 7639 | 49.1% | 30.6 | 41.8 | 54.9 |
| | **CongestionAware** | **79.85** | **99.10** | **17.09** | **2.16** | **8060** | **51.7%** | **32.2** | **41.4** | **60.5** |
| c6000 | Nearest | 64.89 | 79.75 | 13.46 | 1.39 | 6348 | 33.8% | 25.4 | 116.9 | 54.4 |
| | Fastest | 66.26 | 81.49 | 13.76 | 1.47 | 6489 | 34.4% | 26.0 | 97.0 | 56.9 |
| | MinChgOpT | 66.65 | 81.96 | 13.83 | 1.48 | 6522 | 34.5% | 26.1 | 93.2 | 55.9 |
| | DynaThreshold | 69.50 | 85.71 | 14.52 | 1.70 | 6847 | 35.6% | 27.4 | 75.8 | 69.9 |
| | OptChg | 76.49 | 94.58 | 16.17 | 1.92 | 7628 | 39.1% | 30.5 | 42.3 | 54.8 |
| | **CongestionAware** | **80.02** | **99.22** | **17.05** | **2.15** | **8043** | **41.1%** | **32.2** | **35.7** | **58.6** |

Remark: cXX means XX requests. See Table 5 for the results of c3000 and c4000.

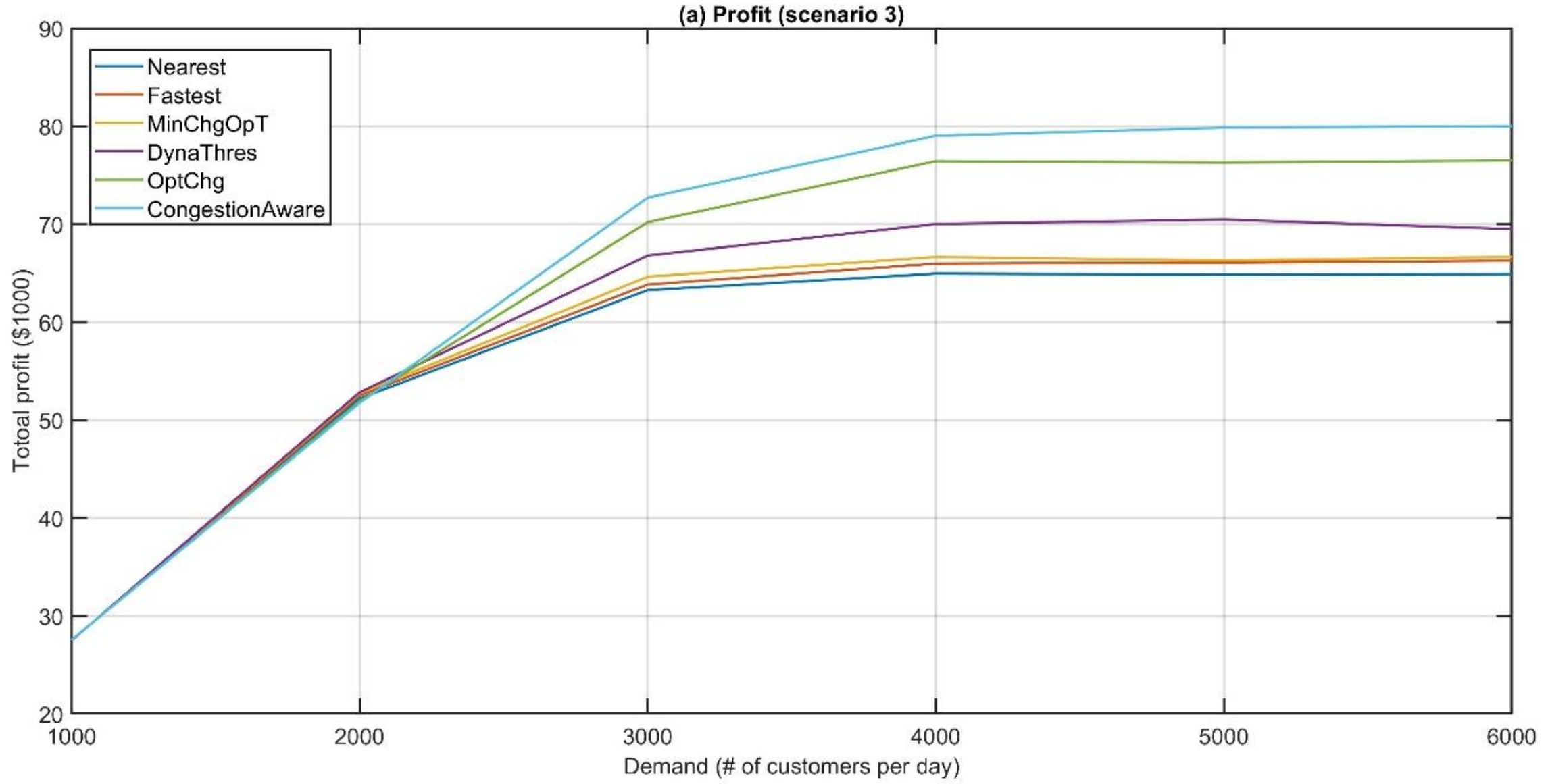

Figure D1. Profits of different charging policies with a variation of customer demand from 1000 to 6000 customers/day.

**Appendix E. Impact of problem size on the computational times of different models (P1-P3) for the CongestionAware policy.**

In this appendix, we report detailed computational times for solving P1-P3 problems of the CongestionAware policy and simulation times for the experiments in Section 4. To further evaluate the computational bottleneck when scaling up the problem size in terms of the number of vehicles, number of customers per day, and number of chargers, we generate four larger test datasets with the number of customers and fleet size up to 14000 customers/day and 500 vehicles. The number of fast and slow chargers is increasing accordingly (over the four charging stations). For each dataset, we randomly generate 15 test instances using the same approach (Section 4.1), in which 10 test instances are used to estimate the model parameters and the remaining 5 test instances are used to evaluate the performance of the different charging policies. The characteristics of the larger test datasets are shown in Table E1.

Table E1. Characteristics of the four larger test datasets.

| Demand (number of request/day) | Fleet size | Battery (kWh) | $\lvert S^{fast} \rvert$ | $\lvert S^{slow} \rvert$ |
|---|---|---|---|---|
| 8000 | 200 | 62 | 12 | 20 |
| 10000 | 300 | 62 | 20 | 30 |
| 12000 | 400 | 62 | 30 | 40 |
| 14000 | 500 | 62 | 40 | 50 |

Remark: Battery size = 62kWh. The number and location of fast and slow charging stations remain the same (2 fast and 2 slow charging stations). The other parameters are based on Table 3.

Table E2 reports the CPU time (in seconds) and the relative gaps to the lower bound for solving the P1 problem. We use a one-hour computational time to solve approximately P1 for the computational studies in Section 4. If no feasible solutions can be found, we extend the computational time with additional hours. The relative gaps range from around 0% to 19.07% for the c1000-c4000 scenarios. We can observe that given the same number of customer demand and fleet size to schedule, the availability of chargers has a significant impact on the gaps to the lower bound, given the same computational time limits, in particular when the number of chargers is extremely low/high. For the three largest test datasets, we report the gaps to the lower bound with 1 and 4 hours of computational time. The results show that the gaps to the lower bound depend on the problem instances. The c12000 scenario presents a particular difficult case where the optimality gap is relatively high (47.59%) with a 4-hour computational time limit. The results highlight the need for further research to enhance the computational scalability of P1 by strengthening its mathematical formulation and developing tailored decomposition algorithms and heuristics. Moreover, the impact of P1 solution quality on overall operational profit remains to be investigated. Note that given an unknown and stochastic customer demand, the day-ahead charging schedule P1 aims to obtain approximate charging schedules for the fleet to determine where and how much energy to charge, given the operator's available computational time and resource limits.

Table E2. CPU times and gaps to the lower bounds for solving P1 for the computational studies in Section 4 and larger test datasets.

| # of requests | $V$ | Battery (kWh) | $\lvert S^{fast}\rvert$ | $\lvert S^{slow}\rvert$ | $\Delta_\ell$ | CPU (sec.) | Gap to lower bound |
|---|---|---|---|---|---|---|---|
| 1000 | 100 | 62 | 6 | 6 | 30 | 0 | 0.00% |
| 2000 | 100 | 62 | 6 | 6 | 30 | 0 | 0.00% |
| 3000 | 100 | 62 | 6 | 6 | 30 | 3600 | 14.06% |
| 3000 | 40* | 62 | 2 | 2 | 30 | 3600 | 19.07% |
| 3000 | 100 | 62 | 20 | 20 | 30 | 3600 | 2.02% |
| 3000 | 90* | 62 | 2 | 10 | 30 | 3600 | 12.66% |
| 3000 | 100 | 62 | 4 | 20 | 30 | 3600 | 4.96% |
| 3000 | 100 | 62 | 8 | 8 | 30 | 3600 | 7.09% |
| 3000 | 100 | 62 | 10 | 10 | 30 | 3600 | 4.59% |
| 3000 | 100 | 72 | 6 | 6 | 30 | 3600 | 8.54% |
| 3000 | 100 | 82 | 6 | 6 | 30 | 3600 | 0.0% |
| 4000 | 100 | 62 | 6 | 6 | 20 | 14400 | 3.93% |
| 4000 | 100 | 72 | 6 | 6 | 20 | 3600 | 6.75% |
| 4000 | 100 | 82 | 6 | 6 | 20 | 3600 | 8.20% |
| 4000 | 100 | 62 | 8 | 8 | 20 | 3600 | 3.18% |
| 4000 | 100 | 62 | 10 | 10 | 20 | 3600 | 2.42% |
| 5000 | 90* | 62 | 6 | 6 | 20 | 3600 | 4.52% |
| 6000 | 80* | 62 | 6 | 6 | 20 | 3600 | 3.90% |
| 8000 | 200 | 62 | 12 | 20 | 20 | 3600 | 5.13% |
| 10000 | 300 | 62 | 20 | 30 | 20 | 3600 | 12.70% |
| 10000 | 300 | 62 | 20 | 30 | 20 | 14400 | 11.48% |
| 12000 | 400 | 62 | 30 | 40 | 20 | 3600 | 47.90% |
| 12000 | 400 | 62 | 30 | 40 | 20 | 14400 | 47.59% |
| 14000 | 500 | 62 | 40 | 50 | 20 | 3600 | 3.98% |
| 14000 | 500 | 62 | 40 | 50 | 20 | 14400 | 3.98% |

Remark: The fleet size that can be solved to obtain feasible solutions, given the energy consumption rate, the battery size of vehicles, and the number of fast and slow chargers for that test instance. We set $\Delta_\ell$=20 for the problem with a number of customers greater than 3000, instead of 30 minutes, to obtain feasible solutions.

Table E3 reports the CPU time to solve vehicle dispatching (every minute) and online vehicle-charger assignment (see Table 1) for different test instances with the number of customers ranging from 3000 to 14000 from 6:00-20:00 per day. For the c3000 to c6000 test instance, the fleet size is 100 vehicles, and the number of fast and slow chargers is 6 each. The setting for c8000-c14000 is shown in Table E1. For c3000, the average number of customers $|R_t|$ and vehicles $|V_t|$ is 9.7 and 28.2, respectively; while for c14000, $|R_t|$ and $|V_t|$ become 24.7 and 166.3, respectively. The standard deviations are relatively high for both variables, too. The average CPU time is 0.0018 seconds for c3000 and 0.0189 seconds for c14000, respectively, fast enough for real-time applications with the problem size of P2 up to hundreds of customers and hundreds of vehicles to dispatch. For P3, the problem size depends on the number of vehicles $|\widetilde{\Omega}_{t(i)}|$ (changing over time) and that of chargers $|S|$ (fixed). For the largest case of c14000, $|S|$ is 90, and the average and standard deviation of $|V_t|$ are 1.9 and 1.1, respectively. The average and standard deviation of the CPU time are 0.0063 and 0.0060 seconds, respectively, showing that solving the P3 problem is very fast, suitable for real-time applications.

Table E3. Problem size and CPU times for solving P2 (vehicle dispatching) and P3 (online vehicle-charger assignment) for different test instances.

| # of requests | $P2$ | | | | | | $P3$ | | | | |
|---|---|---|---|---|---|---|---|---|---|---|---|
| | $\vert R_t\vert$ | | $\vert V_t\vert$ | | CPU (sec.) | | $\vert S\vert$ | $\vert\widetilde{\Omega}_{t(i)}\vert$ | | CPU (sec.) | |
| | avg. | s.d. | avg. | s.d. | avg. | s.d. | | avg. | s.d. | avg. | s.d. |
| c3000 | 9.7 | 7.7 | 28.2 | 19.3 | 0.0018 | 0.0034 | 12 | 1.4 | 0.7 | 0.0013 | 0.0025 |
| c4000 | 21.8 | 13.2 | 19.4 | 18.1 | 0.0021 | 0.0037 | 12 | 7.7 | 3.2 | 0.0036 | 0.0041 |
| c5000 | 34.0 | 17.9 | 19.7 | 18.8 | 0.0032 | 0.0045 | 12 | 8.2 | 2.8 | 0.0043 | 0.0045 |
| c6000 | 46.2 | 21.6 | 20.6 | 19.5 | 0.0041 | 0.0054 | 12 | 8.5 | 2.7 | 0.0044 | 0.0045 |
| c8000 | 41.3 | 24.6 | 32.5 | 36.5 | 0.0054 | 0.0068 | 32 | 18.9 | 8.0 | 0.0125 | 0.0083 |
| c10000 | 32.9 | 27.7 | 68.9 | 53.0 | 0.0090 | 0.0089 | 50 | 31.7 | 11.6 | 0.0401 | 0.0185 |
| c12000 | 21.3 | 14.5 | 99.9 | 75.2 | 0.0102 | 0.0088 | 70 | 2.0 | 1.2 | 0.0058 | 0.0045 |
| c14000 | 24.7 | 21.1 | 166.3 | 86.1 | 0.0189 | 0.0224 | 90 | 1.9 | 1.1 | 0.0063 | 0.0060 |

Table E4 further reports the number of times that P2 and P3 were evoked during the simulation when solving a test instance in question. For P2, this number is bounded by the planning horizon (840, one-minute decision epoch for an 840-minute planning horizon from 6:00-20:00). For P3, the evoked number of executions varies substantially across instances, depending on the interactions of the charging supply and demand. As shown in Table 1 (Step 13), when there are no feasible vehicle-to-charger assignments, a series of attempts is executed by removing one vehicle with the highest SoC at a time until feasible solutions are found. Table E5 reports the CPU time for running the simulation using different charging policies for different test instances. For the CongestionAware policy, it takes around 3.7 minutes to finish a simulation for the most difficult c10000 test instance with the P3 model (Online vehicle-charger assignment) solved more frequently (see Table E4). Table E6 reports the KPIs using different charging policies for the larger test datasets. The results show that CongestionAware generally outperforms the benchmark charging policies, with a marginal exception for the OptChg for the c14000 scenario.

Table E4. The average number of executions for P2 and P3 for different test instances.

| # of requests | $P2$ | $P3$ |
|---|---|---|
| c3000 | 797 | 165 |
| c4000 | 827 | 934 |
| c5000 | 817 | 1473 |
| c6000 | 801 | 1574 |
| c8000 | 839 | 2092 |
| c10000 | 839 | 3327 |
| c12000 | 839 | 388 |
| c14000 | 839 | 259 |

Table E5. Average computational time of one simulation run for different charging policies (in seconds).

| # of requests | Nearest | Fastest | MinChgOpT | DynaThreshold | OptChg | CongestionAware |
|---|---|---|---|---|---|---|
| c1000 | 1.7 | 0.9 | 0.9 | 0.8 | 2.6 | 4.5 |
| c2000 | 2.3 | 1.4 | 1.5 | 1.4 | 3.6 | 3.6 |
| c3000 | 3.4 | 2.6 | 2.9 | 2.0 | 5.9 | 4.1 |
| c4000 | 5.8 | 4.9 | 4.8 | 3.7 | 11.7 | 8.5 |
| c5000 | 8.3 | 8.0 | 8.4 | 6.0 | 19.6 | 13.5 |
| c6000 | 11.1 | 9.7 | 9.8 | 10.0 | 22.8 | 15.1 |
| c8000 | 23.8 | 24.4 | 21.8 | 18.6 | 127.0 | 61.7 |
| c10000 | 34.2 | 33.3 | 33.6 | 26.2 | 302.8 | 222.3 |
| c12000 | 50.2 | 47.7 | 42.4 | 28.4 | 32.9 | 25.0 |
| c14000 | 59.7 | 52.3 | 50.5 | 31.8 | 38.9 | 41.7 |

Table E6. Comparison of the KPIs for different charging policies for larger test datasets.

| Demand | Charging policy | PF | TR | TC | CC | ENG | SR | KMT | WT | CT |
|---|---|---|---|---|---|---|---|---|---|---|
| c8000 | Nearest | 130.6 | 160.4 | 26.90 | 2.8 | 12695 | 55.2% | 50.8 | 171.9 | 136.3 |
| | Fastest | 133.0 | 163.4 | 27.40 | 3.0 | 12944 | 55.9% | 51.8 | 135.5 | 138.3 |
| | MinChgOpT | 134.5 | 165.2 | 27.70 | 3.0 | 13075 | 56.5% | 52.3 | 128.5 | 135.1 |
| | DynaThreshold | 140.5 | 173.3 | 29.40 | 3.5 | 13847 | 58.1% | 55.4 | 133.6 | 165.0 |
| | OptChg | 155.7 | 192.3 | 32.8 | 3.9 | 15463 | 64.7% | 61.9 | 83.7 | 123.6 |
| | **CongestionAware** | **163.2** | **202.6** | **34.9** | **4.5** | **16447** | **68.5%** | **65.8** | **102.3** | **148.1** |
| c10000 | Nearest | 196.6 | 240.2 | 39.40 | 4.2 | 18565 | 66.8% | 74.3 | 187.4 | 196.2 |
| | Fastest | 201.0 | 245.9 | 40.30 | 4.5 | 19032 | 68.2% | 76.1 | 141.1 | 194.5 |
| | MinChgOpT | 205.8 | 251.7 | 41.30 | 4.6 | 19458 | 69.8% | 77.8 | 85.1 | 172.6 |
| | DynaThreshold | 214.9 | 264.2 | 44.00 | 5.3 | 20758 | 73.1% | 83.0 | 129.6 | 239.6 |
| | OptChg | 233.8 | 287.8 | 48.2 | 5.8 | 22756 | 81.2% | 91.0 | 120.7 | 167.5 |
| | **CongestionAware** | **239.0** | **294.2** | **49.3** | **5.9** | **23239** | **83.7%** | **93.0** | **103.7** | **166.4** |
| c12000 | Nearest | 260.6 | 317.6 | 51.50 | 5.5 | 24298 | 74.0% | 97.2 | 215.9 | 245.3 |
| | Fastest | 269.6 | 329.1 | 53.50 | 6.0 | 25229 | 76.9% | 100.9 | 149.1 | 238.6 |
| | MinChgOpT | 277.8 | 339.3 | 55.10 | 6.3 | 25988 | 79.5% | 104.0 | 79.8 | 223.3 |
| | DynaThreshold | 291.1 | 357.6 | 59.10 | 7.4 | 27887 | 84.6% | 111.5 | 154.8 | 329.6 |
| | OptChg | 307.9 | 378.4 | 62.8 | 7.6 | 29642 | 91.0% | 118.6 | 155.9 | 217.0 |
| | **CongestionAware** | **319.5** | **393.0** | **65.3** | **8.2** | **30808** | **95.9%** | **123.2** | **99.6** | **195.4** |
| c14000 | Nearest | 322.1 | 392.3 | 63.40 | 6.7 | 29916 | 78.9% | 119.7 | 233.9 | 293.9 |
| | Fastest | 332.4 | 405.4 | 65.60 | 7.4 | 30942 | 82.0% | 123.8 | 170.9 | 267.9 |
| | MinChgOpT | 345.1 | 421.1 | 68.10 | 7.9 | 32131 | 85.8% | 128.5 | 72.5 | 246.4 |
| | DynaThreshold | 359.1 | 440.3 | 72.20 | 9.0 | 34047 | 90.9% | 136.2 | 145.0 | 399.7 |
| | OptChg | 374.3 | 457.9 | 74.9 | 8.7 | 35333 | 95.7% | 141.3 | 165.5 | 196.8 |
| | **CongestionAware** | **373.7** | **455.9** | **73.9** | **8.3** | **34862** | **95.3%** | **139.4** | **76.9** | **169.4** |

## Appendix F. The day-ahead charging planning model for the OptChg policy

We adapt the day-ahead charging plan model of Ma (2021) to our problem. In that model, charging station capacity, heterogeneous charging stations, and minimum charging duration requirement are not considered for the day-ahead charging planning. Furthermore, the charging waiting time $W_{vh}$ and energy consumption $\delta_{vh}$ are vehicle-specific based on each vehicle's historical driving experiences for the OptChg policy, instead of using charger-type-specific waiting time $W_{hs}$ and average energy consumption $\delta_h$ in Model P1 of the CongestionAware policy. The decision variable $x_h^v$ determines whether vehicle $v$ charges on epoch $h$ or not, and $y_h^v$ determines the amount of charged energy of vehicle v on h. The objective function minimizes total charging costs for a 24-hour planning horizon. Eq. (F2) tracks a vehicle's SoC changes on h. Eq. (F3) sets the boundary of vehicles' battery levels. Eq. (F4) sets the initial battery levels of vehicles. Eq. (F5) binds $x_h^v$ and $y_h^v$. Eq. (F6) states that vehicles can be recharged only when its battery level is below 80% ($E_v^{max}$) of its battery capacity $B$. Eq. (F7) states that the maximum amount of energy can be charged during one epoch cannot exceed the maximum energy delivery limit of a fast charger. Note that the notation remains the same as indicated in the notation table, except $\bar{\varphi}$ here denotes the average charging power (i.e., average charging power of fast and slow chargers), since the day-ahead model of our previous study does not distinguish charger types.). Furthermore, $\widehat{H}^\ell$ here is the entire planning horizon (not the shifted planning horizon as used in Model P1) as $W_{vh}$ and $\delta_{vh}$ is vehicle-specific. We set $M_1 = Y_{fast}^{max}$ and $M_2 = B$ to tighten the MILP formulation.

$$\min Z_1 = \sum_{v\in V}\sum_{h\in\widehat{H}^\ell}\left((p_h+\frac{\gamma}{\bar{\varphi}})y_h^v + (C+\gamma W_{vh})x_h^v\right) + \sum_{v\in V}\bar{p}(B_v - e_{v,n_{H^\ell}+1}) \quad \text{(F1)}$$

Subject to

$$e_{v,h+1} = e_{vh} - \delta_{vh} + y_h^v, \forall v \in V, h \in \widehat{H}^\ell \quad \text{(F2)}$$

$$E_v^{min} \le e_{vh} \le B, \forall v \in V, h \in H^\ell \cup \{n_{\widehat{H}^\ell} + 1\} \quad \text{(F3)}$$

$$e_{v1} = E_v^0, \forall v \in V \quad \text{(F4)}$$

$$y_h^v \le M_1 x_h^v, \forall v \in V, h \in \widehat{H}^\ell \quad \text{(F5)}$$

$$e_{vh} \le E_v^{max} + M_2(1 - x_h^v), \forall v \in V, h \in \widehat{H}^\ell \quad \text{(F6)}$$

$$0 \le y_h^v \le Y_{fast}^{max}, \forall v \in V_h, h \in \widehat{H}^\ell \quad \text{(F7)}$$

$$x_h^v \in \{0,1\}, \forall v \in V, h \in \widehat{H}^\ell \quad \text{(F8)}$$

We compare the day-ahead charging plans obtained from the OptChg policy with that of the CongestionAware policy for c3000 and c4000 scenarios. As a reminder, the charging infrastructure consists of 12 chargers, including 6 fast chargers and 6 slow chargers. For c3000 scenario (see Figure F1), the number of scheduled charging vehicles increases from the 7th hour to its peak at around the 9th hour (40 vehicles scheduled to charge on one charging decision epoch). For the c4000 scenario, the number of scheduled charging vehicles has multiple peaks between the 5th hour and the 8th hour (with its highest demand of 90 vehicles scheduled to charge), then drops to one third (i.e., 30 vehicles to charge). This plan is in response to vehicles' charging waiting times $W_{vh}$ (see Figure F2) to minimize its objective function value (vehicles are scheduled to charge when charging waiting time is zero or low). As the charging stations are uncapacitated, this results in high charging demand compared to the capacitated case obtained from the CongestionAware policy. Figure F3 presents the distribution of $W_{hs}$ on fast chargers and slow chargers used for the CongestionAware policy. Note that both day-ahead charging plans are derived with the same available input data and information, but the applied MILP models differ.

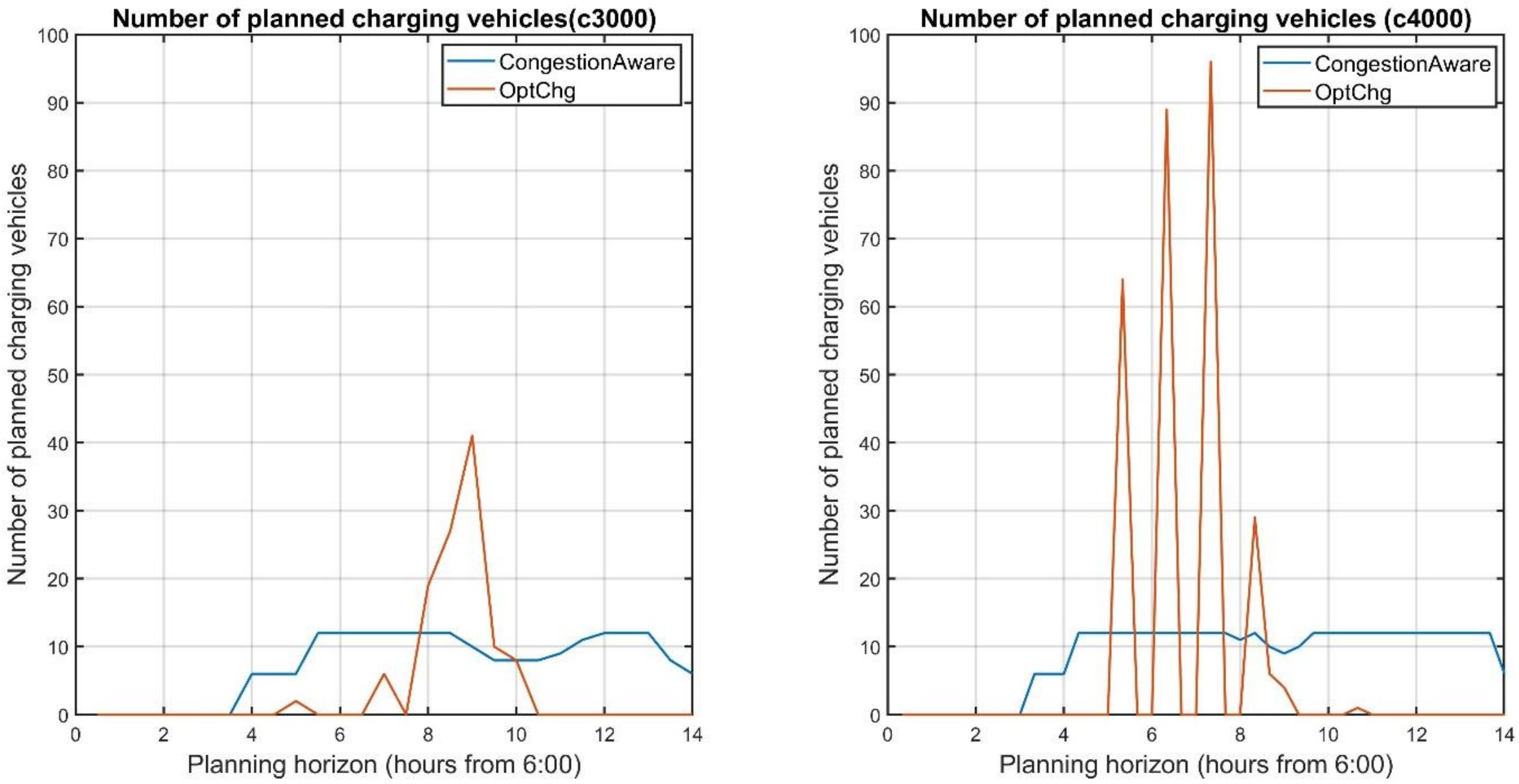

Figure F1. Comparison of the day-ahead charging plans obtained from the OptChg and from the CongestionAware policies for c3000 scenario (left) and c4000 scenario (right).

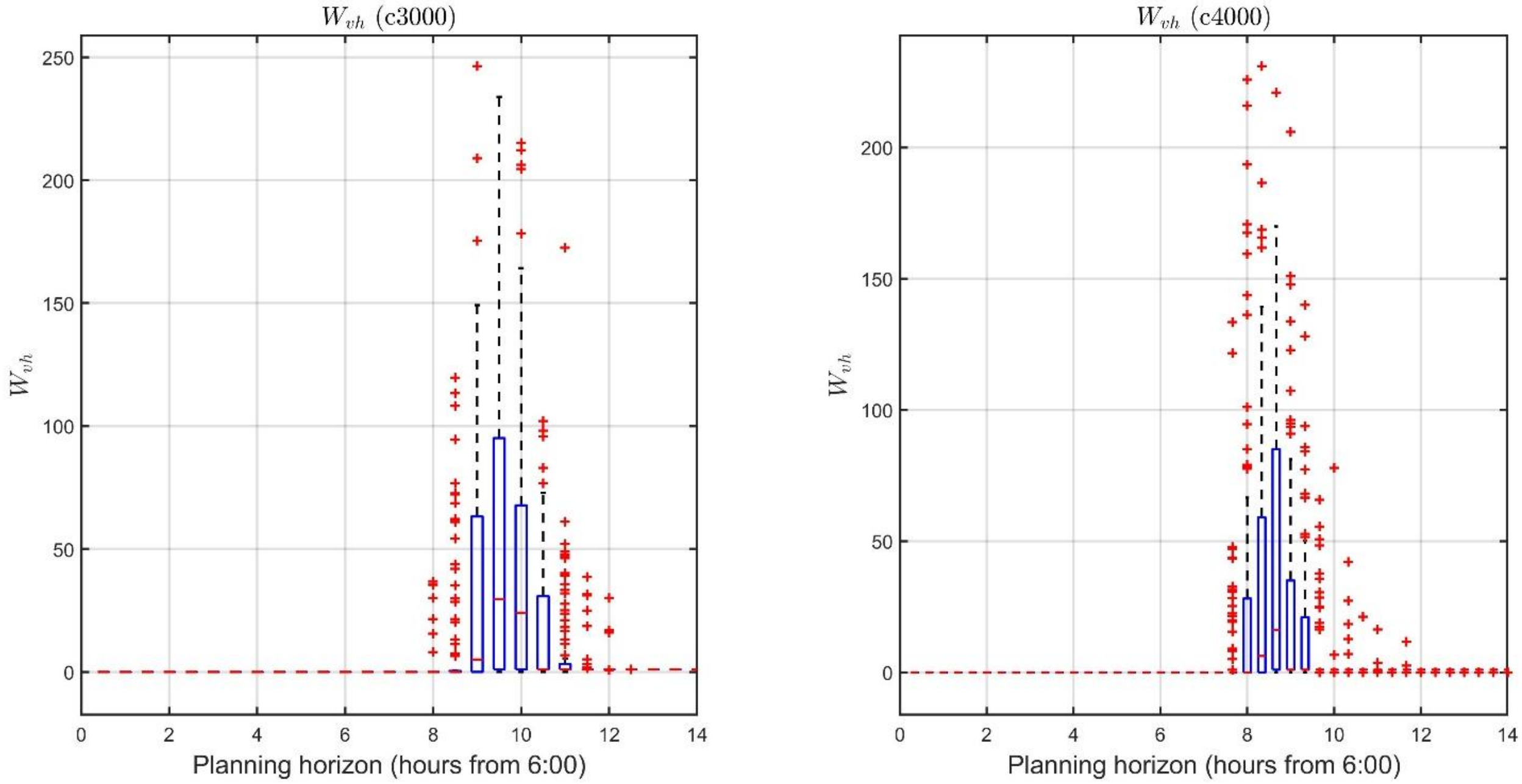

Figure F2. Histogram of $W_{vh}$ used for the OptChg policy for c3000 (left) and c4000 scenarios (right).

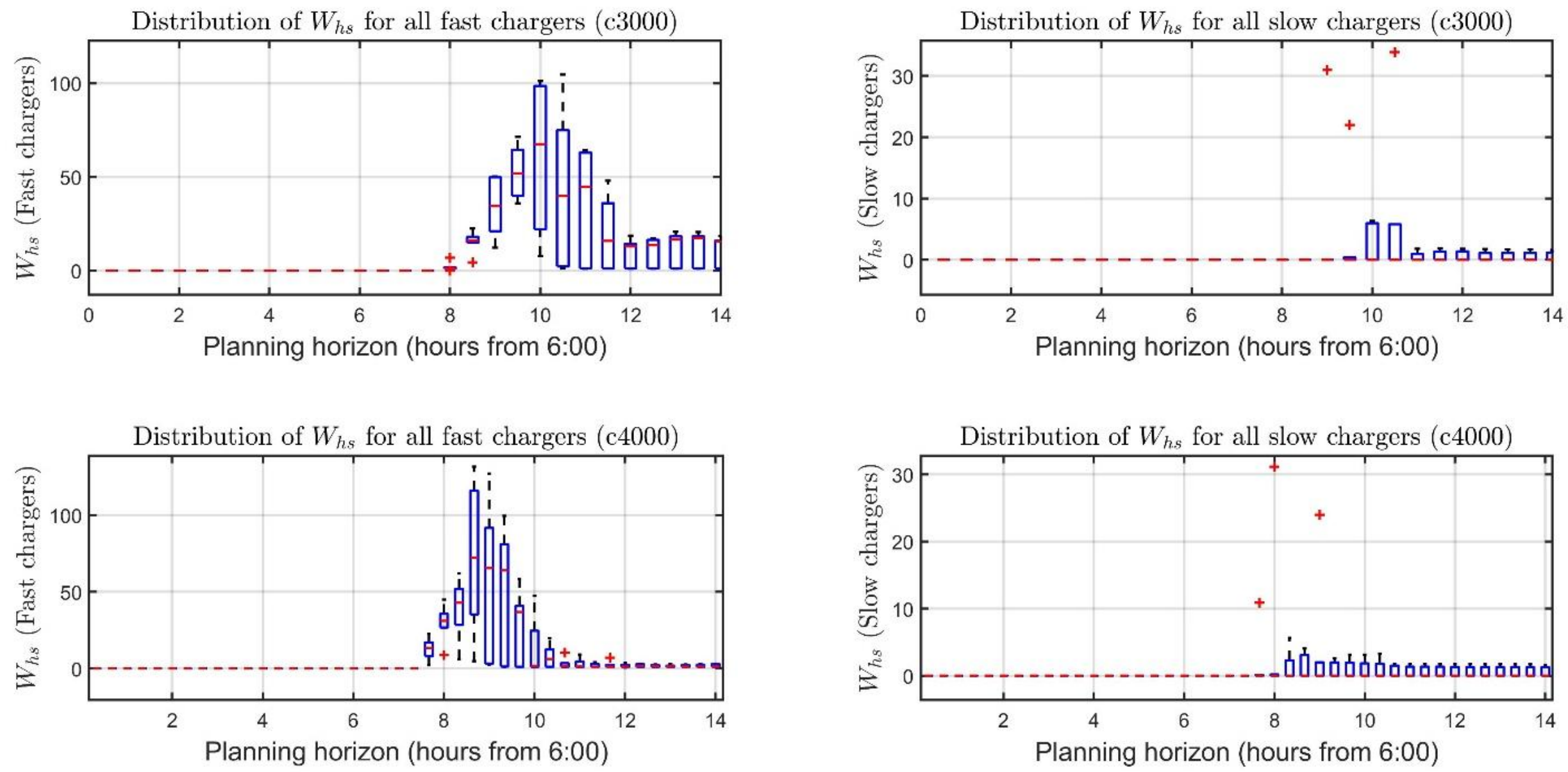

Figure F3. Histogram of $W_{hs}$ on fast chargers (left) and slow chargers (right) used for the CongestionAware policy for c3000 (upper part) and c4000 scenarios (lower part).

**Appendix G. Other detailed results for the c3000 and c4000 scenarios.**

a. Standard deviations of the KPIs for the c3000 and c4000 scenarios

The randomly generated customer demands give rise to variability in the performance of each charging policy. Standard deviations of the KPIs are shown in Table G1, calculated over the 5 validation datasets. In general, the CongestionAware policy and the OptChg policy exhibit lower variability than rule-based policies , demonstrating that it is more stable in the face of demand variability. This is particularly notable for those KPIs where CongestionAware has higher absolute values: profit, revenue and service rate. Additionally, the optimal charging focus of CongestionAware is highlighted by the low variability (by comparison) in total charged energy and charger waiting time. This implies drivers can expect lower and more stable queuing times with the presence of demand variability. It improves drivers' charging experience/stress and the operational performance. For the fleet operator, this predictability is important for driver scheduling and revenue forecasting, reducing the operational/financial risk associated with demand surges.

Table G1. Standard deviations of the KPIs for different charging policies.

| Scenario | Charging policy | PF | TR | TC | CC | ENG | SR | KMT | WT | CT |
|---|---|---|---|---|---|---|---|---|---|---|
| c3000 (3000 requests) | Nearest | 0.52 | 0.66 | 0.11 | 0.03 | 52.7 | 0.71% | 0.2 | 11.2 | 1.4 |
| | Fastest | 0.51 | 0.62 | 0.10 | 0.03 | 47.8 | 0.59% | 0.2 | 7.5 | 2.1 |
| | MinChgOpT | 0.52 | 0.65 | 0.11 | 0.03 | 51.1 | 0.72% | 0.2 | 4.3 | 1.7 |
| | DynaThreshold | 0.24 | 0.32 | 0.13 | 0.03 | 60.9 | 0.36% | 0.2 | 16.1 | 3.5 |
| | OptChg | 0.38 | 0.46 | 0.08 | 0.02 | 36.5 | 0.56% | 0.1 | 1.3 | 1.8 |
| | **CongestionAware** | **0.33** | **0.41** | **0.09** | **0.02** | **40.5** | **0.41%** | **0.2** | **4.5** | **3.0** |
| c4000 (4000 requests) | Nearest | 0.60 | 0.75 | 0.11 | 0.06 | 52.9 | 0.51% | 0.2 | 13.2 | 3.9 |
| | Fastest | 0.71 | 0.85 | 0.11 | 0.05 | 52.6 | 0.65% | 0.2 | 18.4 | 2.4 |
| | MinChgOpT | 0.66 | 0.79 | 0.11 | 0.05 | 53.4 | 0.57% | 0.2 | 16.8 | 1.9 |
| | DynaThreshold | 0.31 | 0.26 | 0.06 | 0.05 | 27.3 | 0.33% | 0.1 | 17.0 | 3.4 |

| | | | | | | | | | |
|---|---|---|---|---|---|---|---|---|---|
| OptChg | 0.30 | 0.35 | 0.05 | 0.01 | 24.4 | 0.34% | 0.1 | 1.1 | 1.8 |
| **CongestionAware** | **0.31** | **0.37** | **0.06** | **0.01** | **27.0** | **0.31%** | **0.1** | **3.1** | **2.0** |

b. The histograms of the SoC of vehicles at the end of service hours for c3000 and c4000

The histograms of the SoC of vehicles at the end of service hours for c3000 and c4000 are shown in Figures G1 and G2. For the CongestionAware and OptChg policies, the SoC of vehicles at the end of service hours is strongly concentrated between $E_v^{min}$ (6.2 kWh) and 10 kWh for both cases (80% or more of the fleet). For the rule-based benchmark policies, the distributions of the SoC of vehicles at the end of service hours for both cases show higher concentration (40% for c3000 and 50% for c4000) within the range of around $E_v^{min}$ and 10 kWh.

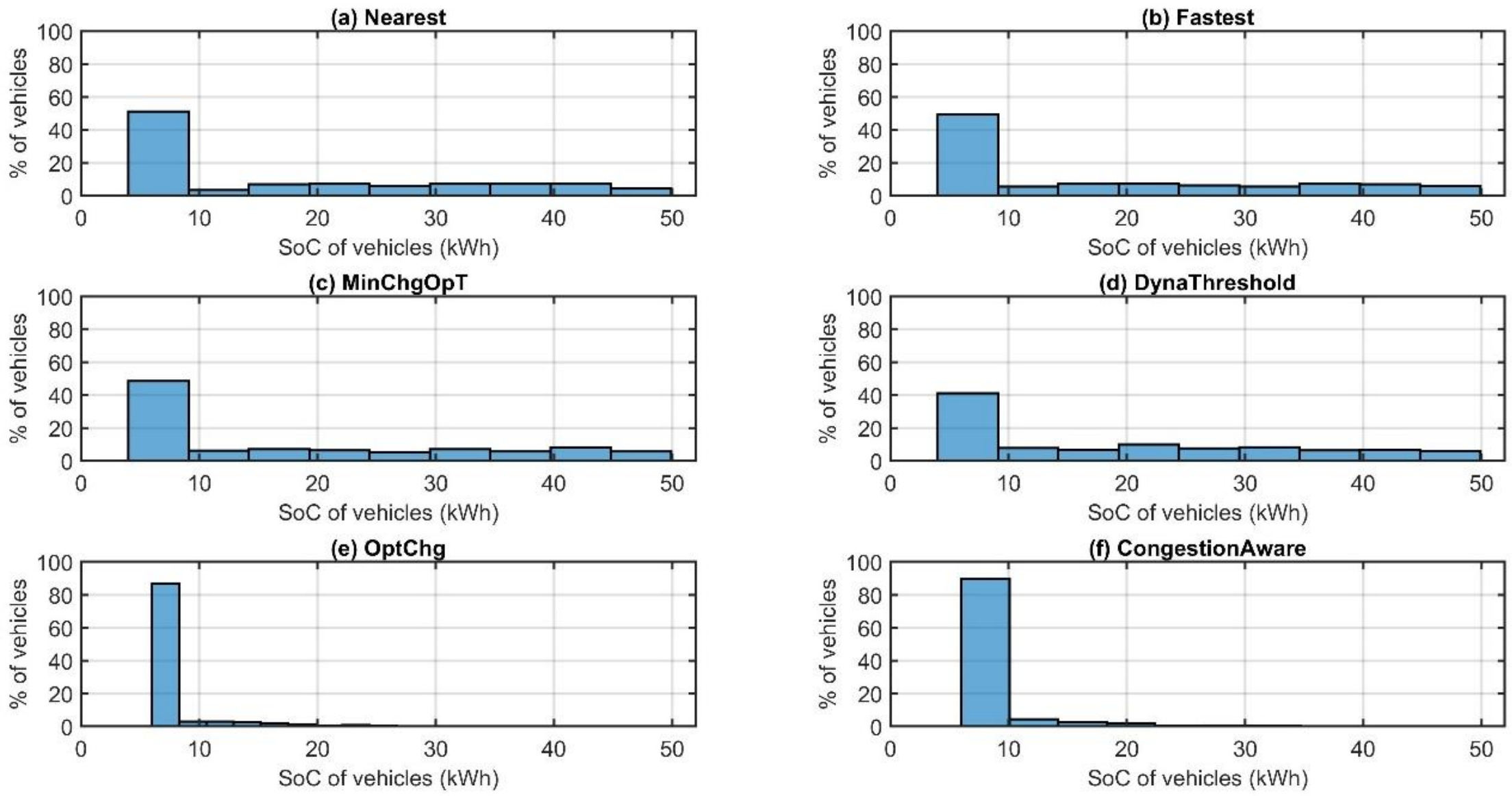

Figure G1. Histograms of the SoCs of vehicles at the end of service hours (20:00) for different charging policies (c3000 scenario).

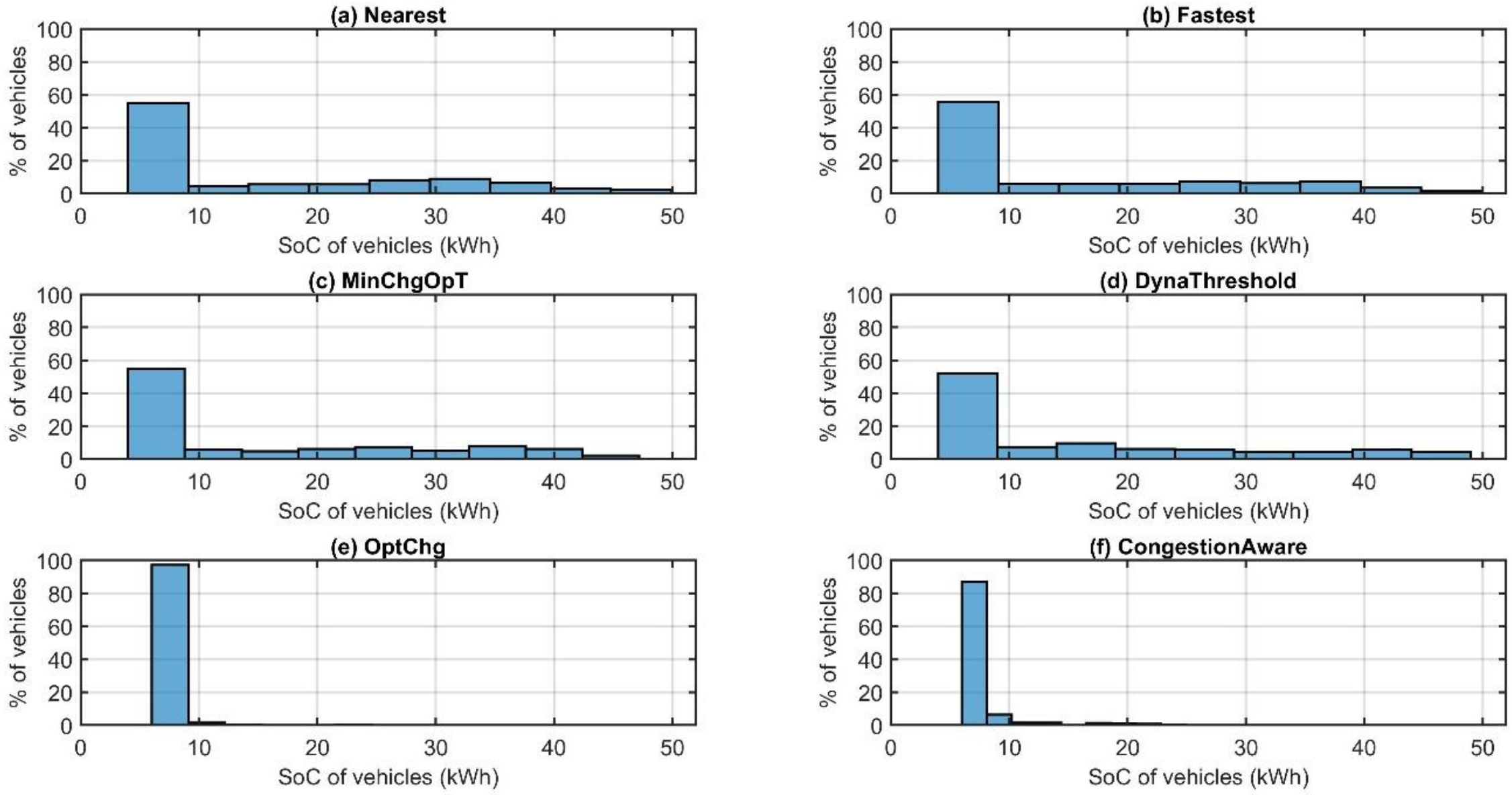

Figure G2. Histograms of the SoCs of vehicles at the end of service hours (20:00) for different charging policies (c4000 scenario).

c. A representative vehicle's energy profile and its state changes using different charging policies during the service hours.

Figure G3 illustrates a representative vehicle's energy profile and its state changes using different charging policies during the service hours. The two left subfigures in Figure G3 show the vehicle energy profiles, while the two right subfigures present the vehicle state changes for CongestionAware and DynaThres policies. We can observe that the vehicle charges to 80% of its battery capacity when using a rule-based charging policy, while using the CongestionAware and OptChg policies the vehicle charges less amount energy either based on its charging plan or by anticipating its energy need to the end of the day. The CongestionAware policy schedules the vehicle to charge much earlier compared to the OptChg because its day-ahead charging plan considers more realistic charging operation constraints (see Appendix F for a more detailed explanation). When using the Nearest policy, the vehicle charges at a slow charger to 80% instead of using fast chargers as the other charging policies. The data tips show the moment of the vehicle's state change. Regarding the CongestionAware policy, we can observe that for the second charging operation the vehicle arrives at a fast charger at the 615th minute, waits for 29.3 minutes, and charges 19.0 minutes with 15.83 kWh energy charged.

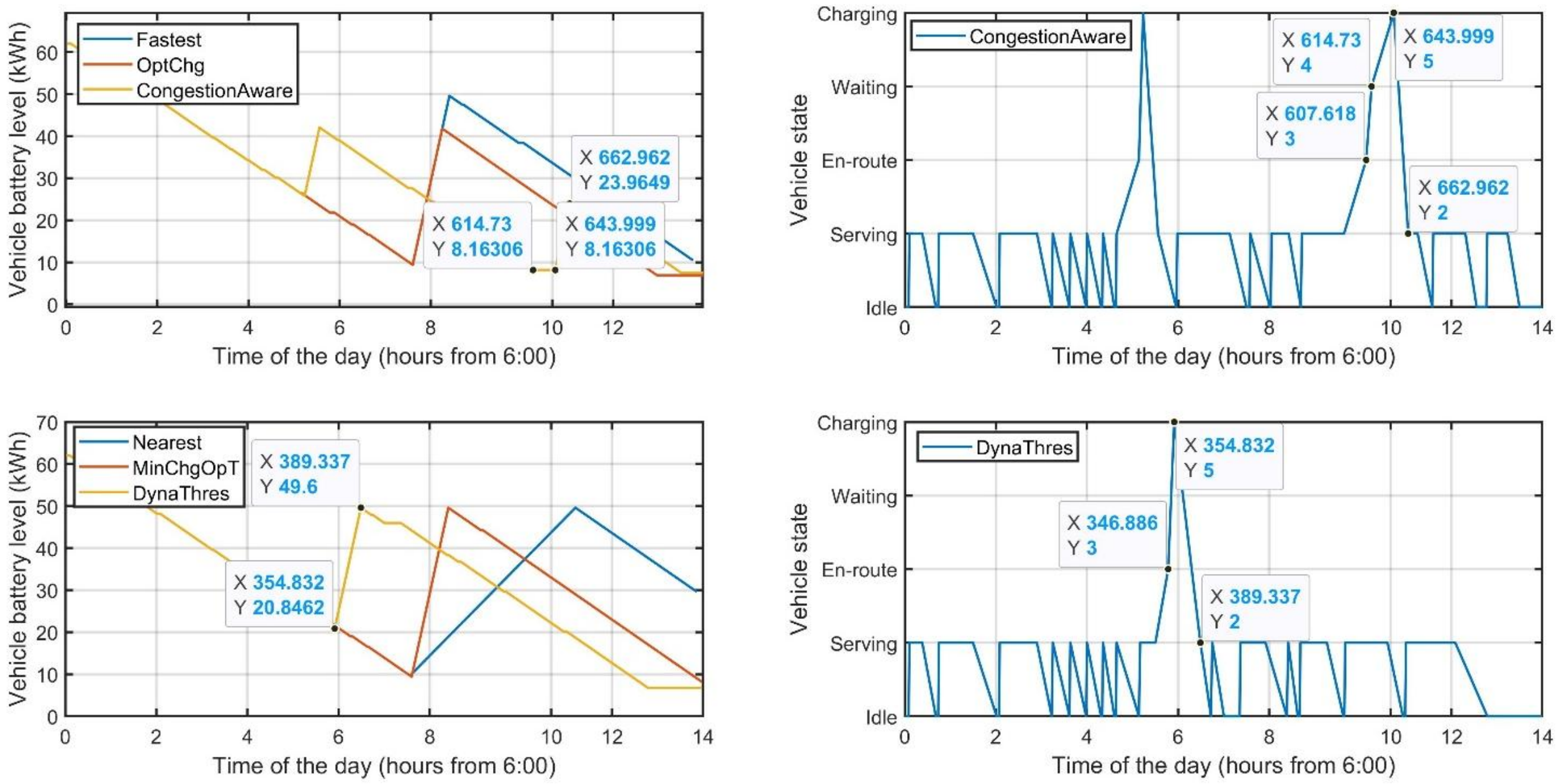

Figure G3. Vehicle state-of-charge (SoC) profile over the service period under different charging policies (left), and the corresponding vehicle state trajectory under the CongestionAware policy (upper right) and the DynaThres policy (lower right). The x-coordinate values displayed in the data tips are expressed in minutes.

## Appendix H. The impact of using re-estimated $\delta_h$ on the performance of the CongestionAware policy

The day-ahead charging plan P1 relies on the average energy-consumption rate $\delta_h$ to anticipate vehicles' time-dependent energy needs. In the main experiments these $\delta_h$ values are estimated from simulations under the Fastest policy. In our day-ahead charging plan model P1, we use this average energy consumption rate $\delta_h$ to capture time-dependent energy consumption of vehicles in response to customer demand fluctuation. We conduct a one-iteration demonstration of consistency by re-estimating $\delta_h$ based on the realized charging operations under the CongestionAware policy for the five test instances. Then we re-run the simulation (re-run the models P1-P3) with the same test instances for a 24-hour period and then compare the result with the case without re-estimating it. As a reminder, we generate 15 test instances with 10 for model parameter estimation, and the remaining 5 instances used to test the proposed methodology. Figure H1 shows the profile of the average vehicle energy consumption $\delta_h$ and that of the re-estimated one. We can observe that $\delta_h$ is higher than the re-estimated one near the end of the service hours for the c3000 scenario. For the c4000 scenario, there is a significant drop from the 10th hour to the end of service hours. Figure H2 compares the day-ahead charging plans of 10 representative vehicles using $\delta_h$ and re-estimated $\delta_h$, respectively, for c3000 (subfigures (a) and (b)) and c4000 scenarios (subfigures (c) and (d)). We can observe that in response to the respective estimated vehicle energy consumption $\delta_h$, the timing and frequency of scheduled charging operations change. Some vehicles have more scheduled charging operations if the re-estimated $\delta_h$ is used, in particular for the c4000 case. However, the performance of the model remains stable using $\delta_h$ or re-estimated $\delta_h$ for both demand scenarios (see Table H1).

The day-ahead plan being built on imperfect estimates of $\delta_h$ has limited sensitivity to the one-iteration based post re-estimation. This is because the online anticipation layer re-optimises vehicle's target SoC, whenever its current SoC falls below the charging threshold, against the realised system state at every epoch, absorbing any gap between planned and realised energy consumption. Note that the accuracy of anticipated charging waiting times $\overline{W}_{hs}$ may affect the day-ahead charging plan. Underestimating charging waiting times during charging peak hours may lead to excessive charging congestion. In our online vehicle-to-charger assignment model P3, we use the current charger occupancy and reservation information to estimate charging waiting times and optimize charger assignments, which allows us to partially mitigate these estimation errors. The impact of waiting-time estimation errors on the system performance remains to be investigated.

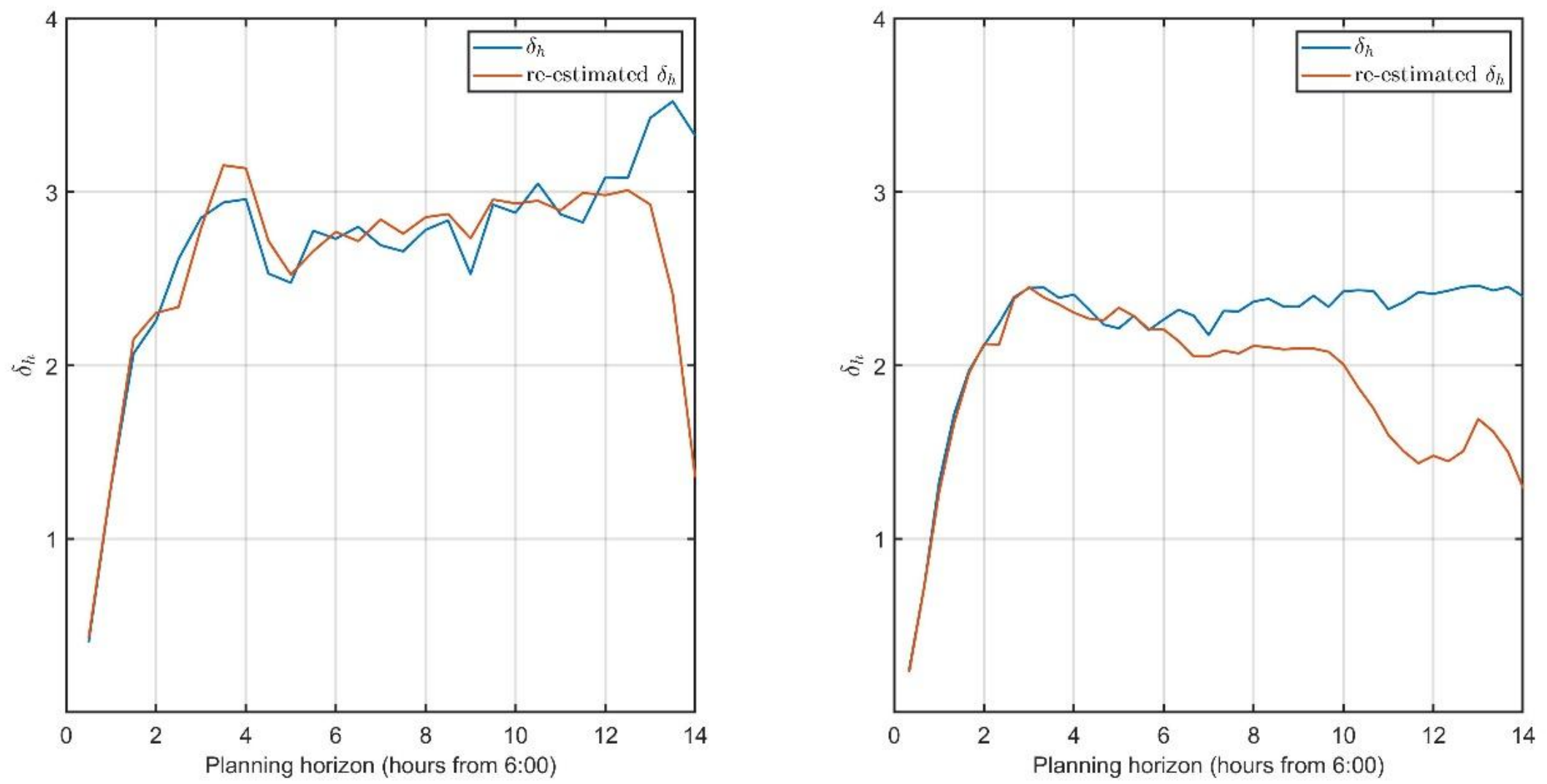

Figure H1. $\delta_h$ and post re-estimated $\delta_h$ for c3000 and c4000 scenarios.

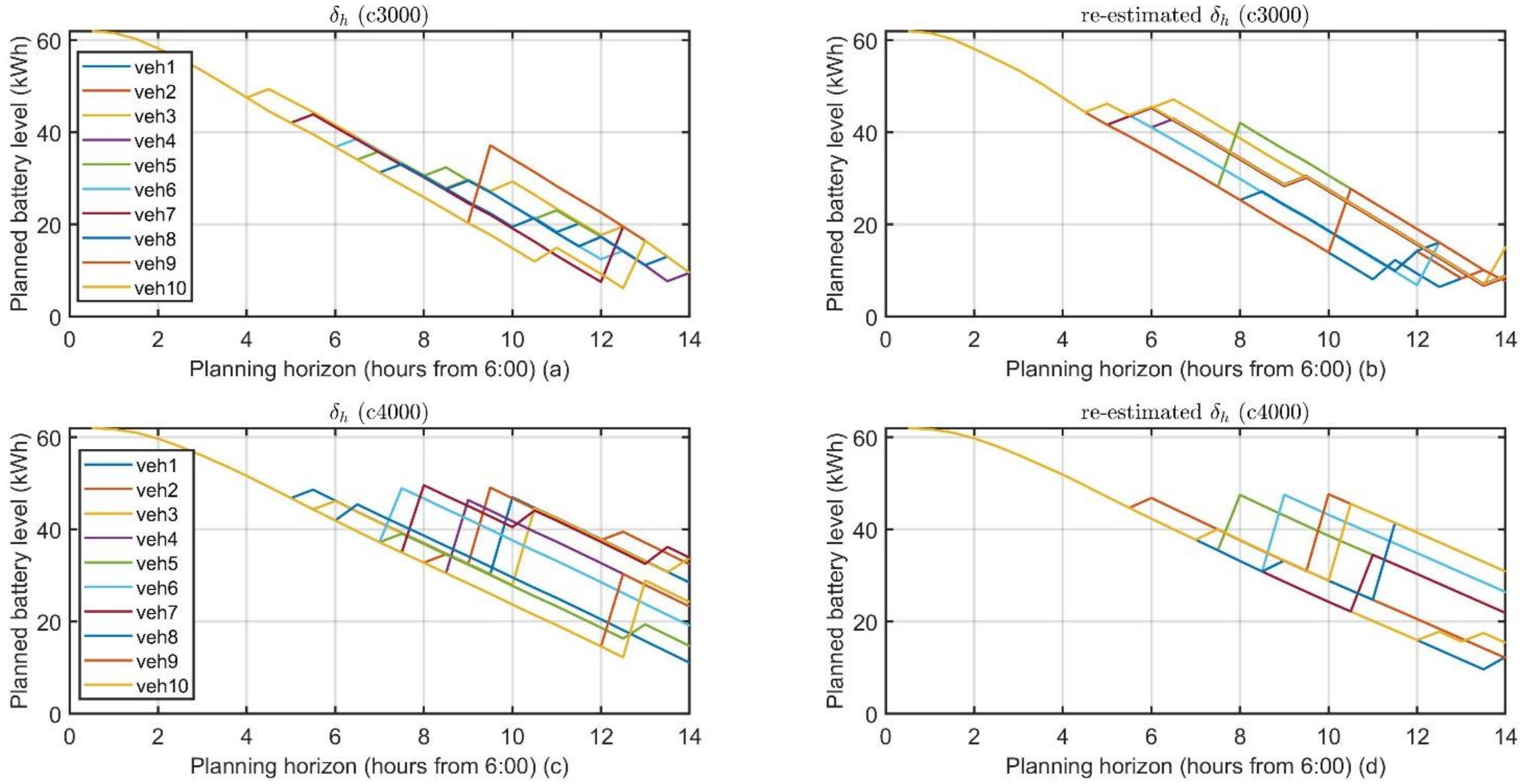

Figure H2. The comparison of vehicle charging plans (i.e., planned vehicle battery levels) during service hours using $\delta_h$ (on the left column) and re-estimated $\delta_h$ (on the right column) for Vehicle 1 to Vehicle 10 for c3000 and c4000 scenarios.

Table H1. Comparison of the KPIs for CongestionAware policy with and without re-estimated $\delta_h$.

| Scenario | Re-estimated $\delta_h$ | PF | TR | TC | CC | ENG | SR | KMT | WT | CT |
|---|---|---|---|---|---|---|---|---|---|---|
| c3000 | No | **72.70** | **89.88** | **15.41** | **1.77** | **7267** | **86.53%** | **29.1** | **26.8** | **41.2** |
| | Yes | 72.67 | 89.82 | 15.39 | 1.77 | 7257 | 86.43% | 29.0 | 28.0 | 37.8 |
| c4000 | No | **79.02** | **98.06** | **16.93** | **2.11** | **7984** | **67.13%** | **31.9** | **43.0** | **62.8** |
| | Yes | 78.71 | 97.68 | 16.88 | 2.09 | 7961 | 67.10% | 31.8 | 45.9 | 60.6 |

**Appendix I. Impact of decision interval length of P1 and online charging-first/service-first decision-order sensitivity analysis for the performance of the CongestionAware policy.**

In this appendix, we first analyze the impact of decision interval lengths on the total profit for the 5 test instances. The tested values for the interval length are 10, 20, 30, 40, and 50 minutes for the c3000 and c4000 scenarios. The results are reported in Table I1. For the c3000 scenario, a shorter time interval from 30 minutes to 20 minutes improves the system performance (+2.10% profit), while further reducing the time interval to 10 minutes provides no improvement. However, using a longer interval from 30 minutes to 40 minutes results in lower system performance (-2.81%). For the c4000 scenario, slight improvement (+1.78% profit) can be observed for the system performance when reducing the interval from 20 minutes to 10 minutes. The results show that in general, shorter intervals (10 or 20 minutes) increase total charged energy with a higher service rate. The results imply the temporal resolution affects the operational performance. In practice, this parameter can be tuned to achieve better operational performance.

Table I1. Impact of the decision interval length of P1 on the CongestionAware policy.

| Scenario | Interval length of P1 (minute) | PF | TR | TC | CC | ENG | SR | KMT | WT | CT |
|---|---|---|---|---|---|---|---|---|---|---|
| c3000[1] | 10 | 73.76 | 91.54 | 15.82 | 1.96 | 7461 | 87.89% | 29.8 | 41.6 | 58.8 |
| | 20 | 74.23 | 91.96 | 15.84 | 1.89 | 7474 | 88.54% | 29.9 | 34.5 | 48.7 |
| | 30 | 72.70 | 89.88 | 15.41 | 1.77 | 7267 | 86.53% | 29.1 | 26.8 | 41.2 |
| | 40 | 70.66 | 87.40 | 14.86 | 1.88 | 7009 | 83.36% | 28.0 | 24.6 | 48.0 |
| c4000[1] | 10 | 80.43 | 100.17 | 17.46 | 2.27 | 8236 | 68.44% | 32.9 | 58.0 | 78.7 |
| | 20 | 79.02 | 98.06 | 16.93 | 2.11 | 7984 | 67.13% | 31.9 | 43.0 | 62.8 |

Remark: 1. For the c3000 (c4000) scenario, no feasible solution can be obtained when the decision interval length is set to 50 (30) minutes or longer.

We compare the performance of the alternative by inverting the order of the charging decision (P3) and service decision (P2). For this **service-first implementation**, the new steps 1-6 (Table 1) remain unchanged. The new step 7 updates the lists of unserved customers $R_{t(i)}$ and the list of idle vehicles $V_{t(i)}$. The new step 8 solves the P2 based on the obtained solutions to assign vehicles to serve customers. The new step 9 then removes the assigned vehicles from $V_{t(i)}$ and updates their states. The new steps 9-16 remain the same as the original Steps 7-14 in Table 1 for vehicles' charging operations. The new steps 17-19 remain unchanged to advance the system to t(i+1), which compute the overnight charging costs, and output the profit.

The computational experiments are conducted with respect to the c3000 and c4000 scenarios with three vehicle battery capacity settings of 62, 72 and 82 kWh, respectively. The results in Table I2 show that service-first-and-charging-second online decision order results in systematically lower profit and service rate compared with the charging-first-and-service-second baseline for high demand scenarios (c4000). This is because a charging-first decision order facilitates the implementation of the day-ahead charging plan, helping to mitigate charging congestion while ensuring that low-SoC vehicles are recharged in a timely manner to maintain better fleet availability. However, service-first can perform better under moderate demand (c3000), particularly with larger battery capacity.

Table I2. Impact of the online charging-first/service-first decision order on the performance of the CongestionAware policy.

| Scenario | Battery (kWh) | Charging-first | PF | TR | TC | CC | ENG | SR | KMT | WT | CT |
|---|---|---|---|---|---|---|---|---|---|---|---|
| c3000 | 62 | **Yes** | **72.70** | **89.88** | **15.41** | **1.77** | **7267** | **86.5%** | **29.1** | **26.8** | **41.2** |
| | | no | 72.46 | 89.47 | 15.27 | 1.73 | 7204 | 87.1% | 28.8 | 23.8 | 36.5 |
| | 72 | **Yes** | **75.34** | **92.66** | **15.71** | **1.61** | **7411** | **89.6%** | **29.6** | **7.3** | **25.8** |
| | | no | 75.68 | 93.07 | 15.78 | 1.61 | 7443 | 90.7% | 29.8 | 9.6 | 25.7 |
| | 82 | **Yes** | **78.08** | **95.79** | **16.17** | **1.54** | **7627** | **93.0%** | **30.5** | **1.0** | **19.3** |
| | | no | 78.45 | 96.19 | 16.19 | 1.55 | 7638 | 94.1% | 30.6 | 2.0 | 19.7 |
| c4000 | 62 | **Yes** | **79.02** | **98.06** | **16.93** | **2.11** | **7984** | **67.1%** | **31.9** | **43.0** | **62.8** |
| | | no | 76.59 | 94.82 | 16.27 | 1.96 | 7677 | 65.2% | 30.7 | 32.7 | 48.7 |
| | 72 | **Yes** | **84.14** | **104.33** | **18.05** | **2.13** | **8515** | **72.3%** | **34.1** | **34.2** | **54.7** |
| | | no | 83.30 | 102.97 | 17.68 | 2.00 | 8339 | 71.6% | 33.4 | 27.8 | 44.9 |
| | 82 | **Yes** | **88.10** | **108.89** | **18.76** | **2.04** | **8847** | **76.0%** | **35.4** | **22.7** | **37.5** |
| | | no | 87.39 | 107.82 | 18.50 | 1.93 | 8727 | 75.6% | 34.9 | 18.3 | 36.2 |

**Appendix J. Paired statistical evaluation of the performance of CongestionAware and the benchmark policies across 30 independent test instances.**

We further evaluate the performance comparison from five to 30 independent test instances, with the day-ahead charging plan obtained from 10 separate instances for the c3000 and c4000 scenarios with battery sizes of 62 kWh and 72 kWh, respectively. For each benchmark policy, we calculate the instance-level difference as CongestionAware minus the benchmark. Table J1 reports mean paired differences, 95% confidence intervals, and standardized paired effect sizes. The two-sided paired t-tests are used to test the statistical significance. The effect size statistics $d_z$, measuring the pairwise gains relative to their between-instance variability, are defined as follows.

$$d_z = \frac{\text{mean paired improvement}}{\text{standard deviation of the paired improvements}}$$

As shown in Table J1, the CongestionAware policy achieves higher profit and service rate on all 30 instances for each scenario. All reported paired t-tests have p-values < 0.001. The corresponding profit gains range from 1.89% (OptChg for the c3000 scenario with a 72 kWh battery) to 21.07% (Nearest for the c4000 scenario with a 62 kWh battery). The service-rate gains range from 2.43% to 21.37%, for the same two comparisons, respectively. The results provide reliable performance gains of the CongestionAware policy.

Table J1. Performance differences and paired comparisons between CongestionAware and the five benchmark policies across 30 independent test instances.

| **Scenario** | **Battery (kWh)** | **Benchmark** | **% of Profit difference [95% CI]** | $d_z$ | **% of Service-rate difference [95% CI]** | $d_z$ |
|---|---|---|---|---|---|---|
| c3000 | 62 | Nearest | 15.50 [15.18, 15.81] | 18.21 | 18.74 [18.41, 19.06] | 21.59 |
| | | Fastest | 14.27 [13.90, 14.65] | 14.28 | 17.25 [16.85, 17.65] | 16.12 |
| | | MinChgOpT | 12.89 [12.53, 13.24] | 13.38 | 15.83 [15.46, 16.21] | 15.66 |
| | | DynaThreshold | 8.44 [8.10, 8.79] | 9.17 | 10.70 [10.31, 11.09] | 10.18 |
| | | OptChg | 3.66 [3.36, 3.95] | 4.62 | 4.59 [4.26, 4.91] | 5.24 |
| | 72 | Nearest | 11.92 [11.54, 12.29] | 11.85 | 13.59 [13.18, 14.01] | 12.27 |
| | | Fastest | 11.04 [10.65, 11.43] | 10.53 | 12.53 [12.09, 12.96] | 10.80 |
| | | MinChgOpT | 10.29 [9.92, 10.65] | 10.51 | 11.79 [11.40, 12.18] | 11.28 |
| | | DynaThreshold | 8.90 [8.49, 9.31] | 8.04 | 10.16 [9.68, 10.64] | 7.89 |
| | | OptChg | 1.89 [1.65, 2.13] | 2.97 | 2.43 [2.12, 2.73] | 2.99 |
| c4000 | 62 | Nearest | 21.07 [20.69, 21.44] | 21.11 | 21.37 [21.02, 21.71] | 23.27 |
| | | Fastest | 19.35 [18.97, 19.73] | 19.21 | 19.76 [19.44, 20.08] | 23.13 |
| | | MinChgOpT | 18.23 [17.86, 18.61] | 18.20 | 18.77 [18.42, 19.12] | 20.07 |
| | | DynaThreshold | 12.87 [12.48, 13.26] | 12.41 | 14.84 [14.46, 15.22] | 14.64 |
| | | OptChg | 3.25 [3.08, 3.42] | 7.15 | 3.88 [3.69, 4.06] | 7.78 |
| | 72 | Nearest | 20.46 [20.19, 20.72] | 29.07 | 19.72 [19.41, 20.04] | 23.49 |
| | | Fastest | 19.43 [19.12, 19.74] | 23.65 | 18.84 [18.49, 19.19] | 20.02 |
| | | MinChgOpT | 18.72 [18.41, 19.03] | 22.71 | 18.30 [17.95, 18.65] | 19.32 |
| | | DynaThreshold | 15.28 [14.77, 15.79] | 11.24 | 16.55 [16.03, 17.08] | 11.78 |
| | | OptChg | 7.13 [6.85, 7.40] | 9.77 | 7.75 [7.44, 8.06] | 9.30 |

Remark: The Holm-adjusted two-sided paired t-test p-values are <0.001 for all paired comparisons.